# ON FRAGMENTS WITHOUT IMPLICATIONS OF BOTH THE FULL LAMBEK LOGIC AND SOME OF ITS SUBSTRUCTURAL EXTENSIONS

ÀNGEL GARCÍA-CERDAÑA[1,2] AND VENTURA VERDÚ[3]

ABSTRACT. In this paper we study some fragments without implications of the (Hilbert) full Lambek logic **HFL** and also some fragments without implications of some of the substructural extensions of that logic. To do this, we perform an algebraic analysis of the Gentzen systems defined by the substructural calculi $\mathbf{FL}_\sigma$. Such systems are extensions of the full Lambek calculus **FL** with the rules codified by a subsequence, $\sigma$, of the sequence $ew_lw_lc$; where $e$ stands for *exchange*, $w_l$ for *left weakening*, $w_r$ for *right weakening*, and $c$ for *contraction*. We prove that these Gentzen systems (in languages without implications) are algebraizable by obtaining their equivalent algebraic semantics. All these classes of algebras are varieties of pointed semilatticed monoids and they can be embedded in their ideal completions. As a consequence of these results, we reveal that the fragments of the Gentzen systems associated with the calculi $\mathbf{FL}_\sigma$ are the restrictions of them to the sublanguages considered, and we also reveal that in these languages, the fragments of the external systems associated with $\mathbf{FL}_\sigma$ are the external systems associated with the restricted Gentzen systems (i.e., those obtained by restriction of $\mathbf{FL}_\sigma$ to the implication-less languages considered). We show that all these external systems without implication have algebraic semantics but they are not algebraizable (and are not even protoalgebraic). Results concerning fragments without implication of intuitionistic logic without contraction were already reported in Bou et al. (2006) and Adillon et al. (2007).

## 1. Introduction and Outline of the Paper.

This paper contributes to the study of fragments without implications of the (Hilbert) full Lambek logic **HFL** and some of its substructural extensions. Hilbert-style axiomatizations for fragments *with* implication of classical and intuitionistic logics are already given in Hosoi (1966a,b). Fragments with implications of some substructural logics have been studied in Ono & Komori (1985), van Alten & Raftery (2004) and Galatos & Ono (2010), where strong separation is proved with respect to the implication connectives for some axiomatizations of $\mathbf{HFL}_e$, $\mathbf{HFL}_{ew_l}$ and **HFL**. As fuzzy logics can also be considered to be substructural logics (see Esteva et al. (2003)), we mention the study of fragments with implication of nine prominent fuzzy logics (MTL, IMTL, SMTL, ΠMTL, BL, SBL, G, Ł and Π) as reported in Cintula et al. (2007).

Concerning the study of fragments *without* implication of intuitionistic logic, we mention the study of fragments in the languages $\langle\vee\rangle$, $\langle\wedge\rangle$, $\langle\vee,\wedge\rangle$ and $\langle\vee,\wedge,\neg\rangle$ (see Porębska & Wroński (1975); Dyrda & Prucnal (1980); Font et al. (1991); Font & Verdú (1991) for fragments in the languages $\langle\vee\rangle$, $\langle\wedge\rangle$, $\langle\vee,\wedge\rangle$ and Rebagliato & Verdú (1993, 1994) for the $\langle\vee,\wedge,\neg\rangle$-fragment). The study of fragments without implication of intuitionistic logics without contraction can be found in: Bou et al. (2006); Adillon et al. (2007). In Bou et al. (2006), the authors analyzed fragments in the languages $\langle\vee,*,\neg,0,1\rangle$ and $\langle\vee,\wedge,*,\neg,0,1\rangle$ of the logic associated with the calculus without contraction, $\mathbf{FL}_{ew}$ . In Adillon et al. (2007), fragments corresponding to the languages $\langle\vee,*,0,1\rangle$, $\langle\vee,\wedge,*,0,1\rangle$ and $\langle\vee,\wedge,0,1\rangle$ were studied.

In this paper we study many fragments without implications of both **HFL** and its substructural extensions, as well as of the Gentzen systems $\mathcal{FL}_\sigma$, and this work can therefore be seen as a continuation of Bou et al. (2006) and Adillon et al. (2007). The systems $\mathcal{FL}_\sigma$ (see Part 1) are defined by the calculi $\mathbf{FL}_\sigma$; that is, the full Lambek calculus extended with the structural rules codified by a subsequence, $\sigma$, of the sequence $ew_lw_lc$, where $e$ stands for *exchange*, $w_l$ for *left weakening*, $w_r$ for *right weakening*, and $c$ for *contraction*. We are interested in fragments of $\mathcal{FL}_\sigma$ in all the languages containing the connectives of additive disjunction ($\vee$), multiplicative conjunction ($*$), 0 and 1. Let $\Psi$ be one of the languages $\langle\vee,*,0,1\rangle$, $\langle\vee,\wedge,*,0,1\rangle$, $\langle\vee,*,{}^{\backprime},{}^{\prime},0,1\rangle$, or $\langle\vee,\wedge,*,{}^{\backprime},{}^{\prime},0,1\rangle$ (the connectives ${}^{\backprime}$ and ${}^{\prime}$ are the right negation and the left negation, respectively). Fix a sequence, $\sigma$ (possibly empty); let $\mathcal{FL}_\sigma[\Psi]$ be the Gentzen systems defined by the calculi $\mathbf{FL}_\sigma[\Psi]$ obtained by dropping from $\mathbf{FL}_\sigma$ the rules for the connectives that are not in $\Psi$; and let $\mathfrak{e}\mathcal{FL}_\sigma[\Psi]$ be its external system. We prove that the external systems $\mathfrak{e}\mathcal{FL}_\sigma[\Psi]$ associated with the Gentzen systems $\mathcal{FL}_\sigma[\Psi]$ are the $[\Psi]$-fragments of $\mathfrak{e}\mathcal{FL}_\sigma$. The steps required to prove this are the following:

A) In Part 2 we introduce, for every $\sigma$, the classes of algebras $\mathring{\mathbb{M}}^{s\ell}_\sigma$ (semilatticed pointed $\sigma$-monoids), $\mathring{\mathbb{M}}^{\ell}_\sigma$ (latticed pointed $\sigma$-monoids), $\mathbb{PM}^{s\ell}_\sigma$ (semilatticed pseudocomplemented $\sigma$-monoids) and $\mathbb{PM}^{\ell}_\sigma$ (latticed pseudocomplemented $\sigma$-monoids). In these acronyms, subindex $\sigma$ is a subsequence of $ew_lw_rc$ and the symbols $e$, $w_l$, $w_r$ and $c$ codify what we refer to as (algebraic) exchange, right-weakening, left-weakening and contraction properties, respectively. Such properties, which are expressed by quasi-inequations, are equivalent, respectively, to the algebraic properties of commutativity, integrality, 0-boundedness, and increasing idempotency, which can all be expressed by equations. Therefore, these classes of algebras are varieties.

B) In Part 3 we prove that the subsystems of $\mathcal{FL}_\sigma$:

$$\mathcal{FL}_\sigma[\vee,*,0,1], \mathcal{FL}_\sigma[\vee,\wedge,*,0,1], \mathcal{FL}_\sigma[\vee,*,{}^{\backprime},{}^{\prime},0,1] \text{ and } \mathcal{FL}_\sigma[\vee,\wedge,*,{}^{\backprime},{}^{\prime},0,1] \tag{1}$$

are algebraizable, and have the varieties $\mathring{\mathbb{M}}^{s\ell}_\sigma$, $\mathring{\mathbb{M}}^{\ell}_\sigma$, $\mathbb{PM}^{s\ell}_\sigma$ and $\mathbb{PM}^{\ell}_\sigma$ as equivalent algebraic semantics, respectively. We also prove that the systems $\mathcal{FL}_\sigma$ are algebraizable with the varieties $\mathbb{FL}_\sigma$ as their equivalent algebraic semantics. We recap on the notion of the algebraization of Gentzen systems in Part 1.

C) In Part 4, by applying the ideal completion method to the classes $\mathring{\mathbb{M}}^{s\ell}_\sigma$, $\mathring{\mathbb{M}}^{\ell}_\sigma$, $\mathbb{PM}^{s\ell}_\sigma$ and $\mathbb{PM}^{\ell}_\sigma$, we prove that for every $\sigma$, they are the subreducts in the corresponding languages of the $\mathbb{FL}_\sigma$ class, i.e., of the variety of pointed residuated lattices defined by the equations codified by $\sigma$.

D) Finally, in Part 5, by using the results of the algebraization in Part 3 and those concerning subreducts obtained in Part 4, we show that, for every $\sigma$, the systems (1) are fragments of $\mathcal{FL}_\sigma$; and that the corresponding external systems are fragments of $\mathfrak{e}\mathcal{FL}_\sigma$.

We also show that each system $\mathcal{FL}_\sigma$ is equivalent to its associated external system. However, it is shown that the fragments considered are not equivalent to any Hilbert system. Moreover, we show that $\mathfrak{e}\mathcal{FL}_\sigma[\vee, *, 0, 1]$, $\mathfrak{e}\mathcal{FL}_\sigma[\vee, \wedge, *, 0, 1]$, $\mathfrak{e}\mathcal{FL}_\sigma[\vee, *, `, ', 0, 1]$ and $\mathfrak{e}\mathcal{FL}_\sigma[\vee, \wedge, *, `, ', 0, 1]$ are not algebraizable (in fact they are not even protoalgebraic) but they have, respectively, the varieties $\mathring{\mathbb{M}}^{s\ell}_\sigma$, $\mathring{\mathbb{M}}^{\ell}_\sigma$, $\mathbb{PM}^{s\ell}_\sigma$ and $\mathbb{PM}^{\ell}_\sigma$ as algebraic semantics with defining equation $1 \vee p \approx p$.

## Part 1. THE LOGICAL SYSTEMS

This part consists of three sections. In Section 2 we recall some notions and results on Gentzen systems that will be used throughout the paper. In Section 3 we present the calculus **FL** and its substructural extensions. Finally, in Section 4 we recap some Hilbert style axiomatizations for the external systems associated with calculi $\mathbf{FL}_\sigma$ that are found in the literature.

### 2. Basic Notions and Results Concerning Gentzen Systems.

Most of the literature on Gentzen systems, and on Hilbert systems, only focuses on their derivable sequents, i.e., on the sequents that can be derived without any hypotheses. In our approach, we analyze the full consequence relation admitting hypotheses in the proofs. The reader should bear in mind this difference between our approach and the one commonly considered in the literature. In this section we recall the notions and results which, from the more general perspective, will be needed later. The results in this section are presented without any proofs; the reader interested in the proofs (or a more detailed presentation) can check Rebagliato & Verdú (1993); Gil et al. (1997) or the monograph by Rebagliato and Verdú on Gentzen Systems: Rebagliato & Verdú (1995). See also Raftery (2006).

2.1. **Sequences.** Given a set $A$ and $n \in \omega$, we can define a finite *sequence*, or *list*, of $n$ elements in $A$ (also called *n-tuple*) as an element in the set $A^n$; that is, as a function from $\{0, \ldots, n-1\}$ to $A$. We will represent finite sequences of $n$ elements in $A$ by $\emptyset$ when $n = 0$ (the *empty sequence*, the sequence of 0 elements in $A$), and by $\langle a_0, \ldots, a_{n-1}\rangle$, or simply by $a_0, \ldots, a_{n-1}$, when $n > 0$. Sometimes we also use the abbreviation $\vec{a}$ for $\langle a_0, \ldots, a_{n-1}\rangle$. The set of all finite sequences of elements in $A$ will be denoted by $A^{<\omega}$. Given a finite sequence $\varsigma$, we denote by $\mathrm{Set}(\varsigma)$ the set whose elements are those appearing in $\varsigma$.

2.2. **Formulas, sequences of formulas.** Let $\mathfrak{L}$ be a propositional language. Let $Var$ be a fixed infinite countable set of propositional variables. The absolutely free algebra of type $\mathfrak{L}$ and generators $Var$ will be denoted by $\mathbf{Fm}_\mathfrak{L}$ and its carrier by $Fm_\mathfrak{L}$. We refer to the elements of $Fm_\mathfrak{L}$ as $\mathfrak{L}$-*formulas*. An $\mathfrak{L}$-*substitution* is an endomorphism on $\mathbf{Fm}_\mathfrak{L}$. We will use the Greek uppercase letters $\Gamma, \Lambda, \Theta, \Delta, \Sigma, \Pi$, possibly with a subindex, for arbitrary elements of $Fm_\mathfrak{L}^{<\omega}$. We use the coma ‘,’ for the operation of concatenation of sequences;

thus given the sequences $\Gamma_1 = \langle \varphi_0, \dots, \varphi_{m-1} \rangle$, $\Gamma_2 = \langle \psi_0, \dots, \psi_{n-1} \rangle$ we write $\Gamma_1, \Gamma_2$ for the sequence

$$\langle \varphi_0, \dots, \varphi_{m-1}, \psi_0, \dots, \psi_{n-1} \rangle.$$

Given an $\mathfrak{L}$-substitution $\sigma$ and a sequence $\Gamma = \varphi_1, \dots, \varphi_m$, we denote by $\sigma\Gamma$, or by $\sigma\vec{\varphi}$, the sequence $\sigma\varphi_, \dots, \sigma\varphi_m$. The set of finite sequences of $\mathfrak{L}$-formulas can be identified with the carrier of the *free monoid* on the set $Fm_{\mathfrak{L}}$ provided with the operation of concatenation of sequences and having the empty sequence as the unit element.

2.3. **Sequents.** Given $m, n \in \omega$, a *sequent of type* $\langle m, n \rangle$, or $\langle m, n \rangle$-sequent, on a set $A$ is an expression of the form $\vec{a} \Rightarrow \vec{b}$, where $\vec{a}$ and $\vec{b}$ are two finite sequences of elements of $A$ such that the length of $\vec{a}$ is $m$ and the length of $\vec{b}$ is $n$. A *trace* is any non-empty subset $\mathcal{T}$ of $\omega \times \omega$. Thus, a trace is a set of possible types, and a sequent of type $\langle m, n \rangle$ on a set $A$ will be called a $\mathcal{T}$-sequent of $A$ if $\langle m, n \rangle \in \mathcal{T}$. We denote by $Seq_A^{\mathcal{T}}$ the set of all $\mathcal{T}$-sequents on $A$, and by $Seq_{\mathfrak{L}}^{\mathcal{T}}$ the set of all $\mathcal{T}$-sequents on $Fm_{\mathfrak{L}}$. The elements of $Seq_{\mathfrak{L}}^{\mathcal{T}}$ will be called *formal sequents* and we refer to them as $\langle \mathfrak{L}, \mathcal{T} \rangle$-sequents. Given a set, $\Phi$, of formal sequents and an $\mathfrak{L}$-substitution, $\sigma$, we denote by $\sigma[\Phi]$ the set of all the sequents $\sigma\vec{\varphi} \Rightarrow \sigma\vec{\psi}$, where $\vec{\varphi} \Rightarrow \vec{\psi} \in \Phi$.

2.4. **Consequence relations.** Let $A$ be a set and let $\mathcal{P}(A)$ be its power set. A *closure operator* on $A$ is a mapping $C : \mathcal{P}(A) \to \mathcal{P}(A)$ such that, for each $X, Y \subseteq A$, the following conditions are satisfied:

(1) $X \subseteq C(X)$,
(2) if $X \subseteq Y$, then $C(X) \subseteq C(Y)$,
(3) $C(C(X)) \subseteq C(X)$.

A closure operator $C$ on $A$ is *finitary* if $C(X) = \bigcup\{C(F) : F \subseteq X,\ F \text{ finite}\}$ for every $X \subseteq A$. We say that a subset $X \subseteq A$ is *$C$-closed* if $C(X) = X$.

Usually, a *consequence relation* on set $A$ is defined as a relation $\vdash$ between subsets of $A$ and elements of $A$ (we write $X \vdash a$ instead of $\langle X, a \rangle \in \vdash$) such that the mapping $\mathcal{P}(A) \to \mathcal{P}(A)$ defined by $C_{\vdash}(X) = \{a \in A : X \vdash a\}$ is a closure operator (called the *closure operator associated with* $\vdash$). This is equivalent to saying that the relation $\vdash$ satisfies the following two properties:

- *Generalized reflexivity:* $a \in X$ implies $X \vdash a$.
- *Transitivity:* $X \vdash a$ and $Y \vdash b$ for every $b \in X$ implies $Y \vdash a$.

Notice that these properties imply the following additional one:

- *Monotonicity:* $X \vdash a$ and $X \subseteq Y$ implies $Y \vdash a$,

We have that $X$ is $C_{\vdash}$-closed if and only if $X = \{a \in A : X \vdash a\}$. A consequence relation $\vdash$ on $A$ is *finitary* if the associated closure operator is finitary; that is, if for every $X \subseteq A$, $C_{\vdash}(X) = \bigcup\{C_{\vdash}(X') : X' \subseteq X,\ X' \text{ finite}\}$ or, equivalently,

$$X \vdash a \ \text{ if and only if } \ X' \vdash a \ \text{ for some } X' \subseteq X, \ \ X' \text{ finite.}$$

To say that $X \vdash a$ is not satisfied, we write $X \nvdash a$. As usual, we write $X,\, a \vdash b$ instead of $X \cup \{a\} \vdash b$, and sometimes $\vdash b$ instead of $\emptyset \vdash b$. The elements of the set $C_{\vdash}(\emptyset)$ are called the *theorems* of the consequence relation. Also, we usually use the sequential notation $a_0, \dots, a_{n-1} \vdash b$ instead of $\{a_0, \dots, a_{n-1}\} \vdash b$ on the implicit understanding that such sequences behave as sets (roughly speaking, the comma is associative, commutative and idempotent).

2.5. **Hilbert systems, Gentzen systems.** Let us recall that, given a propositional language $\mathcal{L}$, a *Hilbert system* on $\mathcal{L}$ is a pair $\mathcal{H} = \langle \mathcal{L}, \vdash_{\mathcal{H}} \rangle$, where $\vdash_{\mathcal{H}}$ is a consequence relation on $Fm_{\mathcal{L}}$ which is invariant under substitutions; that is, for every set of formulas $\Phi \cup \{\varphi\}$ and every substitution $\sigma$, if $\Phi \vdash_{\mathcal{H}} \varphi$, then $\sigma[\Phi] \vdash_{\mathcal{H}} \sigma\varphi$. The relation $\vdash_{\mathcal{H}}$ is called a *Hilbert relation.*

**Definition 1** (Logic in $Seq^{\mathcal{T}}_{\mathcal{L}}$)**.** *A* Gentzen system *is a triple* $\mathcal{G} = \langle \mathcal{L}, \mathcal{T}, \vdash_{\mathcal{G}} \rangle$*, where* $\mathcal{L}$ *is a propositional language,* $\mathcal{T}$ *is a trace, and* $\vdash_{\mathcal{G}}$ *is a consequence relation on* $Seq^{\mathcal{T}}_{\mathcal{L}}$ *which is invariant under substitutions, that is to say that, given a set* $\Phi \cup \{\vec{\varphi} \Rightarrow \vec{\psi}\} \subseteq Seq^{\mathcal{T}}_{\mathcal{L}}$ *of formal sequents,*

$$\text{if } \Phi \vdash_{\mathcal{G}} \vec{\varphi} \Rightarrow \vec{\psi}, \text{ then } \sigma[\Phi] \vdash_{\mathcal{G}} \sigma\vec{\varphi} \Rightarrow \sigma\vec{\psi}, \text{ for every } \mathfrak{L}\text{-substitution, } \sigma.$$

*The relation* $\vdash_{\mathcal{G}}$ *is called a* Gentzen relation of trace $\mathcal{T}$*.* □

Notice that Gentzen systems can be seen as generalizations of Hilbert systems by identifying $\mathcal{L}$-formulas with $\mathcal{L}$-sequents of trace $\{\langle 0, 1 \rangle\}$. A Gentzen system is *finitary* when, for every $\Phi \cup \{\vec{\varphi} \Rightarrow \vec{\psi}\} \subseteq Seq^{\mathcal{T}}_{\mathcal{L}}$ such that $\Phi \vdash_{\mathcal{G}} \vec{\varphi} \Rightarrow \vec{\psi}$, there exists a finite subset $\Phi' \subseteq \Phi$ such that $\Phi' \vdash_{\mathcal{G}} \vec{\varphi} \Rightarrow \vec{\psi}$. In this paper we will deal only with finitary systems.

2.6. **Sequent calculi.** An $\langle \mathcal{L}, \mathcal{T} \rangle$*-rule* is a set $r \subseteq \mathcal{P}_{fin}(Seq^{\mathcal{T}}_{\mathcal{L}}) \times Seq^{\mathcal{T}}_{\mathcal{L}}$ [1] that is obtained as the closure under substitutions of a pair $\langle \Phi, \varsigma \rangle$ such that $\Phi$ is a finite subset of $\mathcal{L}$-sequents (i.e., $\Phi \in \mathcal{P}_{fin}(Seq^{\mathcal{T}}_{\mathcal{L}})$) and $\varsigma$ is an $\mathfrak{L}$-sequent. We will use the pair $\langle \Phi, \varsigma \rangle$ as a name for the rule that it generates. The rules $\langle \emptyset, \varsigma \rangle$ are called *axioms*, and then $\varsigma$ is called an *instance* of the axiom. A rule $r = \langle \Phi, \varsigma \rangle$ is *derivable* in a Gentzen system $\mathcal{G} = \langle \mathcal{L}, \mathcal{T}, \vdash_{\mathcal{G}} \rangle$ if $\Phi \vdash_{\mathcal{G}} \varsigma$; in this case it is also said that $\mathcal{G}$ *satisfies* the rule $r$. Generally we will write a rule $\langle \Phi, \varsigma \rangle$ as $\dfrac{\Phi}{\varsigma}$. An $\langle \mathcal{L}, \mathcal{T} \rangle$*-sequent calculus* is a set of $\langle \mathcal{L}, \mathcal{T} \rangle$-rules. Every $\langle \mathcal{L}, \mathcal{T} \rangle$-sequent calculus, $\mathcal{C}$, determines a Gentzen system $\mathcal{G}_{\mathcal{C}} = \langle \mathcal{L}, \mathcal{T}, \vdash_{\mathcal{C}} \rangle$ in the following way: given $\Phi \cup \{\varsigma\} \subseteq Seq^{\mathcal{T}}_{\mathcal{L}}$, $\Phi \vdash_{\mathcal{C}} \varsigma$ if and only if there is a finite sequence $\varsigma_0, \dots, \varsigma_{n-1}$ of $Seq^{\mathcal{T}}_{\mathcal{L}}$ (which is called a *proof* in $\mathcal{C}$ of $\varsigma$ from $\Phi$) such that $\varsigma_{n-1} = \varsigma$ and for each $i < n$ one of the following conditions hold:

- $\varsigma_i$ is an instance of an axiom of $\mathcal{C}$,
- $\varsigma_i \in \Phi$,
- $\varsigma_i$ is obtained from $\{\varsigma_j : j < i\}$ by using a rule, $r$, of $\mathcal{C}$.

In this case we will say that $\mathcal{G}_{\mathcal{C}}$ is the *Gentzen system determined by the sequent calculus* $\mathcal{C}$.

[1] Given a set $A$, we denote by $\mathcal{P}_{fin}(A)$ the set of finite subsets of $A$.

In this work we consider two kinds of rules: *structural rules* and *rules of introduction for the connectives*. In the *structural rules* the connectives of the language do not appear explicitly; their common characteristic is that they admit instances formed only by variables.

The calculi that we will consider have their set of types between $\omega \times \{1\}$, $\omega \times \{0,1\}$ or $\omega \times \omega$. Moreover, all the calculi considered will have, among their structural rules,[2], the axiom

$$\varphi \Rightarrow \varphi \quad (\mathit{Axiom})$$

and the version of the rule $(Cut)$ appropriate for the corresponding set of types. Thus,

- the $\langle \mathcal{L}, \omega \times \omega \rangle$-rule $(Cut)$ has the form:
$$\frac{\Gamma \Rightarrow \Lambda, \varphi, \Theta \qquad \Sigma, \varphi, \Pi \Rightarrow \Delta}{\Sigma, \Gamma, \Pi \Rightarrow \Lambda, \Delta, \Theta}\,,$$
where all the Greek capital letters represent finite sequences of formulas of any length;
- the $\langle \mathcal{L}, \omega \times \{0,1\} \rangle$-rule $(Cut)$ has the form:
$$\frac{\Gamma \Rightarrow \varphi \qquad \Sigma, \varphi, \Pi \Rightarrow \Delta}{\Sigma, \Gamma, \Pi \Rightarrow \Delta}\,,$$
where $\Gamma$, $\Sigma$ and $\Pi$ are finite sequences of formulas of any length and $\Delta$ is a sequence of one formula at most;
- the $\langle \mathcal{L}, \omega \times \{1\} \rangle$-rule $(Cut)$ has the form:
$$\frac{\Gamma \Rightarrow \varphi \qquad \Sigma, \varphi, \Pi \Rightarrow \psi}{\Sigma, \Gamma, \Pi \Rightarrow \psi}\,,$$
where $\Gamma$, $\Sigma$ and $\Pi$ are finite sequences of formulas of any length.

Other structural rules that we will consider are the following:

| | Left | Right |
|---|---|---|
| Exchange: | $\dfrac{\Gamma, \varphi, \psi, \Pi \Rightarrow \Delta}{\Gamma, \psi, \varphi, \Pi \Rightarrow \Delta}\ (e \Rightarrow)$ | $\dfrac{\Gamma \Rightarrow \Lambda, \varphi, \psi, \Theta}{\Gamma \Rightarrow \Lambda, \psi, \varphi, \Theta}\ (\Rightarrow e)$ |
| Weakening: | $\dfrac{\Gamma, \Pi \Rightarrow \Delta}{\Gamma, \varphi, \Pi \Rightarrow \Delta}\ (w \Rightarrow)$ | $\dfrac{\Gamma \Rightarrow \Lambda, \Theta}{\Gamma \Rightarrow \Lambda, \psi, \Theta}\ (\Rightarrow w)$ |
| Contraction: | $\dfrac{\Gamma, \varphi, \varphi, \Pi \Rightarrow \Delta}{\Gamma, \varphi, \Pi \Rightarrow \Delta}\ (e \Rightarrow)$ | $\dfrac{\Gamma \Rightarrow \Lambda, \varphi, \varphi, \Theta}{\Gamma \Rightarrow \Lambda, \varphi, \Theta}\ (\Rightarrow e)$ |

Note that if the set of types is $\omega \times \{0,1\}$ or $\omega \times \{1\}$, the rules $(\Rightarrow e)$ and $(\Rightarrow c)$ are not expressible. Observe also that if the set of types is $\omega \times \{1\}$, the rule $(\Rightarrow w)$ is not expressible. If the set of types is $\omega \times \{0,1\}$, the rule $(\Rightarrow w)$ takes the form:

[2] Strictly speaking, the rules that we present in what follows are families of rules and not only individual rules.

$$\frac{\Gamma \Rightarrow \emptyset}{\Gamma \Rightarrow \psi}\ (\Rightarrow w)$$

A *rule of introduction* is a rule $\langle \Phi, \varsigma \rangle$ such that every formula in $\Phi$ is a subformula of some formula in $\varsigma$ and, moreover, there is a formula in $\Phi$ that is a proper subformula of some formula in $\varsigma$.

**Definition 2** (Regular sequent calculus). *We say that an $\langle \mathcal{L}, \mathcal{T} \rangle$-sequent calculus, $\mathcal{C}$, is* regular *if and only if $\langle 1, 1 \rangle \in \mathcal{T}$ and $\mathcal{C}$ contains only structural rules together with two families of rules of introduction for each connective such that, if k is the arity of a connective $\iota$ of $\mathcal{L}$, the generators of the rules of these two families are of the form:*

$$\frac{S}{\Gamma, \iota(\varphi_1, \ldots, \varphi_k), \Pi \Rightarrow \Delta}, \qquad \frac{T}{\Gamma \Rightarrow \Lambda, \iota(\varphi_1, \ldots, \varphi_k), \Theta}$$

*in such a way that the sequences in the sequents of the sets $S$ and $T$ and the sequences $\Gamma$, $\Pi$, $\Delta$, $\Lambda$ and $\Theta$ contain only variables. For these two families of rules we use the labels $(\iota \Rightarrow)$ and $(\Rightarrow \iota)$ and we call them respectively* rule of introduction to the left *and* rule of introduction to the right *for the connective $\iota$. We will say that a Gentzen system is* regular *if it is defined by a sequent calculus which is regular.*

**Definition 3** (Restrictions of a regular system). *Given a regular Gentzen system $\mathcal{G} = \langle \mathcal{L}, \mathcal{T}, \vdash_{\mathcal{G}} \rangle$, i.e., one defined by a regular sequent calculus $\mathcal{C}$, if $\Psi$ is a sublanguage of $\mathcal{L}$, we will denote by $\mathcal{C}[\Psi]$ the calculus in the language $\Psi$ obtained by dropping from the calculus $\mathcal{C}$ the introduction rules for the connectives that do not belong to $\Psi$. We denote by $\mathcal{G}[\Psi]$ the Gentzen system defined by the calculus $\mathcal{C}[\Psi]$. If $\Psi = \langle \iota_1, \ldots, \iota_n \rangle$, we will use the explicit notations $\mathcal{C}[\iota_1, \ldots, \iota_n]$ and $\mathcal{G}[\iota_1, \ldots, \iota_n]$ for $\mathcal{C}[\Psi]$ and $\mathcal{G}[\Psi]$, respectively. If $\Psi' = \langle \Psi, \iota_1, \ldots, \iota_n \rangle$ is a sublanguage of $\mathcal{L}$ which is obtained by adding to $\Psi$ the connectives $\iota_1, \ldots, \iota_n$, we will use the notations $\mathcal{C}[\Psi, \iota_1, \ldots, \iota_n]$ and $\mathcal{G}[\Psi, \iota_1, \ldots, \iota_n]$ for $\mathcal{C}[\Psi']$ and $\mathcal{G}[\Psi']$, respectively.*

2.7. **Expansions, extensions, fragments.** Let $\mathcal{G} = \langle \mathcal{L}, \mathcal{T}, \vdash_{\mathcal{G}} \rangle$ and $\mathcal{G}' = \langle \mathcal{L}', \mathcal{T}', \vdash_{\mathcal{G}'} \rangle$ be such that $\mathcal{L} \leq \mathcal{L}'$ and $\mathcal{T} \subseteq \mathcal{T}'$. We will say that $\mathcal{G}'$ is an *expansion* of $\mathcal{G}$ if and only if $\vdash_{\mathcal{G}} \subseteq \vdash_{\mathcal{G}'}$. In the case that $\mathcal{L} = \mathcal{L}'$ and $\mathcal{T} = \mathcal{T}'$ we will say that $\mathcal{G}'$ is an *extension* of $\mathcal{G}$. An expansion $\mathcal{G}' = \langle \mathcal{L}', \mathcal{T}', \vdash_{\mathcal{G}'} \rangle$ of $\mathcal{G} = \langle \mathcal{L}, \mathcal{T}, \vdash_{\mathcal{G}} \rangle$ is an *axiomatic expansion* if and only if $\mathcal{T} = \mathcal{T}'$ and the systems $\mathcal{G}$ and $\mathcal{G}'$ are respectively defined by two sequent calculi, $\mathcal{C}$ and $\mathcal{C}'$, such that: a) they have the same proper rules, and b) every axiom of $\mathcal{C}$ is an axiom of $\mathcal{C}'$ (i.e., $\mathcal{C}'$ is obtained from $\mathcal{C}$ by adding only axioms). If, moreover, $\mathcal{L} = \mathcal{L}'$, then we say that $\mathcal{G}'$ is an *axiomatic extension* of $\mathcal{G}$.

An expansion $\mathcal{G}' = \langle \mathcal{L}', \mathcal{T}', \vdash_{\mathcal{G}'} \rangle$ of $\mathcal{G} = \langle \mathcal{L}, \mathcal{T}, \vdash_{\mathcal{G}} \rangle$ is *conservative* if and only if:

$$\text{for every } \Phi \cup \{\varsigma\} \subseteq Seq_{\mathcal{L}}^{\mathcal{T}},\ \Phi \vdash_{\mathcal{G}} \varsigma \ \text{ iff } \ \Phi \vdash_{\mathcal{G}'} \varsigma.$$

If $\mathcal{G}'$ is a conservative expansion of $\mathcal{G}$, we also say that $\mathcal{G}$ is the $\langle \mathcal{L}, \mathcal{T} \rangle$-*fragment* of $\mathcal{G}'$. If $\mathcal{T} = \mathcal{T}'$, we say that $\mathcal{G}$ is the $\mathcal{L}$-*fragment* of $\mathcal{G}'$ and, if $\mathcal{L} = \mathcal{L}'$, we say that $\mathcal{G}$ is the $\mathcal{T}$-*fragment* of $\mathcal{G}'$. We will use the notations $\langle \mathcal{L}, \mathcal{T} \rangle$-$\mathcal{G}'$, $\mathcal{L}$-$\mathcal{G}'$ and $\mathcal{T}$-$\mathcal{G}'$ for the $\langle \mathcal{L}, \mathcal{T} \rangle$-fragment, the $\mathcal{L}$-fragment, and the $\mathcal{T}$-fragment of $\mathcal{G}'$, respectively.

2.8. **Equivalence of Gentzen systems.** Let $\mathcal{G}$ and $\mathcal{G}'$ be two Gentzen systems $\langle \mathcal{L}, \mathcal{T}, \vdash \rangle$ and $\langle \mathcal{L}', \mathcal{T}', \vdash' \rangle$ such that $\mathcal{L} = \mathcal{L}'$. An $\langle \mathcal{L}, \mathcal{T}, \mathcal{T}' \rangle$-*translation* is a map

$$\tau : Seq_{\mathcal{L}}^{\mathcal{T}} \longrightarrow \mathcal{P}_{fin}(Seq_{\mathcal{L}}^{\mathcal{T}'})$$

such that:

- for every $\langle m, n \rangle \in \mathcal{T}$, if $\langle m, n \rangle \neq \langle 0, 0 \rangle$ then the sequents in the set
$$\tau(p_0, \ldots, p_{m-1} \Rightarrow q_0, \ldots, q_{n-1})$$
only use variables in $\{p_0, \ldots, p_{m-1}, q_0, \ldots, q_{n-1}\}$,
- for every $\langle m, n \rangle \in \mathcal{T}$, $\varphi_0, \ldots, \varphi_{m-1} \in Fm_{\mathcal{L}}$ and $\psi_0, \ldots, \psi_{n-1} \in Fm_{\mathcal{L}}$, if $\langle m, n \rangle \neq \langle 0, 0 \rangle$ then
$$\tau(\varphi_0, \ldots, \varphi_{m-1} \Rightarrow \psi_0, \ldots, \psi_{n-1}) = e[\tau(p_0, \ldots, p_{m-1} \Rightarrow q_0, \ldots, q_{n-1})],$$
where $e$ is the substitution such that $e(p_i) = \varphi_i$ and $e(q_i) = \psi_i$.

¿From the above conditions it is clear that the map $\tau$ is determined[3] by the restriction of $\tau$ to the set:

$$\{\tau(p_0, \ldots, p_{m-1} \Rightarrow q_0, \ldots, q_{n-1}) : \langle m, n \rangle \in \mathcal{T}\}.$$

It is said that the Gentzen systems $\mathcal{G}$ and $\mathcal{G}'$ are *equivalent* if there is an $\langle \mathcal{L}, \mathcal{T}, \mathcal{T}' \rangle$-translation $\tau : Seq_{\mathcal{L}}^{\mathcal{T}} \longrightarrow \mathcal{P}_{fin}(Seq_{\mathcal{L}}^{\mathcal{T}'})$ and an $\langle \mathcal{L}, \mathcal{T}', \mathcal{T} \rangle$-translation $\rho : Seq_{\mathcal{L}}^{\mathcal{T}'} \longrightarrow \mathcal{P}_{fin}(Seq_{\mathcal{L}}^{\mathcal{T}})$ such that:

1) for all $\Phi \cup \{\varsigma\} \subseteq Seq_{\mathcal{L}}^{\mathcal{T}}$, it holds that $\Phi \vdash \varsigma$ iff $\tau[\Phi] \vdash' \tau(\varsigma)$,
2) for all $\Phi \cup \{\varsigma\} \subseteq Seq_{\mathcal{L}}^{\mathcal{T}'}$, it holds that $\Phi \vdash' \varsigma$ iff $\rho[\Phi] \vdash \rho(\varsigma)$,
3) for all $\varsigma \in Seq_{\mathcal{L}}^{\mathcal{T}'}$, it holds that $\varsigma \dashv\vdash' \tau\rho(\varsigma)$,
4) for all $\varsigma \in Seq_{\mathcal{L}}^{\mathcal{T}}$, it holds that $\varsigma \dashv\vdash \rho\tau(\varsigma)$.

It is known that the previous definition is redundant because the conjunction of 1) and 3) is equivalent to the conjunction of 2) and 4) (Rebagliato & Verdú, 1995, Proposition 2.1).

2.9. **Theories, matrices, models and filters in a Gentzen system.** Assume that a Gentzen system $\mathcal{G} = \langle \mathcal{L}, \mathcal{T}, \vdash \rangle$ is fixed. A subset $\Phi$ of $Seq_{\mathcal{L}}^{\mathcal{T}}$ is a $\mathcal{G}$-*theory* if $\Phi \vdash \varsigma$ implies $\varsigma \in \Phi$. The set of all $\mathcal{G}$-theories is denoted by $Th\,\mathcal{G}$. An $\langle \mathcal{L}, \mathcal{T} \rangle$-*matrix*, or simply a *matrix*, is a pair $\langle \mathbf{A}, R \rangle$ where $\mathbf{A}$ is an $\mathcal{L}$-algebra and $R \subseteq \bigcup\{A^m \times A^n : \langle m, n \rangle \in \mathcal{T}\}$. Every $\langle \mathcal{L}, \mathcal{T} \rangle$-matrix allows us to introduce a Gentzen system with language $\mathcal{L}$ and type $\mathcal{T}$: given a matrix $\langle \mathbf{A}, R \rangle$, just consider the consequence relation $\models_{\langle \mathbf{A}, R \rangle}$ defined by $\Phi \models_{\langle \mathbf{A}, R \rangle} \varsigma$ if and only if:

for every homomorphism, $h$, from $\mathbf{Fm}_{\mathcal{L}}$ into $\mathbf{A}$; if $h[\Phi] \subseteq R$, then $h(\varsigma) \in R$.

It is easily verified that this consequence relation is invariant under substitutions. Hence, we have that $\langle \mathcal{L}, \mathcal{T}, \models_{\langle \mathbf{A}, R \rangle} \rangle$ is a Gentzen system (possibly not finitary). In the case that $\vdash \leq \models_{\langle \mathbf{A}, R \rangle}$ (i.e., for every set $\Phi \cup \{\varsigma\}$ of sequents, if $\Phi \vdash \varsigma$ then $\Phi \models_{\langle \mathbf{A}, R \rangle} \varsigma$) then it is said that $\langle \mathbf{A}, R \rangle$ is a $\mathcal{G}$-*model* and that $R$ is a $\mathcal{G}$-*filter* of $\mathbf{A}$. It is well known that whenever $\mathcal{G}$ is defined by means of a sequent calculus, then $R$ is a $\mathcal{G}$-filter if and only if $R$ contains all the interpretations of the axioms and is closed under each of the rules. Another easy remark is that $\Phi$ is a $\mathcal{G}$-theory if and only if $\Phi$ is a $\mathcal{G}$-filter of $\mathbf{Fm}_{\mathcal{L}}$.

[3]Strictly speaking, the map $\tau$ is *quasi*-determined because the value $\tau(\Gamma \Rightarrow \Delta)$ is determined with the exception of the case in which $\Gamma$ and $\Delta$ are both the empty sequence.

2.10. **Algebraization of Gentzen systems.** If $\mathbb{K}$ is a class of $\mathcal{L}$-algebras, then the equational logic $\models_{\mathbb{K}}$ can be seen as a Gentzen system with language $\mathcal{L}$ and set of types $\mathcal{T} = \{1\} \times \{1\}$ where we identify an equation $\varphi \approx \psi \in Eq_{\mathcal{L}}$ with the sequent $\varphi \Rightarrow \psi$. The class $\mathbb{K}$ is an *algebraic semantics* for a Gentzen system $\mathcal{G} = \langle \mathcal{L}, \mathcal{T}, \vdash \rangle$ if there is a translation $\tau : Seq_{\mathcal{L}}^{\mathcal{T}} \longrightarrow \mathcal{P}(Eq_{\mathcal{L}})$ such that for all $\Phi \cup \{\varsigma\} \subseteq Seq_{\mathcal{L}}^{\mathcal{T}}$:

$$\Phi \vdash \varsigma \quad \text{iff} \quad \tau[\Phi] \models_{\mathbb{K}} \tau(\varsigma).$$

A Gentzen system, $\mathcal{G}$, is said to be *algebraizable with equivalent algebraic semantics* $\mathbb{K}$ if $\mathcal{G}$ and $\models_{\mathbb{K}}$ are equivalent Gentzen systems.

It holds that if $\mathbb{K}$ is an equivalent algebraic semantics for $\mathcal{G}$, then so is the quasivariety $\mathbf{Q}(\mathbb{K})$ generated by $\mathbb{K}$ (Gil et al., 1997, Corollary 4.2). It is also known that if $\mathbb{K}$ and $\mathbb{K}'$ are equivalent algebraic semantics for $\mathcal{G}$, then $\mathbb{K}$ and $\mathbb{K}'$ generate the same quasivariety (Gil et al., 1997, Corollary 4.4). This quasivariety is called *the equivalent quasivariety semantic* for $\mathcal{G}$.

It should be noticed that if $\mathcal{H}$ is a Hilbert system, then the fact that it is algebraizable in the sense of Blok & Pigozzi (1989) with the set of equivalence formulas $\Delta(p, q)$ and the set of defining equations $\Theta(p)$ coincides precisely with the fact that it is algebraizable in the above sense under the translations $\tau(p) := \Theta(p)$ and $\rho(p \approx q) := \Delta(p, q)$. Hence, the algebraization of Gentzen systems generalizes the algebraization of Hilbert systems introduced in Blok & Pigozzi (1989).

Now we state a result that we will need in Section 9 It provides a sufficient condition to prove the algebraization of a Gentzen system (Rebagliato & Verdú, 1995, Lemma 2.5) (see also (Gil et al., 1997, Lemma 4.5)). In fact, it is also known that this condition is necessary (see (Rebagliato & Verdú, 1995, Lemma 2.24)).

**Lemma 1.** *Let $\mathcal{G}$ be a Gentzen system $\langle \mathcal{L}, \mathcal{T}, \vdash \rangle$ and $\mathbb{K}$ a quasivariety. Suppose that there are two translations $\tau : Seq_{\mathcal{L}}^{\mathcal{T}} \longrightarrow \mathcal{P}_{fin}(Eq_{\mathcal{L}})$ and $\rho : Eq_{\mathcal{L}} \longrightarrow \mathcal{P}_{fin}(Seq_{\mathcal{L}}^{\mathcal{T}})$ such that:*

1) *for all $\varsigma \in Seq_{\mathcal{L}}^{\mathcal{T}}$, $\varsigma \dashv\vdash_{\mathcal{G}} \rho\tau(\varsigma)$,*
2) *for all $\varphi \approx \psi \in Eq_{\mathcal{L}}$, $\varphi \approx \psi \ \ =\!\mid\models_{\mathbb{K}} \ \tau\rho(\varphi \approx \psi)$,*
3) *for all $\mathbf{A} \in \mathbb{K}$, the set*

$$R := \{\langle \bar{a}, \bar{b} \rangle \in A^m \times A^n : \langle m, n \rangle \in \mathcal{T},\ \mathbf{A} \models \tau(p_0, \dots, p_{m-1} \Rightarrow q_0, \dots, q_{n-1})[\![\bar{a}, \bar{b}]\!]\}$$

*is a $\mathcal{G}$-filter,*[4]

4) *for all $\Phi \in Th\,\mathcal{G}$, the relation*

$$\theta_\Phi := \{\langle \varphi, \psi \rangle \in Fm_{\mathcal{L}}^2 : \rho(\varphi \approx \psi) \subseteq \Phi\}$$

*is a congruence relative to the quasivariety $\mathbb{K}$, i.e., $\mathbf{Fm}_{\mathcal{L}}/\theta_\Phi \in \mathbb{K}$.*

*Then, $\mathcal{G}$ is algebraizable with equivalent algebraic semantics $\mathbb{K}$.*

[4] If $\bar{a} = \langle a_0, \dots, a_{m-1} \rangle$ and $\bar{b} = \langle b_0, \dots, b_{n-1} \rangle$, we denote by $\tau(p_0, \dots, p_{m-1} \Rightarrow q_0, \dots, q_{n-1})[\![\bar{a}, \bar{b}]\!]$ the result obtained by applying the assignation $v$ such that $v(p_0) = a_0, \dots, v(p_{m-1}) = a_{m-1}, v(q_0) = b_0, \dots, v(q_{n-1}) = b_{n-1}$ to the equations of the set $\tau(p_0, \dots, p_{m-1} \Rightarrow q_0, \dots, q_{n-1})$.

2.11. **The Leibniz operator.** One interesting property of algebraizable Gentzen systems with respect to quasivarieties is the existence of a characterization of congruences relative to the quasivariety. To describe this characterization we need the notion of the *Leibniz operator*. Let $\mathbf{A}$ be an $\mathcal{L}$-algebra, and let $\mathcal{T}$ be a set of types. If $m, n \in \omega$, $\langle \bar{x}, \bar{y} \rangle \in A^m \times A^n$ and $a, b \in A$, then $\langle \bar{x}, \bar{y} \rangle(a|b)$ will denote the result of replacing one occurrence (if it exists) of $a$ in $\langle \bar{x}, \bar{y} \rangle$ with $b$. Given an $\langle \mathcal{L}, \mathcal{T} \rangle$-matrix, $\langle \mathbf{A}, R \rangle$, the *Leibniz congruence*, $\Omega_{\mathbf{A}} R$, of the matrix $\langle \mathbf{A}, R \rangle$ is the equivalence relation on $A$ defined in the following way: $\langle a, b \rangle \in \Omega_{\mathbf{A}} R$ if and only if, for every $\langle m, n \rangle \in \mathcal{T}$, $\langle \bar{x}, \bar{y} \rangle \in A^m \times A^n$, $k \in \omega$, $\varphi(p, q_0, \ldots, q_{k-1}) \in Fm_{\mathcal{L}}$ and $c, c_0, \ldots, c_{k-1} \in A$:

$$\langle \bar{x}, \bar{y} \rangle(c|\varphi^{\mathbf{A}}(a, c_0, \ldots, c_{k-1})) \in R \quad \Leftrightarrow \quad \langle \bar{x}, \bar{y} \rangle(c|\varphi^{\mathbf{A}}(b, c_0, \ldots, c_{k-1})) \in R.$$

We emphasize that this definition does not depend on any Gentzen system. It holds that $\Omega_{\mathbf{A}} : \bigcup\{A^m \times A^n : \langle m, n \rangle \in \mathcal{T}\} \longrightarrow Con(\mathbf{A})$. This map is known as the *Leibniz operator* on $\mathbf{A}$. It is easy to show that $\Omega_{\mathbf{A}} R$ is characterized by the fact that it is the largest congruence of $\mathbf{A}$ that is compatible with $R$ (i.e., if $\langle \bar{x}, \bar{y} \rangle \in R$ and $\langle a, b \rangle \in \theta$, then $\langle \bar{x}, \bar{y} \rangle(a|b) \in R$). Let $\mathcal{G} = \langle \mathcal{L}, \mathcal{T}, \vdash_{\mathcal{G}} \rangle$ be a Gentzen system and let $\Phi$ be a $\mathcal{G}$-theory. Then, $\langle \mathbf{Fm}_{\mathcal{L}}, \Phi \rangle$ is an $\langle \mathcal{L}, \mathcal{T} \rangle$-matrix that is a $\mathcal{G}$-model. In this case we denote the Leibniz operator by $\Omega$ instead of $\Omega_{\mathbf{Fm}_{\mathcal{L}}}$.

Next, we give the result mentioned previously which concerns the congruences (Rebagliato & Verdú, 1993, Theorem 2.23) (see also (Gil et al., 1997, Theorem 4.7)).

**Theorem 1.** *Let $\mathcal{G}$ be a Gentzen system and $\mathbb{K}$ a quasivariety. The following statements are equivalent:*

1) *$\mathcal{G}$ is algebraizable with equivalent algebraic semantics $\mathbb{K}$.*
2) *For every, $\mathcal{L}$-algebra $\mathbf{A}$, the Leibniz operator, $\Omega_{\mathbf{A}}$, is an isomorphism between the lattice of $\mathcal{G}$-filters of $\mathbf{A}$ and the lattice of $\mathbb{K}$-congruences of $\mathbf{A}$.*
3) *The Leibniz operator, $\Omega$, is a lattice isomorphism between $Th\,\mathcal{G}$ and $Con_{\mathbb{K}}\mathbf{Fm}_{\mathcal{L}}$.*

2.12. **Hilbert systems associated with a Gentzen system.** Let $\mathcal{G}$ be a Gentzen system $\langle \mathcal{L}, \mathcal{T}, \vdash \rangle$. There are at least two methods in the literature of associating a Hilbert system with $\mathcal{G}$. The most common method is based on considering the derivable sequents. Specifically, $\Sigma \vdash_{\mathcal{I}(\mathcal{G})} \varphi$ holds when:

$$\text{there is a finite subset } \{\varphi_0, \ldots, \varphi_{n-1}\} \text{ of } \Sigma \text{ such that } \emptyset \vdash \varphi_0, \ldots, \varphi_{n-1} \Rightarrow \varphi.$$

This approach yields a Hilbert system, called *internal* if and only if the Gentzen system satisfies the following structural rules: (*Axiom*), (*Cut*), exchange, left weakening and contraction. Another method, which works for all Gentzen systems such that $\langle 0, 1 \rangle \in \mathcal{T}$, yields the so-called *external* system.[5] The *external system* associated with $\mathcal{G}$ is defined as the system $\mathcal{E}(\mathcal{G})$ such that:

$$\Sigma \vdash_{\mathcal{E}(\mathcal{G})} \varphi \quad \text{iff} \quad \{\, \emptyset \Rightarrow \psi : \psi \in \Sigma\} \vdash \emptyset \Rightarrow \varphi.$$

Since we have restricted ourselves to finitary Gentzen systems, it is clear that $\mathcal{E}(\mathcal{G})$ is finitary.

[5] We use the names *internal* and *external* following Avron (see Avron (1988)).

## 3. Basic Intuitionistic Substructural Systems.

In this section we recap the definitions of the basic substructural Gentzen systems, $\mathcal{FL}_\sigma$, and their associated external Hilbert systems, $\mathfrak{e}\mathcal{FL}_\sigma$. The calculi, $\mathbf{FL}_\sigma$, defining systems $\mathcal{FL}_\sigma$ are presented with sequents of type $\omega \times \{0,1\}$ and in the language $\langle \vee, \wedge, *, \backslash, /, ^{\backprime}, ', 0, 1\rangle$ which includes two negations $^{\backprime}$ and $'$ and a constant symbol 0. The systems $\mathcal{FL}_\sigma$ also allow a definitionally equivalent presentation in sequents of type $\omega \times \{1\}$ and in the language without negations $\langle \vee, \wedge, *, \backslash, /, 0, 1\rangle$ (see for instance); but we have decided to adopt the former since we want to study certain fragments without implications that contain such negations. In Section 3.2 we define the notion of a *mirror image of a sequent* and set out the *Law of Mirror Images* for the systems $\mathcal{FL}_\sigma[\Psi]$ such that $\Psi$ contains the two implications or the two negations. In Section 3.3 we characterize the sequential Leibniz congruence of the theories of the systems $\mathcal{FL}[\Psi]$ and the Leibniz congruence of the theories of their associated external systems $\mathfrak{e}\mathcal{FL}[\Psi]$ for the case in which the language $\Psi$ contains one of the implication connectives. We go on to prove that these external systems are protoalgebraic. We arrive at results related to characterizations as a consequence of other more general results established in larger classes of Gentzen systems that $\mathcal{FL}_\sigma[\Psi]$ systems belong to.

3.1. **The FL calculus and its substructural extensions.** We will now recap the definition of the intuitionistic basic substructural calculus **FL** in its version with sequents of type $\omega \times \{0,1\}$ and with the two negations (denoted $^{\backprime}$ and $'$ hereafter) as primitive connectives. It is well known that the calculus **FL** has the variety $\mathbb{FL}$ of the pointed residuated lattices as its algebraic counterpart. As we will see later, this relationship can be established in a strong sense: the system defined by **FL** and the equational system associated with the variety $\mathbb{FL}$ are equivalent as Gentzen systems. The definition of **FL**, introduced by Hiroakira Ono, appears for the first time in Ono (1990) (cf. also Ono (1990, 1993, 1998, 2003b); Galatos & Ono (2006) ).

**Definition 4** (Full Lambek calculus). *Let $\mathfrak{L}$ be the propositional language $\langle \vee, \wedge, *, \backslash, /, ^{\backprime}, ', 0, 1\rangle$ of type $\langle 2,2,2,2,2,1,1,0,0\rangle$. Let $\varphi, \psi$ be $\mathfrak{L}$-formulas; $\Gamma, \Pi, \Sigma$ be finite sequences (possibly empty) of $\mathfrak{L}$-formulas and $\Delta$ be a sequence with at the most one formula. The* Full Lambek calculus **FL** *is the $\langle \mathfrak{L}, \omega \times \{0,1\}\rangle$-calculus defined by the following axioms and rules:*

$$\varphi \Rightarrow \varphi \quad (Axiom) \qquad\qquad \frac{\Gamma \Rightarrow \varphi \qquad \Sigma, \varphi, \Pi \Rightarrow \Delta}{\Sigma, \Gamma, \Pi \Rightarrow \Delta} \quad (Cut)$$

$$\frac{\Sigma, \varphi, \Gamma \Rightarrow \Delta \qquad \Sigma, \psi, \Gamma \Rightarrow \Delta}{\Sigma, \varphi \vee \psi, \Gamma \Rightarrow \Delta} \ (\vee \Rightarrow) \quad \frac{\Gamma \Rightarrow \varphi}{\Gamma \Rightarrow \varphi \vee \psi} \ (\Rightarrow \vee_1) \quad \frac{\Gamma \Rightarrow \psi}{\Gamma \Rightarrow \varphi \vee \psi} \ (\Rightarrow \vee_2)$$

$$\frac{\Sigma, \varphi, \Gamma \Rightarrow \Delta}{\Sigma, \varphi \wedge \psi, \Gamma \Rightarrow \Delta} \ (\wedge_1 \Rightarrow) \quad \frac{\Sigma, \psi, \Gamma \Rightarrow \Delta}{\Sigma, \varphi \wedge \psi, \Gamma \Rightarrow \Delta} \ (\wedge_2 \Rightarrow) \quad \frac{\Gamma \Rightarrow \varphi \qquad \Gamma \Rightarrow \psi}{\Gamma \Rightarrow \varphi \wedge \psi} \ (\Rightarrow \wedge)$$

$$\frac{\Sigma, \varphi, \psi, \Gamma \Rightarrow \Delta}{\Sigma, \varphi * \psi, \Gamma \Rightarrow \Delta} \ (* \Rightarrow) \qquad\qquad \frac{\Gamma \Rightarrow \varphi \qquad \Pi \Rightarrow \psi}{\Gamma, \Pi \Rightarrow \varphi * \psi} \ (\Rightarrow *)$$

Table 1. Connectives of **FL**

| **Symbol** | **Name** |
|---|---|
| $\vee$ | Disjunction |
| $\wedge$ | Conjunction (or additive conjunction) |
| $*$ | Fusion (or multiplicative conjunction) |
| $\backslash$ | Right implication |
| $/$ | Left implication |
| $`$ | Right negation |
| $'$ | Left negation |
| 0 | Falsity (or Zero) |
| 1 | Truth (or One) |

$$\frac{\Gamma \Rightarrow \varphi \qquad \Sigma, \psi, \Pi \Rightarrow \Delta}{\Sigma, \Gamma, \varphi\backslash\psi, \Pi \Rightarrow \Delta}\ (\backslash \Rightarrow) \qquad \frac{\varphi, \Gamma \Rightarrow \psi}{\Gamma \Rightarrow \varphi\backslash\psi}\ (\Rightarrow \backslash)$$

$$\frac{\Gamma \Rightarrow \varphi \quad \Sigma, \psi, \Pi \Rightarrow \Delta}{\Sigma, \psi/\varphi, \Gamma, \Pi \Rightarrow \Delta}\ (/ \Rightarrow) \qquad \frac{\Gamma, \varphi \Rightarrow \psi}{\Gamma \Rightarrow \psi/\varphi}\ (\Rightarrow /)$$

$$\frac{\Gamma \Rightarrow \varphi}{\Gamma, \varphi^{`} \Rightarrow \emptyset}\ ({}^{`} \Rightarrow) \qquad \frac{\varphi, \Gamma \Rightarrow \emptyset}{\Gamma \Rightarrow \varphi^{`}}\ (\Rightarrow {}^{`})$$

$$\frac{\Gamma \Rightarrow \varphi}{{}'\varphi, \Gamma \Rightarrow \emptyset}\ ({}' \Rightarrow) \qquad \frac{\Gamma, \varphi \Rightarrow \emptyset}{\Gamma \Rightarrow {}'\varphi}\ (\Rightarrow {}')$$

$$\frac{\Sigma, \Gamma \Rightarrow \Delta}{\Sigma, 1, \Gamma \Rightarrow \Delta}\ (1 \Rightarrow) \qquad \emptyset \Rightarrow 1\ \ (\Rightarrow 1)$$

$$0 \Rightarrow \emptyset\ \ (0 \Rightarrow) \qquad \frac{\Gamma \Rightarrow \emptyset}{\Gamma \Rightarrow 0}\ (\Rightarrow 0)$$

The connectives of the language $\mathfrak{L}$ of **FL** are denoted by the names given in Table 1.

We now define the extensions of calculus **FL** with different combinations for the structural rules of exchange, weakening (right and left) and contraction.

**Definition 5.** (The calculi $\mathbf{FL}_\sigma$).
$\mathbf{FL}_e$ *is the calculus obtained by adding to the rules of* **FL** *the* $\langle \mathfrak{L}, \omega \times \{0,1\}\rangle$*-rule of exchange:*

$$\frac{\Gamma, \varphi, \psi, \Pi \Rightarrow \Delta}{\Gamma, \psi, \varphi, \Pi \Rightarrow \Delta}\ (e \Rightarrow).$$

$\mathbf{FL}_{w_l}$ *is the calculus obtained by adding to the rules of* **FL** *the* $\langle \mathfrak{L}, \omega \times \{0,1\}\rangle$*-rule of left-weakening:*

$$\frac{\Sigma, \Gamma \Rightarrow \Delta}{\Sigma, \varphi, \Gamma \Rightarrow \Delta}\ (w \Rightarrow).$$

$\mathbf{FL}_{w_r}$ *is the calculus obtained by replacing in* $\mathbf{FL}$ *the rule* $(\Rightarrow 0)$ *with the* $\langle \mathfrak{L}, \omega \times \{0,1\}\rangle$*-rule of right weakening:*

$$\frac{\Gamma \Rightarrow \emptyset}{\Gamma \Rightarrow \varphi}\ (\Rightarrow w).$$

$\mathbf{FL}_w$ *is the calculus obtained by adding to* $\mathbf{FL}$ *the rules* $(w \Rightarrow)$ *and* $(\Rightarrow w)$*.*
$\mathbf{FL}_c$ *is the calculus obtained by adding to* $\mathbf{FL}$ *the* $\langle \mathfrak{L}, \omega \times \{0,1\}\rangle$*-rule of contraction:*

$$\frac{\Sigma, \varphi, \varphi, \Gamma \Rightarrow \Delta}{\Sigma, \varphi, \Gamma \Rightarrow \Delta}\ (c \Rightarrow).$$

*Let* $\sigma$ *be a subsequence of the sequence* $ew_lw_rc$*. For short, the sequence* $w_lw_r$ *will be denoted by* $w$*. We will use* $\mathbf{FL}_\sigma$ *to refer to the calculus obtained by adding to the rules of* $\mathbf{FL}$ *the* $\langle \mathfrak{L}, \omega \times \{0,1\}\rangle$*-rules codified by the letters appearing in* $\sigma$*; when* $\sigma$ *is the empty sequence, then* $\mathbf{FL}_\sigma$ *is the calculus* $\mathbf{FL}$*. We will call the* $\mathbf{FL}_\sigma$ *calculi* basic substructural calculi*.*

So, for example, $\mathbf{FL}_{w_rc}$ is the calculus obtained by adding to the rules of $\mathbf{FL}$ the structural rules $(\Rightarrow w)$ and $(c \Rightarrow)$.

**Remark 1.** *Note that in the calculi with the rule of left-weakening,* $(1 \Rightarrow)$ *is an instance of* $(w \Rightarrow)$*; and that in calculi with the rule of right-weakening,* $(\Rightarrow 0)$ *is an instance of* $(\Rightarrow w)$*.*

**Definition 6.** (Systems $\mathcal{FL}_\sigma$ and $\mathfrak{e}\mathcal{FL}_\sigma$). *Let* $\sigma \leq ew_lw_rc$.[6] *The Gentzen system defined by a calculus* $\mathbf{FL}_\sigma$ *will be denoted by* $\mathcal{FL}_\sigma$*. The systems* $\mathcal{FL}_\sigma$ *will be called* basic substructural Gentzen systems*. The external Hilbert system associated with an* $\mathcal{FL}_\sigma$ *will be denoted by* $\mathfrak{e}\mathcal{FL}_\sigma$*.*

Let us recall that a calculus with the cut rule $(Cut)$ satisfies the *Cut Elimination Theorem* if, for any sequent $\varsigma$ derivable in the calculus, there exists a proof of $\varsigma$ in which the cut rule is not used. A calculus satisfies the *subformula property* if, for any sequent $\varsigma$ derivable in the calculus, there exists a derivation of $\varsigma$ in which every formula included in it is a subformula of any formula in $\varsigma$.

It is well known that the calculi $\mathbf{FL}_\sigma$, except for $\mathbf{FL}_c$, satisfy the cut elimination theorem and the subformula property (see, for example, (Ono, 1998, Theorem 6, Lemma 7) or also, Ono (2003b)).

**Theorem 2** (Ono (1998))**.** *The calculi* $\mathbf{FL}_\sigma$*, with* $\sigma \neq c$*, satisfy the Cut Elimination Theorem. The calculus* $\mathbf{FL}_c$ *does not satisfy the Cut Elimination Theorem. The calculi* $\mathbf{FL}_\sigma$*, with* $\sigma \neq c$*, satisfy the Subformula Property. The calculus* $\mathbf{FL}_c$ *does not satisfy the Subformula Property.*

As a direct consequence, we have that the calculi $\mathbf{FL}_\sigma[\Psi]$,[7] $\Psi$ being a sublanguage of $\mathfrak{L}$, also satisfy the Cut Elimination Theorem and the Subformula Property, except for $\sigma = c$.

**Corollary 1.** *Let* $\Psi$ *be any sublanguage of* $\mathfrak{L}$ *whatsoever. The calculi* $\mathbf{FL}_\sigma[\Psi]$*, with* $\sigma \neq c$*, satisfy the Cut Elimination Theorem and and the Subformula Property.*

[6]In order to denote that a sequence of symbols $\sigma_1$ is a subsequence of another sequence $\sigma_2$, we will write $\sigma_1 \leq \sigma_2$.

[7]See notations described in Definition 3

**Remark 2** (About subsystems and fragments)**.** *At this point, we stress that we do not know a priori if* $\mathcal{FL}_\sigma[\Psi]$ *coincides or not with* $\Psi$-$\mathcal{FL}_\sigma$ *(i.e., the* $\Psi$*-fragment of* $\mathcal{FL}_\sigma$*). The Cut Elimination Theorem cannot be applied to prove that the subsystems referred to are fragments in the appropriate sublanguages, because this theorem and, therefore, the Subformula Property, are solely applicable to derivable sequents in such calculi. This fact also holds for the external systems* $\mathfrak{e}\mathcal{FL}_\sigma[\Psi]$ *and the* $\Psi$*-fragments of* $\mathfrak{e}\mathcal{FL}_\sigma$*. It is, however, clear that if* $\mathcal{FL}_\sigma[\Psi]$ *and* $\Psi$-$\mathcal{FL}_\sigma$ *coincide, then* $\mathfrak{e}\mathcal{FL}_\sigma[\Psi]$ *and* $\Psi$-$\mathfrak{e}\mathcal{FL}_\sigma$ *will coincide as well.*

3.2. **Mirror Images.** The notion of a *mirror image* comes from the framework of the residuated structures (see, for example, Jipsen & Tsinakis (2002)), where we have partial order $\leq$, an operation $*$ compatible with such order and two operations $\backslash$ and $/$ satisfying the Law of Residuation: for every $a, b, c$ from the universe, $a*b \leq c$ iff $b \leq a\backslash c$ iff $a \leq c/b$. The *mirror image of a formula* is the formula obtained by replacing the terms in the form $x * y$ with $y * x$ and the terms in the form $x\backslash y$ and $x/y$ with $y/x$ i $y\backslash x$, respectively. An essential fact regarding these residuated structures is that, according to the *Law of Mirror Images*, if $\varphi \leq \psi$ is satisfied in a residuated structure, then $\varphi' \leq \psi'$ will also be satisfied, where $\varphi'$ and $\psi'$ are the mirror images of $\varphi$ and $\psi$, respectively.

In this section, we introduce the notion of mirror image in the context of $\mathfrak{L}$-sequents and establish a result similar to the Law of Mirror Images applied to derivations in the calculi $\mathbf{FL}_\sigma$.

**Definition 7.** *Let* $\varphi \in Fm_{\mathfrak{L}}$*. We define a mapping* $\mu$ *from* $Fm_{\mathfrak{L}}$ *to* $Fm_{\mathfrak{L}}$ *by the following:*

$$\mu(\varphi) := \begin{cases} \varphi, & \text{if } \varphi \in Var \text{ } or \text{ } \varphi \in \{0,1\} \\ \mu(\alpha) \vee \mu(\beta), & \text{if } \varphi = \alpha \vee \beta, \\ \mu(\alpha) \wedge \mu(\beta), & \text{if } \varphi = \alpha \wedge \beta, \\ \mu(\beta) * \mu(\alpha), & \text{if } \varphi = \alpha * \beta, \\ \mu(\alpha)\backslash\mu(\beta), & \text{if } \varphi = \beta/\alpha, \\ \mu(\beta)/\mu(\alpha), & \text{if } \varphi = \alpha\backslash\beta, \\ \mu(\alpha)^{\backprime}, & \text{if } \varphi = {}'\alpha, \\ {}'\mu(\alpha), & \text{if } \varphi = \alpha^{\backprime}. \end{cases}$$

*Suppose that* $\mu(\varphi)$ *is the* mirror image *of* $\varphi$*. Let* $\Gamma$ *be a finite sequence of* $\mathfrak{L}$*-formulas. We define the* mirror image*,* $\mu(\Gamma)$*, of* $\Gamma$ *as follows:*

$$\mu(\Gamma) := \begin{cases} \mu(\varphi_{m-1}), \ldots, \mu(\varphi_0), & \text{if } \Gamma = \varphi_0, \ldots, \varphi_{m-1} \\ \emptyset, & \text{if } \Gamma = \emptyset. \end{cases}$$

*We define the* mirror image of a sequent $\Gamma \Rightarrow \Delta$ *as the sequent*

$$\mu(\Gamma \Rightarrow \Delta) := \mu(\Gamma) \Rightarrow \mu(\Delta).$$

*Note that* $\mu \circ \mu$ *is the identity function.*

In the following lemma we prove that the mirror image of a rule related to a calculus $\mathbf{FL}_\sigma$ is a derived rule derived from that calculus. This result will allow us to establish the *mirror image principle* which states that the mirror image of a derivation in $\mathbf{FL}_\sigma$ is a derivation as well and, in particular, that the mirror image of a derivable sequent in $\mathbf{FL}_\sigma$ is also derivable.

**Lemma 2.** *If $\langle \Phi, s\rangle$ is an instance of a rule of $\mathbf{FL}_\sigma$, then its mirror image, $\langle \mu[\Phi], \mu(s)\rangle$, is an instance of a rule of $\mathbf{FL}_\sigma$. More specifically:*

1) *If $\langle \Phi, s\rangle$ is an instance of a structural rule, then $\langle \mu[\Phi], \mu(s)\rangle$ is also an instance of such a rule.*
2) *If $\langle \Phi, s\rangle$ is an instance of an introduction rule for any of the connectives in $\{\vee, \wedge, *, 0, 1\}$, then $\langle \mu[\Phi], \mu(s)\rangle$ is also an instance of such a rule.*
3) *If $\langle \Phi, s\rangle$ is an instance of an introduction rule for the connective $\backslash$ (from the connective $/$), then $\langle \mu[\Phi], \mu(s)\rangle$ is an instance of the introduction rule for the connective $/$ (from the connective $\backslash$).*
4) *If $\langle \Phi, s\rangle$ is an instance of the introduction rule for the connective $\backprime$ (for the connective $'$), then $\langle \mu[\Phi], \mu(s)\rangle$ is an instance of the introduction rule for the connective $'$ (for the connective $\backprime$).*

*Proof.* This is a simple and routine proof. As an example we will consider the case concerning the rule $(\backslash \Rightarrow)$:

$$\frac{\Gamma \Rightarrow \varphi \qquad \Sigma, \psi, \Pi \Rightarrow \Delta}{\Sigma, \Gamma, \varphi\backslash\psi, \Pi \Rightarrow \Delta}$$

We have $\mu(\Gamma \Rightarrow \varphi) = \mu(\Gamma) \Rightarrow \mu(\varphi)$ and $\mu(\Sigma, \psi, \Pi \Rightarrow \Delta) = \mu(\Pi), \mu(\psi), \mu(\Sigma) \Rightarrow \mu(\Delta)$. By applying the rules $(/ \Rightarrow)$ to the sequents:

$$\mu(\Gamma) \Rightarrow \mu(\varphi) \text{ and } \mu(\Pi), \mu(\psi), \mu(\Sigma) \Rightarrow \mu(\Delta)$$

we obtain the sequent:

$$\mu(\Pi), \mu(\psi)/\mu(\varphi), \mu(\Gamma), \mu(\Sigma) \Rightarrow \mu(\Delta)$$

and this sequent is exactly $\mu(\Sigma, \Gamma, \varphi\backslash\psi, \Pi \Rightarrow \Delta)$. □

**Theorem 3** (Law of Mirror Images)**.** *Let $\Upsilon \cup \{\varsigma\} \subseteq Seq_{\mathfrak{L}}^{\omega\times\{0,1\}}$. It holds that:*

$$\text{if } \Upsilon \vdash_{\mathbf{FL}_\sigma} \varsigma, \text{ then } \mu[\Upsilon] \vdash_{\mathbf{FL}_\sigma} \mu(\varsigma).$$

*Specifically, if $\varsigma_1, \ldots, \varsigma_n$ is a proof of $\varsigma$ from $\Upsilon$ in $\mathbf{FL}_\sigma$, then $\mu(\varsigma_1), \ldots, \mu(\varsigma_n)$ is a proof of $\mu(\varsigma)$ from $\mu[\Upsilon]$ in $\mathbf{FL}_\sigma$.*
*In particular, if $\varsigma$ is derivable in $\mathbf{FL}_\sigma$, then $\mu(\varsigma)$ is derivable in $\mathbf{FL}_\sigma$ as well.*

*Proof.* By induction on $n$. If $\varsigma$ is an axiom, then $\mu(\varsigma)$ is also an axiom of $\mathbf{FL}_\sigma$, since $\mu(\varphi) \Rightarrow \mu(\varphi)$ is an instance of (*Axiom*), $\mu(0 \Rightarrow \emptyset)$ is $0 \Rightarrow \emptyset$ and $\mu(\emptyset \Rightarrow 1)$ is $\emptyset \Rightarrow 1$. If $\varsigma \in \Upsilon$, it is obviously true.

Suppose that $\varsigma$ is obtained by applying an instance $\langle\{\varsigma_i, \varsigma_j\}, \varsigma\rangle$, with $i \leq j < n$, of a rule of $\mathbf{FL}_\sigma$. Then, by Lemma 2, we have that $\langle\{\mu(\varsigma_i), \mu(\varsigma_j)\}, \mu(\varsigma)\rangle$ is an instance of a rule of $\mathbf{FL}_\sigma$. Moreover, according to the induction hypothesis $\mu(\varsigma_1), \ldots, \mu(\varsigma_i)$ is a proof of $\mu(\varsigma_i)$ and $\mu(\varsigma_1), \ldots, \mu(\varsigma_j)$ is a proof of $\mu(\varsigma_j)$. Therefore, $\mu(\varsigma_1), \ldots, \mu(\varsigma_n)$ is a proof of $\mu(\varsigma)$. □

**Remark 3.** *Note that the same result as that obtained in Theorem 3 is also valid for the calculi* $\mathbf{FL}_\sigma[\Psi]$ *obtained as a restriction on the rules of* $\mathbf{FL}_\sigma$ *to a sublanguage* $\Psi$ *of* $\mathfrak{L}$ *such that* $\langle \backslash, / \rangle \leq \Psi$ *or that* $\langle `, ' \rangle \leq \Psi$ *and, of course, for all sublanguages having neither implications nor negations. So we then have the more general result given in the following theorem.*

**Theorem 4.** *Let* $\Psi$ *be a sublanguage of* $\mathfrak{L}$ *such that* $\langle \backslash, / \rangle \leq \Psi$ *or such that* $\langle `, ' \rangle \leq \Psi$ *or such that it contains neither implications nor negations. Then, for any* $\Upsilon \cup \{\varsigma\} \subseteq Seq_{\Psi}^{\omega\times\{0,1\}}$*, it holds that:*

$$\textit{if } \Upsilon \vdash_{\mathbf{FL}_\sigma[\Psi]} \varsigma, \textit{ then } \mu[\Upsilon] \vdash_{\mathbf{FL}_\sigma[\Psi]} \mu(\varsigma).$$

*In particular, if* $\varsigma$ *is derivable in* $\mathbf{FL}_\sigma[\Psi]$*, then* $\mu(\varsigma)$ *is derivable in* $\mathbf{FL}_\sigma[\Psi]$ *as well.*

3.3. **Characterization of Leibniz Congruence.** In this section, we give a characterization of the sequential Leibniz congruence $\Omega\Phi$ associated with every $\mathcal{FL}[\Psi]$-theory, $\Phi$, for all languages $\Psi \leq \mathfrak{L}$. We also give a characterization of the Leibniz congruence $\Omega T$ associated with every $\mathfrak{c}\mathcal{FL}[\Psi]$-theory, $T$, for all languages $\Psi \leq \mathfrak{L}$ which contain the connective $\backslash$ or the connective $/$. First, we establish some results that encompass more general Gentzen systems. The characterizations presented here will be a consequence of these results.

**Theorem 5.** *Let* $\mathcal{L}$ *be any language,* $\mathcal{T} \subseteq \omega\times\omega$ *and* $\langle \mathbf{A}, R\rangle$ *an* $\langle \mathcal{L}, \mathcal{T}\rangle$*-matrix. We consider the relation:*

$$\Theta_R = \{\langle a,b\rangle \in A^2 : \langle a,b\rangle \in R \textit{ and } \langle b,a\rangle \in R\}.$$

*It holds that:*

i) *If* $\Theta_R$ *is reflexive, then* $\Omega_{\mathbf{A}}R \subseteq \Theta_R$*.*
ii) *If* $\Theta_R \in Con(\mathbf{A})$ *and* $R$ *is closed under the cut rule, then* $\Theta_R = \Omega_{\mathbf{A}}R$*.*

*Proof. i)*: Let $\langle a,b\rangle \in \Omega_{\mathbf{A}}R$. As $\Omega_{\mathbf{A}}$ is compatible with $R$ and $\langle a,a\rangle \in R$, we have $\langle a,a\rangle(a|b) \in R$ and, therefore, $\langle a,b\rangle \in R$ and $\langle b,a\rangle \in R$, i.e., $\langle a,b\rangle \in \Theta_R$.
*ii)*: Suppose $\langle \bar{x},\bar{y}\rangle \in R$ and $\langle a,b\rangle \in \Theta_R$. We assume that $\bar{x} = \bar{x}_1, a, \bar{x}_2$. So, $\langle\langle \bar{x}_1, a, \bar{x}_2\rangle, \bar{y}\rangle \in R$ but, as $\langle b,a\rangle \in R$ and $R$ is closed under the cut rule, we obtain $\langle\langle \bar{x}_1, b, \bar{x}_2\rangle, \bar{y}\rangle \in R$. Similarly, if $\bar{y} = \bar{y}_1, a, \bar{y}_2$, it is proved that $\langle \bar{x}, \langle \bar{y}_1, b, \bar{y}_2\rangle\rangle \in R$. □

**Corollary 2.** *Let* $\mathcal{G} = \langle \mathcal{L}, \mathcal{T}, \vdash_{\mathcal{G}}\rangle$ *be a Gentzen system such that* $\langle 1,1\rangle \in \mathcal{T}$*. Suppose that* $\mathcal{G}$ *satisfies the structural rules* $(Axiom)$ *and* $(Cut)$*. Given a* $\mathcal{G}$*-theory,* $\Phi$*, the set:*

$$\Theta_\Phi = \{\langle \varphi,\psi\rangle \in Fm_{\mathcal{L}}^2 : \varphi \Rightarrow \psi \in \Phi \textit{ and } \psi \Rightarrow \varphi \in \Phi\}$$

*is an equivalence relation and, if* $\Theta_\Phi \in Con(\mathbf{Fm}_{\mathcal{L}})$*, then* $\Omega\Phi = \Theta_\Phi$*.*

*Proof.* $\Theta_\Phi$ is an equivalence relation:

- *Reflexivity:* For each formula $\varphi$, we have that $\langle \varphi,\varphi\rangle \in \Theta_\Phi$, since $\varphi \Rightarrow \varphi$ is an instance of $(Axiom)$.
- *Symmetry:* According to the definition of $\Theta_\Phi$.
- *Transitivity:* If $\varphi \Rightarrow \psi \in \Phi$ and $\psi \Rightarrow \gamma \in \Phi$, then by applying $(Cut)$ we obtain $\varphi \Rightarrow \gamma \in \Phi$.

If $\Theta_\Phi \in Con(\mathbf{Fm}_\mathcal{L})$, by Theorem 5, we have that $\Omega\Phi = \Theta_\Phi$. □

**Corollary 3.** *Let $\mathcal{G}$ be as in Corollary 2. Suppose that, for every $\mathcal{G}$-theory $\Phi$, $\Theta_\Phi \in Con(\mathbf{Fm}_\mathcal{L})$. Thus, it is the case that:*

$$\text{for each } \Phi_1, \Phi_2 \in Th\,\mathcal{G}, \text{ if } \Phi_1 \subseteq \Phi_2, \text{ then } \Omega\Phi_1 \subseteq \Omega\Phi_2.$$

*Proof.* If $\Phi_1, \Phi_2 \in Th\,\mathcal{G}$ and furthermore $\Phi_1 \subseteq \Phi_2$ and also $\langle\varphi, \psi\rangle \in \Theta_{\Phi_1}$, then $\varphi \Rightarrow \psi \in \Phi_2$ and also $\psi \Rightarrow \varphi \in \Phi_2$; and, therefore, $\langle\varphi, \psi\rangle \in \Theta_{\Phi_2}$. So, $\Theta_{\Phi_1} \subseteq \Theta_{\Phi_2}$, that is, by Corollary 2, $\Omega\Phi_1 \subseteq \Omega\Phi_2$. □

**Theorem 6.** *Let $\mathcal{G} = \langle\mathcal{L}, \mathcal{T}, \vdash_\mathcal{G}\rangle$ be a Gentzen system with $\langle 0, 1\rangle, \langle 1, 1\rangle \in \mathcal{T}$. Suppose that $\mathcal{G}$ satisfies the structural rules (Axiom) and (Cut) and that $\mathcal{L}$ has a connective $\backslash$ whereby the following rules are satisfied:*

$$\frac{\Gamma \Rightarrow \varphi \qquad \Sigma, \psi, \Pi \Rightarrow \Delta}{\Sigma, \Gamma, \varphi\backslash\psi, \Pi \Rightarrow \Delta}\ (\backslash \Rightarrow) \qquad\qquad \frac{\varphi, \Gamma \Rightarrow \psi}{\Gamma \Rightarrow \varphi\backslash\psi}\ (\Rightarrow \backslash)$$

*If $T$ is a theory of the external system $\mathcal{E}(\mathcal{G})$, let $\Phi_T$ be the $\mathcal{G}$-theory generated by $\{\emptyset \Rightarrow \alpha : \alpha \in T\}$ and let $\Theta_{\Phi_T}$ be the equivalence relation of Corollary 2 defined by $\Phi_T$. It holds that:*

i) $\Theta_{\Phi_T} = \{\langle\varphi, \psi\rangle \in Fm^2_\mathcal{L} : \varphi\backslash\psi \in T \text{ and } \psi\backslash\varphi \in T\}$.

ii) *If $\Theta_{\Phi_T} \in Con(\mathbf{Fm}_\mathcal{L})$, then $\Omega T = \Omega\Phi_T$.*

*Proof. i)*: $\langle\varphi, \psi\rangle \in \Theta_{\Phi_T}$ is equivalent to $\{\varphi \Rightarrow \psi, \varphi \Rightarrow \psi\} \subseteq \Phi_T$. By using (*Axiom*), (*Cut*), $(\backslash \Rightarrow)$ i $(\Rightarrow \backslash)$, it is easy to see that:

$$\{\varphi \Rightarrow \psi, \psi \Rightarrow \varphi\} \dashv\vdash_\mathcal{G} \{\emptyset \Rightarrow \varphi\backslash\psi, \emptyset \Rightarrow \psi\backslash\varphi\}.$$

Consequently, $\langle\varphi, \psi\rangle \in \Theta_{\Phi_T}$ is equivalent to: $\{\emptyset \Rightarrow \varphi\backslash\psi, \emptyset \Rightarrow \psi\backslash\varphi\} \subseteq \Phi_T$. However, due to the definition of $\Phi_T$, this is equivalent to

$$\{\emptyset \Rightarrow \alpha : \alpha \in T\} \vdash_\mathcal{G} \{\emptyset \Rightarrow \varphi\backslash\psi, \emptyset \Rightarrow \psi\backslash\varphi\},$$

which, due to the definition of $\mathcal{E}(\mathcal{G})$, is equivalent to $T \vdash_{\mathcal{E}(\mathcal{G})} \{\varphi\backslash\psi, \psi\backslash\varphi\}$ that, as $T$ is a theory of $\mathcal{E}(\mathcal{G})$, is equivalent to $\{\varphi\backslash\psi, \psi\backslash\varphi\} \subseteq T$.
*ii)*: If $\Theta_{\Phi_T}$ is a congruence, due to Corollary 2, we have that $\Omega\Phi_T = \Theta_{\Phi_T}$. If $\langle\varphi, \psi\rangle \in \Omega\Phi_T$ and $\alpha \in T$, then $\emptyset \Rightarrow \alpha \in \Phi_T$ and, by the compatibility of $\Omega\Phi_T$ with $\Phi_T$, we have that $(\emptyset \Rightarrow \alpha)(\varphi|\psi) \in \Phi_T$. So, $\emptyset \Rightarrow \alpha(\varphi|\psi) \in \Phi_T$, that is, $\alpha(\varphi|\psi) \in T$. Therefore, $\Theta_{\Phi_T}$ is compatible with $T$. Now let $\vartheta$ be a congruence compatible with $T$ and suppose $\langle\varphi, \psi\rangle \in \vartheta$. We have that $\varphi\backslash\varphi \in T$, since $\emptyset \Rightarrow \varphi\backslash\varphi$ is obtained from $\varphi \Rightarrow \varphi$ by applying $(\Rightarrow \backslash)$. Therefore, by the compatibility of $\vartheta$ with $T$, we have $\varphi\backslash\psi \in T$ and $\psi\backslash\varphi \in T$, i.e., $\langle\varphi, \psi\rangle \in \Theta_{\Phi_T}$. In short, $\Theta_{\Phi_T}$ is the largest congruence compatible with $T$. □

**Corollary 4.** *Let $\mathcal{G}$ be as in Theorem 6. Suppose that $\Theta_{\Phi_T} \in Con(\mathbf{Fm}_\mathcal{L})$, for every $\mathcal{E}(\mathcal{G})$-theory, $T$. Thus, for each $T_1, T_2 \in Th\,\mathcal{E}(\mathcal{G})$, if $T_1 \subseteq T_2$, then $\Omega T_1 \subseteq \Omega T_2$, that is, $\mathcal{E}(\mathcal{G})$ is protoalgebraic.*

*Proof.* If $T_1 \subseteq T_2$, then clearly $\Phi_{T_1} \subseteq \Phi_{T_2}$. Therefore, by Corollary 3, we have that $\Omega\Phi_{T_1} \subseteq \Omega\Phi_{T_2}$; but this, by *ii)* of Theorem 6, means that $\Omega T_1 \subseteq \Omega T_2$. □

Similarly, we obtain the following results.

**Theorem 7.** *Let $\mathcal{G} = \langle \mathcal{L}, \mathcal{T}, \vdash_{\mathcal{G}} \rangle$ be a Gentzen system with $\langle 0,1\rangle, \langle 1,1\rangle \in \mathcal{T}$. Suppose that $\mathcal{G}$ satisfies the structural rules (Axiom) and (Cut) and that $\mathcal{L}$ contains a connective / whereby the following rules are satisfied:*

$$\frac{\Gamma \Rightarrow \varphi \quad \Sigma, \psi, \Pi \Rightarrow \Delta}{\Sigma, \psi/\varphi, \Gamma, \Pi \Rightarrow \Delta}\ (/\Rightarrow) \qquad \frac{\Gamma, \varphi \Rightarrow \psi}{\Gamma \Rightarrow \psi/\varphi}\ (\Rightarrow /)$$

*If $T$ is a theory of the external system $\mathcal{E}(\mathcal{G})$, let $\Phi_T$ be the $\mathcal{G}$-theory generated by $\{\emptyset \Rightarrow \alpha : \alpha \in T\}$ and let $\Theta_{\Phi_T}$ be the equivalence relation in Corollary 2 defined by $\Phi_T$. The following hold:*

i) $\Theta_{\Phi_T} = \{\langle \varphi, \psi\rangle \in Fm^2_{\mathcal{L}} : \psi/\varphi \in T$ *and* $\varphi/\psi \in T\}$.
ii) *If* $\Theta_{\Phi_T} \in Con(\mathbf{Fm}_{\mathcal{L}})$, *then* $\Omega T = \Omega \Phi_T$.

**Corollary 5.** *Let $\mathcal{G}$ be as in Theorem 7. Suppose that $\Theta_{\Phi_T} \in Con(\mathbf{Fm}_{\mathcal{L}})$, for every $\mathcal{E}(\mathcal{G})$-theory, $T$. So, for each $T_1, T_2 \in Th\,\mathcal{E}(\mathcal{G})$, if $T_1 \subseteq T_2$, then $\Omega T_1 \subseteq \Omega T_2$; that is, $\mathcal{E}(\mathcal{G})$ is protoalgebraic.*

With the above results, we have all we need to establish characterizations of the Leibniz congruence of $\mathcal{FL}[\Psi]$-theories for any $\Psi \leq \mathfrak{L}$, and of the theories of the external systems $\mathfrak{e}\mathcal{FL}[\Psi]$ for all the sublanguages $\Psi$ of $\mathfrak{L}$ that contain one of the implication connectives.

**Lemma 3.** *Let $\Psi$ be any sublanguage of the language $\mathfrak{L}$ of* **FL**. *For every $\mathcal{FL}[\Psi]$-theory, $\Phi$, the set:*

$$\Theta_\Phi = \{\langle \varphi, \psi\rangle \in Fm_\Psi : \varphi \Rightarrow \psi \in \Phi \text{ and } \psi \Rightarrow \varphi \in \Phi\}$$

*is a congruence of* $\mathbf{Fm}_\Psi$.

*Proof.* Suppose $\langle \varphi_1, \psi_1\rangle, \langle \varphi_2, \psi_2\rangle \in \theta_\Phi$, that is:

$$\{\varphi_1 \Rightarrow \psi_1, \psi_1 \Rightarrow \varphi_1, \varphi_2 \Rightarrow \psi_2, \psi_2 \Rightarrow \varphi_2\} \subseteq \Phi.$$

By using the introduction rules of each binary connective $\odot \in \Psi$ it is easy to see that:

$$\{\varphi_1 \odot \varphi_2 \Rightarrow \psi_1 \odot \psi_2, \psi_1 \odot \psi_2 \Rightarrow \varphi_1 \odot \varphi_2\} \subseteq \Phi.$$

In contrast, if $\langle\, ^{\backprime} \rangle \leq \Psi$ or $\langle\, ' \rangle \leq \Psi$, by using the introduction rules for the negations, it is easy to prove that if $\{\varphi \Rightarrow \psi, \psi \Rightarrow \varphi\} \subseteq \Phi$, then:

$$\{\varphi^{\backprime} \Rightarrow \psi^{\backprime}, \psi^{\backprime} \Rightarrow \varphi^{\backprime}\} \subseteq \Phi \quad \text{or} \quad \{{}'\varphi \Rightarrow {}'\psi, {}'\psi \Rightarrow {}'\varphi\} \subseteq \Phi.$$

□

**Corollary 6.** *Let $\Psi$ be any sublanguage of the language $\mathfrak{L}$ of* **FL**. *For every $\mathcal{FL}[\Psi]$-theory, $\Phi$:*

$$\Omega\Phi = \{\langle \varphi, \psi\rangle \in Fm_\Psi : \varphi \Rightarrow \psi \in \Phi \text{ and } \psi \Rightarrow \varphi \in \Phi\}.$$

*Proof.* The set $\Theta_\Phi$ of Lemma 3 is a congruence. Therefore, since the Gentzen system $\mathcal{FL}[\Psi]$ determined by $\mathbf{FL}[\Psi]$ complies with the conditions of the Gentzen systems considered in Corollary 2, we have that $\Theta_\Phi = \Omega\Phi$. □

**Corollary 7.** *Let $\Psi$ be a sublanguage of $\mathfrak{L}$ such that $\langle \backslash \rangle \leq \Psi$. Then, for every $\mathfrak{e}\mathcal{FL}[\Psi]$-theory, $T$, $\Omega T = \{\langle \varphi, \psi\rangle \in Fm^2_\Psi : \varphi\backslash\psi \in T$ and $\psi\backslash\varphi \in T\}$.*

*Proof.* Let $\Phi_T$ be the $\mathfrak{e}\mathcal{FL}[\Psi]$-theory generated by $\{\emptyset \Rightarrow \varphi : \varphi \in T\}$. System $\mathcal{FL}[\Psi]$ complies with the conditions of the Gentzen systems considered in Theorem 6. So, for $i)$ of Theorem 6, we have that:

$$\Theta_{\Phi_T} = \{\langle \varphi, \psi \rangle \in Fm^2_\Psi : \varphi \backslash \psi \in T \text{ and } \psi \backslash \varphi \in T\}.$$

Nevertheless, by Lemma 3, $\Theta_{\Phi_T}$ is a congruence and, by $ii)$ of Theorem 6, $\Omega T = \Omega \Phi_T$. □

Similarly, by using Theorem 7, we obtain the following result:

**Corollary 8.** *Let $\Psi$ be a sublanguage of $\mathfrak{L}$ such that $\langle / \rangle \leq \Psi$. Then, for every $\mathfrak{e}\mathcal{FL}[\Psi]$-theory, $T$, $\Omega T = \{\langle \varphi, \psi \rangle \in Fm^2_\Psi : \psi / \varphi \in T$ and $\varphi / \psi \in T\}$.*

Finally, we obtain that the external systems $\mathfrak{e}\mathcal{FL}[\Psi]$, for a sublanguage of $\mathfrak{L}$ that contains at least one of the two implications, are protoalgebraic.

**Corollary 9.** *If $\Psi$ is a sublanguage of $\mathfrak{L}$ such that $\langle \backslash \rangle \leq \Psi$ or such that $\langle / \rangle \leq \Psi$, then $\mathfrak{e}\mathcal{FL}[\Psi]$ is protoalgebraic.*

*Proof.* By Corollaries 4 and 5 □

## 4. Hilbert-style Axiomatization for $\mathfrak{e}\mathcal{FL}_\sigma$ Systems.

In this section we present some axiomatizations for the external systems, $\mathfrak{e}\mathcal{FL}_\sigma$, that are found in the literature.

**Definition 8.** (Galatos et al., 2007a, Section 2.5.1) **HFL** *is the Hilbert system with language $\langle \vee, \wedge, *, \backslash, /, 0, 1 \rangle$ of type $\langle 2, 2, 2, 2, 2, 0, 0 \rangle$ defined by the following axioms and rules:*

| | | |
|---|---|---|
| $(\backslash\text{-}id)$ | $\varphi \backslash \varphi$ | *(identity)* |
| $(\backslash\text{-}pf)$ | $(\varphi \backslash \psi) \backslash [(\gamma \backslash \varphi) \backslash (\gamma \backslash \varphi)]$ | *(prefixing)* |
| $(\backslash\text{-}as)$ | $\varphi \backslash [(\psi / \varphi) \backslash \psi]$ | *(assertion)* |
| $(a)$ | $[(\psi \backslash \gamma) / \varphi] \backslash [\psi \backslash (\gamma / \varphi)]$ | *(associativity)* |
| $(* \backslash /)$ | $[(\psi * (\psi \backslash \varphi)) / \psi] \backslash (\varphi / \psi)$ | *(fusion implications)* |
| $(* \wedge)$ | $[(\varphi \wedge 1) * (\psi \wedge 1)] \backslash (\varphi \wedge \psi)$ | *(fusion conjunction)* |
| $(\wedge_1\backslash)$ | $(\varphi \wedge \psi)\backslash\varphi$ | *(conjunction implication 1)* |
| $(\wedge_2\backslash)$ | $(\varphi \wedge \psi)\backslash\psi$ | *(conjunction implication 2)* |
| $(\backslash\wedge)$ | $[(\gamma\backslash\varphi) \wedge (\gamma\backslash\psi)] \backslash [\gamma\backslash(\varphi \wedge \psi)]$ | *(implication conjunction)* |
| $(\backslash\vee_1)$ | $\varphi\backslash(\varphi \vee \psi)$ | *(implication disjunction 1)* |
| $(\backslash\vee_2)$ | $\psi\backslash(\varphi \vee \psi)$ | *(implication disjunction 2)* |
| $(\vee\backslash)$ | $[(\varphi\backslash\gamma) \wedge (\psi\backslash\gamma)] \backslash [(\varphi \vee \psi)\backslash\gamma]$ | *(disjunction implication)* |
| $(\backslash *)$ | $\psi\backslash(\varphi\backslash(\varphi * \psi))$ | *(implication fusion)* |
| $(*\backslash)$ | $[\psi\backslash(\varphi\backslash\gamma)] \backslash [(\varphi * \psi)\backslash\gamma]$ | *(fusion implication)* |
| $(1)$ | $1$ | *(unit)* |
| $(1\backslash)$ | $1\backslash(\varphi\backslash\varphi)$ | *(unit implication)* |
| $(\backslash 1)$ | $\varphi\backslash(1\backslash\varphi)$ | *(implication unit)* |
| $(\backslash\text{-}mp)$ | $\langle\{\varphi, \varphi\backslash\psi\}, \psi\rangle$ | *($\backslash$-modus ponens)* |
| $(adj_u)$ | $\langle\{\varphi\}, \varphi \wedge 1\rangle$ | *(adjunction unit)* |
| $(\backslash\text{-}pn)$ | $\langle\{\varphi\}, \psi\backslash(\varphi * \psi)\rangle$ | *($\backslash$-product normality)* |
| $(/\text{-}pn)$ | $\langle\{\varphi\}, (\psi * \varphi)/\psi)\rangle$ | *($/$-product normality)* |

It is easy to see that, if we add to the previous axiomatization the schemata:

| | |
|---|---|
| $(\backslash\text{-}\, def_1)$ | $\varphi^{\backprime}\backslash(\varphi\backslash 0)$ |
| $(\backslash\text{-}\, def_2)$ | $(\varphi\backslash 0)\backslash\varphi^{\backprime}$ |
| $(/\text{-}\, def_1)$ | ${}^{\prime}\varphi/(0/\varphi)$ |
| $(/\text{-}\, def_2)$ | $(0/\varphi)/{}^{\prime}\varphi$ |

then the Hilbert system which defines the new axiomatization is a definitional expansion, in the language $\langle \vee, \wedge, *, \backslash, /, {}^{\backprime}, {}^{\prime}, 0, 1\rangle$, of the system **HFL**.

**Definition 9.** (Galatos et al., 2007a, Section 2.5.1) $\mathbf{HFL}_e$ *is the Hilbert system in the language* $\langle \vee, \wedge, *, \rightarrow, 0, 1\rangle$ *of type* $\langle 2, 2, 2, 2, 0, 0\rangle$ *defined by the following axioms and rules:*[8]

[8]We adopt the convention that the connective $\rightarrow$ establishes a less strong binding than the other connectives.

| | | |
|---|---|---|
| $(id)$ | $\varphi \to \varphi$ | *(identity)* |
| $(pf)$ | $(\varphi \to \psi) \to ((\gamma \to \varphi) \to (\gamma \to \varphi))$ | *(prefixing)* |
| $(per)$ | $(\varphi \to (\psi \to \gamma)) \to (\psi \to (\varphi \to \gamma))$ | *(permutation)* |
| $(* \wedge)$ | $((\varphi \wedge 1) * (\psi \wedge 1)) \to \varphi \wedge \psi$ | *(fusion conjunction)* |
| $(\wedge_1 \to)$ | $\varphi \wedge \psi \to \varphi$ | *(conjunction implication 1)* |
| $(\wedge_2 \to)$ | $\varphi \wedge \psi \to \psi$ | *(conjunction implication 2)* |
| $(\to \wedge)$ | $(\gamma \to \varphi) \wedge (\gamma \to \psi) \to (\gamma \to \varphi \wedge \psi)$ | *(implication conjunction)* |
| $(\to \vee_1)$ | $\varphi \to \varphi \vee \psi$ | *(implication disjunction 1)* |
| $(\to \vee_2)$ | $\psi \to \varphi \vee \psi$ | *(implication disjunction 2)* |
| $(\vee \to)$ | $(\varphi \to \gamma) \wedge (\psi \to \gamma) \to (\varphi \vee \psi \to \gamma)$ | *(disjunction implication)* |
| $(\to *)$ | $\psi \to (\varphi \to \varphi * \psi)$ | *(implication fusion)* |
| $(* \to)$ | $(\psi \to (\varphi \to \gamma)) \to (\varphi * \psi \to \gamma)$ | *(fusion implication)* |
| $(1)$ | $1$ | *(unit)* |
| $(1 \to)$ | $1 \to (\varphi \to \varphi)$ | *(unit implication)* |
| $(mp)$ | $\langle \{\varphi, \varphi \to \psi\}, \psi \rangle$ | *(modus ponens)* |
| $(adj_u)$ | $\langle \{\varphi\}, \varphi \wedge 1 \rangle$ | *(adjunction unit)* |

If we add to the axiomatization of $\mathbf{HFL}_e$ the schemata:

| | |
|---|---|
| $(\neg\text{-}def_1)$ | $\neg\varphi \to (\varphi \to 0)$ |
| $(\neg\text{-}def_2)$ | $(\varphi \to 0) \to \neg\varphi$ |

then the Hilbert system which defines the new axiomatization is a definitional expansion, in language $\langle \vee, \wedge, *, \to, \neg, 0, 1 \rangle$, of the system $\mathbf{HFL}_e$.

**Definition 10** (Systems $\mathbf{HFL}_\sigma$). *We consider the schemata:*

| | |
|---|---|
| $(w_l)$ | $\varphi \backslash (\psi \backslash \varphi)$ |
| $(w_r)$ | $0 \backslash \varphi$ |
| $(c)$ | $(\varphi \backslash (\varphi \backslash \psi)) \backslash (\varphi \backslash \psi)$ |

*If $\sigma \leq w_l w_r c$, then $\mathbf{HFL}_\sigma$ denotes the extension of $\mathbf{HFL}$ with the schemata codified by $\sigma$. If $e \leq \sigma \leq e w_l w_r c$, then $\mathbf{HFL}_\sigma$ denotes the extension of $\mathbf{HFL}$ with the schemata codified by $\sigma$. In this case, we substitute the connective $\backslash$ for the connective $\to$ in the schemata referred to.*

**Theorem 8.** (Cf. (Galatos et al., 2007a, Section 2.5)) *Let $\sigma$ be a subsequence, maybe empty, of $e w_l w_r c$; then the following holds:*

$$\mathfrak{e}\mathcal{FL}_\sigma = \mathbf{HFL}_\sigma.$$

4.1. **Strongly Separable Axiomatizations.** In van Alten & Raftery (2004) a strongly separable axiomatization for the fragment without exponents, without 0 and without additive constants of *intuitionist linear logic* (see (Troelstra, 1992, p.67) ) is presented; that is, for the logic of commutative residuated lattices. The presentation, when we add the symbol 0 to the language without adding any schema where the symbol expressly appears, yields a strongly separable axiomatization with respect to $\to$ of the Hilbert system $\mathbf{HFL}_e$ (see first paragraph of (Galatos et al., 2007a, Section 2.5) ). We transcribe it below:

| | |
|---|---|
| $(id)$ | $\varphi \to \varphi$ |
| $(pf)$ | $(\varphi \to \psi) \to [(\gamma \to \varphi) \to (\gamma \to \varphi)]$ |
| $(per)$ | $[\varphi \to (\psi \to \gamma)] \to [\psi \to (\varphi \to \gamma)]$ |
| $(mp)$ | $\langle\{\varphi, \varphi \to \psi\}, \psi\rangle$ |
| $(\to \vee_1)$ | $\varphi \to \varphi \vee \psi$ |
| $(\to \vee_2)$ | $\psi \to \varphi \vee \psi$ |
| $(dis)$ | $\langle\{\varphi \to \gamma, \psi \to \gamma\}, \varphi \vee \psi \to \gamma\rangle$ |
| $(\wedge_1 \to)$ | $\varphi \wedge \psi \to \varphi$ |
| $(\wedge_2 \to)$ | $\varphi \wedge \psi \to \psi$ |
| $(\to \wedge)$ | $(\gamma \to \varphi) \wedge (\gamma \to \psi) \to (\gamma \to \varphi \wedge \psi)$ |
| $(adj)$ | $\langle\{\varphi, \psi\}, \varphi \wedge \psi\rangle$ |
| $(\to *)$ | $(\psi \to (\varphi \to \gamma)) \to (\varphi * \psi \to \gamma)$ |
| $(* \to)$ | $\psi \to (\varphi \to (\varphi * \psi))$ |
| $(1)$ | $1$ |
| $(1 \to)$ | $1 \to (\varphi \to \varphi)$ |

There is a strongly separable axiomatization for the logic **HFL** with respect to the set of basic connectives $\{\backslash, /\}$ (see Galatos & Ono (2010)).

**Remark 4.** *The system* $\mathbf{HFL}_{ew}$ *is definitionally equivalent to Monoidal Logic (see Höhle (1995); Gottwald (2001); Gottwald et al. (2003)), also called Intuitionistic Logic without contraction in Adillon (2001); Adillon & Verdú (2000); Bou et al. (2006). If we add to the calculus of Van Alten and Raftery the schema* $\varphi \to (\psi \to \varphi)$*, then we obtain an axiomatization for* $\mathbf{HFL}_{ew_l}$ *where one can prove that the rule* $(dis)$ *may be replaced by the schema* $(\vee \to)$*, the rule* $(adj)$*, by the schema* $\varphi \to (\psi \to \varphi \wedge \psi)$*, and schemata* $(1)$ *and* $(1 \to)$*, by the schema* $\varphi \to 1$*. If we add schema* $0 \to \varphi$ *to this axiomatization for* $\mathbf{HFL}_{ew_l}$*, we obtain an axiomatization for Monoidal Logic in language* $\langle\vee, \wedge, *, \to, 0, 1\rangle$*. Ono & Komori (1985) show that this axiomatization is separable. It can be seen that it is also strongly separable.*

## Part 2. THE ALGEBRAS

In this part we present the ordered, latticed and semilatticed algebraic structures that will constitute the semantic core of some of the fragments considered in Part 3. In Section 5, we give some preliminary notions and results. In Section 6.1, we introduce the notion of *pointed monoid*. A pointed monoid (ordered, semilatticed or latticed) is obtained by adding the constant symbol 0 to the type of similarity of a monoid (ordered, semilatticed or latticed): the symbol is interpreted as a fixed, but arbitrary, element of the universe of the structure. We define the varieties of algebras $\mathring{\mathbb{M}}^{s\ell}_{\sigma}$ and $\mathring{\mathbb{M}}^{\ell}_{\sigma}$, where the subindex $\sigma$ is a subsequence of $ew_lw_rc$ and the symbols $e$, $w_l$, $w_r$ and $c$ codify what we refer to as (algebraic) exchange, right-weakening, left-weakening and contraction properties, respectively. These properties, which are expressed by quasi-inequations are equivalent, respectively, to the following properties: commutativity, integrality, 0-boundedness and increasing idempotency. As we see in later parts of the paper, once a sequence $\sigma$ is fixed, the $\mathring{\mathbb{M}}^{s\ell}_{\sigma}$ and $\mathring{\mathbb{M}}^{\ell}_{\sigma}$ classes are equivalent, respectively, to the $\langle \vee, *, 0, 1\rangle$-fragment and to the $\langle \vee, \wedge, *, 0, 1\rangle$-fragment of the Gentzen system $\mathcal{FL}_{\sigma}$, and they are algebraic semantics, respectively, of the $\langle \vee, *, 0, 1\rangle$-fragment and of the $\langle \vee, \wedge, *, 0, 1\rangle$-fragment of the external system $\mathfrak{e}\mathcal{FL}_{\sigma}$.

In Section 7 we recap the notion of residuation and both the definitions and properties of residuated lattices together with those of $\mathbb{FL}$-algebras and of their subvarieties, $\mathbb{FL}_{\sigma}$. In our presentation we use the notion of *relative pseudocomplement.*

In Section 8 we present the notion of *pseudocomplementation* in the framework of semilatticed monoids and we define the classes $\mathbb{PM}^{s\ell}$ and $\mathbb{PM}^{\ell}$ of semilatticed and latticed pseudocomplemented monoids. We show that the $\mathbb{PM}^{s\ell}$ and $\mathbb{PM}^{\ell}$ classes are varieties. We also analyze the case in which the pseudocomplementation is with respect to the minimum element of the monoid. The pseudocomplements constitute the algebraic counterpart of negations: in later parts of the paper we will state the connection between the varieties $\mathbb{PM}^{s\ell}_{\sigma}$ and $\mathbb{PM}^{\ell}_{\sigma}$ on the one hand (subvarieties of $\mathbb{PM}^{s\ell}$ and $\mathbb{PM}^{\ell}$ defined by the equations codified by $\sigma$) and on the other, the fragments of the Gentzen system $\mathcal{FL}_{\sigma}$ and the associated external system $\mathfrak{e}\mathcal{FL}_{\sigma}$ in the languages $\langle \vee, *, {}^{\backprime}, {}', 0, 1\rangle$ and $\langle \vee, \wedge, *, {}^{\backprime}, {}', 0, 1\rangle$.

## 5. Some Preliminaries.

An *algebraic language* is a pair $\mathcal{L} = \langle F, \tau\rangle$, where $F$ is a set of *functional symbols* and $\tau$ is a mapping $\tau : F \to \omega$ (where $\omega$ denotes the set of natural numbers) which is called the *algebraic similarity type.* For every $f \in F$, $\tau(f)$ is called the *arity* of the functional symbol $f$. Functional symbols with an arity of 0 are also called *constant symbols.* If $\mathcal{L} = \langle F, \tau\rangle$ is an algebraic language with a finite number of functionals, we say that $\mathcal{L}$ is *finite.* In this case, if $F = \{f_1, \dots, f_n\}$, we identify $\mathcal{L}$ with the sequence $\langle f_1, \dots, f_n\rangle$ and say that $\mathcal{L}$ is the language $\langle f_1, \dots, f_n\rangle$ of type $\langle \tau(f_1), \dots, \tau(f_n)\rangle$.

If $\mathcal{L} = \langle F, \tau\rangle$ and $\mathcal{L}' = \langle F', \tau'\rangle$ are two algebraic languages such that $F \subseteq F'$ and $\tau = \tau' \restriction F$ (i.e., the restriction of $\tau'$ to $F$), then we say that $\mathcal{L}$ is a *sublanguage* of $\mathcal{L}'$ and we write $\mathcal{L} \leq \mathcal{L}'$. In this case we also say that $\mathcal{L}'$ is an *expansion* of $\mathcal{L}$. If $\mathcal{L}'$ is an expansion of $\mathcal{L}$ with a finite set of new functionals $f_1, \dots, f_n$, sometimes we will denote the language $\mathcal{L}'$ as $\langle \mathcal{L}, f_1, \dots, f_n\rangle$ .

Given an algebraic language $\mathcal{L} = \langle F, \tau\rangle$, an algebra $\mathbf{A}$ of type $\mathcal{L}$ ($\mathcal{L}$-algebra, for short) is a pair $\langle A, \{f^{\mathbf{A}} : f \in F\}\rangle$, where $A$ is a non-empty set called the *universe* of $\mathbf{A}$ and where, for every $f \in F$: if $\tau(f) = n$, then $f^{\mathbf{A}}$ is an $n$-ary operation on $A$ (a 0-ary operation on $A$ is an element of $A$). If $\mathcal{L}$ is a finite language $\langle f_1, \dots, f_n\rangle$ of type $\langle \tau(f_1), \dots, \tau(f_n)\rangle$, we write $\mathbf{A} = \langle A, f_1^{\mathbf{A}}, \dots, f_n^{\mathbf{A}}\rangle$ and we say that $\mathbf{A}$ is an algebra of type $\langle \tau(f_1), \dots, \tau(f_n)\rangle$. The superscripts in the operations will be omitted when they are clear from the context. An algebra, $\mathbf{A}$, is *finite* if $A$ is finite and it is *trivial* if $A$ has only one element. To denote classes of algebras we will use capital letters in blackboard boldface, e.g., $\mathbb{K}, \mathbb{M} \ldots$. The members of a class $\mathbb{K}$ of algebras will sometimes be called $\mathbb{K}$-algebras.

We will call any pair $\mathcal{A} = \langle \mathbf{A}, \leq \rangle$, where $\mathbf{A}$ is an algebra and $\leq$ is a partial order defined in the universe of $\mathbf{A}$, an *algebra with order* or *order-algebra*. We say that $\mathbf{A}$ is the *algebraic reduct* of $\mathcal{A}$. If $\mathcal{L} = \langle F, \tau \rangle$ is finite with $F = \{f_1, \ldots, f_n\}$, we say that $\mathcal{A}$ is an $\langle f_1, \ldots, f_n, \leq \rangle$-algebra of algebraic type $\langle \tau(f_1), \ldots, \tau(f_n) \rangle$. As in the case of algebras, we will use capital letters in blackboard boldface, e.g., $\mathbb{K}, \mathbb{M} \ldots$ to denote the generic classes of order-algebras.

If $\mathcal{L} = \langle F, \tau \rangle$, we consider the first-order language with equality

$$\mathcal{L}^{\preccurlyeq} = \langle F \cup \{\approx, \preccurlyeq\}, \tau' \rangle,$$

where its type of similarity, $\tau'$, is given by $\tau' \restriction F = \tau$ and $\tau'(\preccurlyeq) = 2$. This language, then, besides the equality symbol, contains a proper single binary relational symbol of arity 2, $\preccurlyeq$, and a set of functional symbols. An $\mathcal{L}$-algebra with order can be seen as an $\mathcal{L}^{\preccurlyeq}$-structure, where the functional interpretation is $f^{\mathcal{A}} = f^{\mathbf{A}}$, for every $f \in F$, and the relational interpretation is $\preccurlyeq^{\mathcal{A}} = \leq$. The atomic formulas of this language are of the form $\varphi \approx \psi$ and $\varphi \preccurlyeq \psi$, where $\varphi$ and $\psi$ are $\mathcal{L}$-terms; and the rest of the $\mathcal{L}^{\preccurlyeq}$-formulas, as in the case of the $\mathcal{L}$-algebras, are generated as usual with the help of the quantifiers $\forall$ and $\exists$ together with the boolean connectives $\sqcup, \sqcap, \supset, \sim$.

We will call an $\mathcal{L}$-*inequation* any $\mathcal{L}^{\preccurlyeq}$-atomic formula of the form $\varphi \preccurlyeq \psi$ and will call an $\mathcal{L}$-*quasi-inequation* any $\mathcal{L}^{\preccurlyeq}$-formula of the form:

$$\varphi_0 \preccurlyeq \psi_0 \ \sqcap \ \ldots \ \sqcap \ \varphi_{n-1} \preccurlyeq \psi_{n-1} \supset \varphi_n \preccurlyeq \psi_n,$$

where $\varphi_i, \psi_i \in Fm_{\mathcal{L}}$ for every $i \leq n$. Note that every inequation is a quasi-inequation with $n = 0$.

A *Horn formula* is a first-order formula of one of these three types:

*i*) $\beta$,
*ii*) $\alpha_1 \sqcap \ldots \sqcap \alpha_n \supset \beta$,
*iii*) $\sim (\alpha_1 \sqcup \ldots \sqcup \alpha_n)$,

where $\alpha_1, \ldots, \alpha_n, \beta$ are atomic formulas. A *universal Horn sentence* is a universal closure of a Horn formula and we classify it as *strict* if it is the universal closure of a Horn formula of the form *i*) or *ii*). So, the $\mathcal{L}$-equations and $\mathcal{L}$-quasi-equations as well as the $\mathcal{L}$-inequations and $\mathcal{L}$-quasi-inequations are Horn formulas of the language $\mathcal{L}^{\preccurlyeq}$. We will say that an $\mathcal{L}$-quasi-inequation, $\delta$, *is satisfied*, or *is valid*, in an $\langle \mathcal{L}, \leq \rangle$-algebra $\mathcal{A}$ if the corresponding strict universal Horn sentence is valid in $\mathcal{A}$ in the usual sense; and we will annotate $\mathcal{A} \models \delta$ and $\mathbb{K} \models \delta$ if it is valid for all the members belonging to the $\mathbb{K}$ class of algebras with order.

Given an $\langle \mathcal{L}, \leq \rangle$-algebra $\mathcal{A} = \langle A, \leq \rangle$, the order is *definable* in terms of the functionals if there is an $\mathcal{L}$-formula, $\delta$, of first order with its variables in $\{x_1, x_2\}$ such that:

$$\{\langle a, b \rangle \in A^2 : a \leq b\} = \{\langle a, b \rangle \in A^2 : \mathbf{A} \models \delta^{\mathbf{A}}(a, b)\}.$$

5.1. **Monotonicity, antimonotonicity.** Let $\iota$ be an operation of arity $k \geq 1$ defined in an ordered set $\langle A, \leq \rangle$ and let $i$ be such that $1 \leq i \leq k$. We will say that $\iota$ is *monotonous in the i-th argument* with respect to the order if, for every $a, b, c_1, \ldots, c_k \in A$:

$$\iota(c_1, \ldots, c_{i-1}, a, c_{i+1}, \ldots, c_k) \leq \iota(c_1, \ldots, c_{i-1}, b, c_{i+1}, \ldots, c_k) \text{ when } a \leq b.$$

If the operation $\iota$ is monotonous in all the arguments, we will say, simply, that it is *monotonous*. In that case, we also say that the operation is *compatible* with the order. With this expression we indicate that a structure containing such an operation is not only a structure with operations and an order, but there is a nexus established between one of the operations and the order.

We will say that $\iota$ is *antimonotonous in the i-th* argument with respect to the order if, for every $a, b, c_1, \ldots, c_k \in A$:

$$\iota(c_1, \ldots, c_{i-1}, b, c_{i+1}, \ldots, c_k) \leq \iota(c_1, \ldots, c_{i-1}, a, c_{i+1}, \ldots, c_k) \text{ when } a \leq b.$$

If the operation $\iota$ is antimonotonous in all the arguments, we will say, simply, that it is *antimonotonous.* In that case we say that this operation is *dually compatible* with the order.

## 6. Pointed Semilatticed and Latticed Monoids.

We define below some classes of algebras with order that are characterized by the fact that they contain a binary operation, $*$, which is monotonous in the two arguments with respect to the order. The contents of this section are based on Dubreil-Jacotin et al. (1953) and Birkhoff (1973).

A *partially-ordered groupoid*, or a *po-groupoid* in short, is an structure $\mathcal{A} = \langle A, *, \leq \rangle$, where $\leq$ is a partial order defined in $A$ and $*$ is a binary operation on $A$ which is monotonous with respect to that order; i.e., for every $a, b, c \in A$:

$r)$ if $a \leq b,$ then $a * c \leq b * c,$
$l)$ if $a \leq b,$ then $c * a \leq c * b$

or, equivalently, if $a \leq c$ and $b \leq d,$ then $a * b \leq c * d.$

Let $\langle *, 1 \rangle$ be an algebraic language of type $\langle 2, 0 \rangle$. An algebra with order $\mathcal{A} = \langle A, *, 1, \leq \rangle$ is a *po-monoid* if $\langle A, *, \leq \rangle$ is a po-groupoid and the algebraic reduct is a monoid; that is, the operation $*$ is associative, and 1 is the identity element.[9]

A *commutative po-groupoid (po-monoid)* is a po-groupoid (po-monoid) such that the operation $*$ is commutative. In this case, conditions $r)$ and $l)$ of monotonicity are equivalent.

A *semilatticed groupoid*, or *sℓ-groupoid*, is an algebra $\mathbf{A} = \langle A, \vee, * \rangle$ of type $\langle 2, 2 \rangle$ where the following equations are satisfied:

1. A set of equations defining the $\vee$-semilattices (e.g., idempotency, associativity and commutativity of the operation $\vee$).
2. The distributivity laws of $*$ with respect to $\vee$, i.e.,
   *r)* $(x \vee y) * z \approx (x * z) \vee (y * z),$

[9]The *identity* element or the *unit* element of a groupoid $\mathbf{G} = \langle G, * \rangle$ is an element $1 \in G$ such that, for every $a \in G$, the following condition is met: $a * 1 = 1 * a = 1$. If there is such an element, then it is unique.

*l)* $z * (x \vee y) \approx (z * x) \vee (z * x)$.

A *latticed groupoid*, or *$\ell$-groupoid*, is an algebra $\mathbf{A} = \langle A, \vee, \wedge, * \rangle$ of type $\langle 2, 2, 2 \rangle$ such that in $\mathbf{A}$ a set of equations is satisfied defining the $\langle \vee, \wedge \rangle$-lattices and the equations *r)* and *l)*. A *semilatticed monoid*, or *$s\ell$-monoid*, is an algebra $\mathbf{A} = \langle A, \vee, *, 1 \rangle$ of type $\langle 2, 2, 0 \rangle$ such that $\langle A, \vee, * \rangle$ is an $s\ell$-groupoid and such that its $\langle *, 1 \rangle$-reduct is a monoid. A *latticed monoid*, or *$\ell$-monoid* is an algebra $\mathbf{A} = \langle A, \vee, \wedge, *, 1 \rangle$ of type $\langle 2, 2, 2, 0 \rangle$ such that $\langle A, \vee, \wedge, * \rangle$ is an $\ell$-groupoid and such that its $\langle *, 1 \rangle$-reduct is a monoid. It is easy to prove that equations *r)* and *l)* of the definition of semilatticed groupoids are equivalent to the equation:

$$(x \vee y) * (z \vee t) \approx (x * z) \vee (x * t) \vee (y * z) \vee (y * t). \tag{2}$$

If the operation $*$ is commutative, then the term *commutative* is added to the name of the corresponding class of algebras. Of course, in the commutative classes, the conditions *r)* and *l)* are equivalent. Every semilatticed groupoid $\mathbf{A} = \langle A, \vee, * \rangle$ defines a partially-ordered algebra $\mathcal{A} = \langle \mathbf{A}, \leq \rangle$, where $\leq$ is the order defined by the semilattice by means of the atomic formula $x \vee y \approx y$, in such a way that the structure $\langle A, *, \leq \rangle$ is a po-groupoid, as shown in the following proposition.

**Proposition 1.** *In every $s\ell$-groupoid $\mathbf{A}$ the monotonicity conditions are satisfied.*

*Proof.* Let $a, b, c, d \in A$ and suppose that $a \leq c$ and $b \leq d$. Then, by applying (2) and the properties of semilattices, we have:
$c * d = (a \vee c) * (b \vee d) = (a * b) \vee (a * d) \vee (c * b) \vee (c * d) = (a * b) \vee (c * b) \vee (a * d) \vee (c * d) = (a * b) \vee (c * b) \vee ((a \vee c) * d) = (a * b) \vee (c * b) \vee (c * d) = (a * b) \vee (c * (b \vee d)) = (a * b) \vee (c * d)$ and, therefore, $a * b \leq c * d$. $\square$

Thus, the order-algebra $\langle A, *, \leq \rangle$, where $\leq$ is the order of the semilattice (or of the lattice), is a po-groupoid. So, every $s\ell$-groupoid naturally defines a po-groupoid. It is easy to see that $\mathbf{A} = \langle A, \vee, * \rangle$, $\mathcal{A} = \langle \mathbf{A}, \leq \rangle = \langle A, \vee, *, \leq \rangle$ and $\langle A, *, \leq \rangle$ are definitionally equivalent structures. However, it is not true that every po-groupoid $\mathcal{A} = \langle A, *, \leq \rangle$, which is a semilattice under its relation of partial order, defines a semilatticed groupoid; that is, it is not true that an $s\ell$-groupoid is a po-groupoid that is a semilattice under its relation of partial order. In the following proposition it is shown that distributivity is a stronger condition than monotonicity.

**Proposition 2.** *There are po-monoids (and hence, po-groupoids) $\langle A, *, 1, \leq \rangle$ that are lattices (and hence, semilattices) under their relation of partial order which do not satisfy the distributivity condition of the operation $*$ with respect to the operation $\vee$ associated with the semilattice.*

*Proof.* Let $\mathcal{A} = \langle A, *, \leq \rangle$, where $A = \{0, a, b, 1\}$ and where the order and the operation $*$ are defined in the following diagram and table:

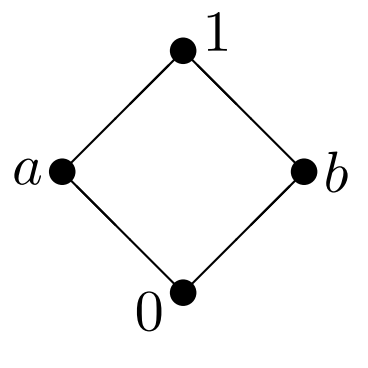

| $*$ | 0 | $a$ | $b$ | 1 |
|---|---|---|---|---|
| 0 | 0 | 0 | 0 | 0 |
| $a$ | 0 | 0 | 0 | $a$ |
| $b$ | 0 | 0 | 0 | $b$ |
| 1 | 0 | $a$ | $b$ | 1 |

Note that the operation $*$ is associative, commutative and that 1 is the unit element. The order is latticed and the operation $*$ is obviously monotonous with respect to the order. In contrast, the distributivity of $*$ with respect to $\vee$ is not satisfied:

$$(a \vee b) * b = 1 * b = b \neq 0 = 0 \vee 0 = (a * b) \vee (b * b).$$

So we have a po-monoid that is a lattice under its relation of partial order but does not satisfy the distributivity condition. □

Observe that all the binary operations of the semilatticed and latticed groupoids and monoids are monotonous. Thus, for instance, if we consider the structure $\langle A, \vee, \wedge, *, \leq \rangle$ associated with a latticed groupoid $\langle A, \vee, \wedge, * \rangle$, we have that their reducts $\langle A, \vee, \leq \rangle$, $\langle A, \wedge, \leq \rangle$ and $\langle A, *, \leq \rangle$ are po-groupoids.

6.1. **The $\mathring{\mathbb{M}}^{s\ell}_{\sigma}$ and $\mathring{\mathbb{M}}^{\ell}_{\sigma}$ varieties.** If we expand the algebraic type of the structures considered in the previous section with the constant symbol 0, which will be interpreted as a fixed element but arbitrary in the universe, we obtain the corresponding classes of structures called *pointed* structures.

- A *pointed po-groupoid* is a partially-ordered algebra $\mathcal{A} = \langle A, *, 0, \leq \rangle$ of algebraic type $\langle 2, 0 \rangle$, where $\langle A, *, \leq \rangle$ is a po-groupoid and 0 is a fixed but arbitrary element of $A$ (a *distinguished* element).
- A *pointed po-monoid* is a partially-ordered algebra $\mathcal{A} = \langle A, *, 0, 1, \leq \rangle$ of algebraic type $\langle 2, 0, 0 \rangle$ such that $\langle A, *, 1, \leq \rangle$ is a po-monoid and 0 is a fixed but arbitrary element of $A$.
- A *pointed $s\ell$-monoid* is an algebra $\mathbf{A} = \langle A, \vee, *, 0, 1 \rangle$ of type $\langle 2, 2, 0, 0 \rangle$ such that $\mathbf{A} = \langle A, \vee, *, 1 \rangle$ is an $s\ell$-monoid and 0 is a fixed but arbitrary element of $A$. We will denote the class of $s\ell$-pointed monoids by $\mathring{\mathbb{M}}^{s\ell}$.
- A *pointed $\ell$-monoid* is an algebra $\mathbf{A} = \langle A, \vee, \wedge, *, 0, 1 \rangle$ of type $\langle 2, 2, 2, 0, 0 \rangle$ such that $\mathbf{A} = \langle A, \vee, \wedge, *, 1 \rangle$ is an $\ell$-monoid and 0 is a fixed but arbitrary element of $A$. We will denote the class of $\ell$-pointed monoids by $\mathring{\mathbb{M}}^{\ell}$.

Obviously, the classes $\mathring{\mathbb{M}}^{s\ell}$ and $\mathring{\mathbb{M}}^{\ell}$ are varieties. This is summarized in the following propositions.

**Proposition 3.** $\mathring{\mathbb{M}}^{s\ell}$ *is the equational class of* $\langle \vee, *, 0, 1 \rangle$*-algebras of type* $\langle 2, 2, 0, 0 \rangle$ *that satisfy: a) a set of equations that define the class of* $\vee$*-semilattices; b) a set of* $\langle *, 1 \rangle$*-equations that define the class of monoids; and c) the distributivity equations of the monoidal operation with respect to the operation* $\vee$.

**Proposition 4.** $\mathring{\mathbb{M}}^{\ell}$ *is the equational class of the* $\langle \vee, \wedge, *, 0, 1 \rangle$*-algebras of type* $\langle 2, 2, 2, 0, 0 \rangle$ *that satisfy: a) a set of* $\langle \vee, \wedge \rangle$*-equations that define the class of the lattices; b) a set of* $\langle *, 1 \rangle$*-equations that define the class of the monoids; and c) the distributivity equations of the monoidal operation with respect to the operation* $\vee$.

The varieties $\mathring{\mathbb{M}}^{s\ell}$ and $\mathring{\mathbb{M}}^{\ell}$, as we will see in Part 3, constitute the algebraic counterpart of the $\langle \vee, *, 0, 1 \rangle$-fragment and the $\langle \vee, \wedge, *, 0, 1 \rangle$-fragment of the Gentzen system $\mathcal{FL}$, respectively. We define below the varieties $\mathring{\mathbb{M}}^{s\ell}_{\sigma}$ and $\mathring{\mathbb{M}}^{\ell}_{\sigma}$, where $\sigma$ is a subsequence of the sequence $ew_lw_rc$, and where the symbols $e$, $w_l$, $w_r$ and $c$ codify the (algebraic) properties that we will call *exchange*, *right-weakening*, *left-weakening* and *contraction*, respectively. As we show below, these properties have the form of quasi-inequations, and their satisfaction in a pointed po-monoid is equivalent, respectively, to the satisfaction of the following properties: commutativity, integrality, 0-boundedness and increasing idempotency. As we also see in Part 3, once a sequence, $\sigma$, is fixed, the $\mathring{\mathbb{M}}^{s\ell}_{\sigma}$ and $\mathring{\mathbb{M}}^{\ell}_{\sigma}$ classes are the algebraic counterpart of the $\langle\vee,*,0,1\rangle$-fragment and the $\langle\vee,\wedge,*,0,1\rangle$-fragment of the Gentzen system $\mathcal{FL}_{\sigma}$, respectively.

**Definition 11** (Exchange property)**.** *We say that a po-groupoid, $\mathcal{A}$, satisfies the* exchange property *if the following quasi-inequation is satisfied:*

$$x * y \preccurlyeq z \supset y * x \preccurlyeq z \qquad (e \preccurlyeq)$$

**Lemma 4.** *Let $\mathcal{A}$ be a po-groupoid and let $u, v, t$ be terms of its language. The following are equivalent:*

i) $\mathcal{A} \models u \preccurlyeq t \ \supset \ v \preccurlyeq t$.
ii) $\mathcal{A} \models v \preccurlyeq u$.

*Proof.* Suppose that the variables of $u$, $v$ and $t$ are in $\{x_1, \ldots, x_m\}$.
*i) $\Rightarrow$ ii)*: Let $a_1, \ldots, a_n \in A$. Since *i)*, as $u^{\mathcal{A}}(a_1, \ldots, a_n) \leq u^{\mathcal{A}}(a_1, \ldots, a_n)$, we have $v^{\mathcal{A}}(a_1, \ldots, a_n) \leq u^{\mathcal{A}}(a_1, \ldots, a_n)$.
*ii) $\Rightarrow$ i)*: Let $a_1, \ldots, a_n \in A$. Suppose that $u^{\mathcal{A}}(a_1, \ldots, a_n) \leq t^{\mathcal{A}}(a_1, \ldots, a_n)$. By *ii)*, $v^{\mathcal{A}}(a_1, \ldots, a_n) \leq u^{\mathcal{A}}(a_1, \ldots, a_n)$. Then we have $v^{\mathcal{A}}(a_1, \ldots, a_n) \leq t^{\mathcal{A}}(a_1, \ldots, a_n)$. □

As we show below, the property $(e \preccurlyeq)$ is equivalent to the commutativity of the groupoid operation and, therefore, the po-groupoids satisfying $(e \preccurlyeq)$ are, precisely, the commutative po-groupoids.

**Proposition 5.** *Let $\mathcal{A}$ be a po-groupoid. The following conditions are equivalent:*

i) *$\mathcal{A}$ satisfies the quasi-inequation $(e \preccurlyeq)$.*
ii) *$\mathcal{A}$ satisfies the inequation $x * y \preccurlyeq y * x$.*
iii) *$\mathcal{A}$ satisfies the equation $x * y \approx y * x$.*

*Proof. i)* and *ii)* are equivalent, due to Lemma 4. The equivalence between *ii)* and *iii)* is evident. □

**Definition 12** (Left-weakening property)**.** *We say that a po-monoid, $\mathcal{A}$, satisfies the* left-weakening property *if the following quasi-inequation is satisfied:*

$$x * y \preccurlyeq z \supset x * t * y \preccurlyeq z \qquad (w \preccurlyeq)$$

**Definition 13** (Integral po-monoid)**.** *We say that a po-monoid, $\mathcal{A}$, is* integral *if the unit element is the maximum with respect to the order; that is, if $\mathcal{A} \models x \preccurlyeq 1$.*

In a po-monoid, the left-weakening property is equivalent to the integrality and it is also equivalent to the fact that the result of operating two elements of the monoid is always less than or equal to either of the two element.

**Proposition 6.** *Let $\mathcal{A}$ be a po-monoid. The following conditions are equivalent:*

i) $\mathcal{A} \models x * y \preccurlyeq x$,
ii) $\mathcal{A} \models x \preccurlyeq 1$,
iii) $\mathcal{A} \models x * y \preccurlyeq y$,
iv) $\mathcal{A} \models x * z * y \preccurlyeq x * y$,
v) $\mathcal{A}$ *satisfies the quasi-inequation* $(w \preccurlyeq)$.

*Proof.* Let $a, b, c \in A$.
$i) \Rightarrow ii)$ : By the fact that 1 is the unit element and condition $i)$, we have $a = 1 * a \leq 1$.
$ii) \Rightarrow iii)$: Since $ii)$, we have $a \leq 1$. By applying monotonicity, we obtain $a * b \leq 1 * b = b$.
$iii) \Rightarrow iv)$: Since $iii)$, we have $c * b \leq b$; and by applying monotonicity, we obtain $a * c * b \leq a * b$.
$iv) \Rightarrow v)$: Due to Lemma 4.
$v) \Rightarrow i)$: As 1 is the unit element, we have $a * 1 = a$ and therefore $a * 1 \leq a$. Hence, by using $v)$, we obtain $a * b * 1 \leq a$; that is, $a * b \leq a$. □

Therefore, the po-monoids satisfying the property $(w \preccurlyeq)$ are, precisely, the integral po-monoids.

**Definition 14** (Right-weakening property)**.** *We say that a pointed po-groupoid satisfies the* right-weakening *property if the following quasi-inequation is satisfied:*

$$x \preccurlyeq 0 \supset x \preccurlyeq y \qquad (\preccurlyeq w)$$

In a pointed po-groupoid, the right-weakening property is equivalent to the fact that the distinguished element 0 is the minimum with respect to the order.

**Proposition 7.** *Let $\mathcal{A}$ be a pointed po-groupoid. The following are equivalent:*

i) $\mathcal{A}$ *satisfies the quasi-inequation* $(\preccurlyeq w)$.
ii) $\mathcal{A} \models 0 \preccurlyeq x$.

*Proof.* Let $a, b \in A$.
$i) \Rightarrow ii)$ : Given that $0 \leq 0$, by applying $i)$ we obtain $0 \leq a$.
$ii) \Rightarrow i)$: If $a \leq 0$, given that $0 \leq b$, due to transitivity we have that $a \leq b$. □

**Definition 15** (Contraction property)**.** *We say that a po-groupoid satisfies the* contraction *property if the following quasi-inequation is satisfied:*

$$x * x \preccurlyeq y \supset x \preccurlyeq y \qquad (c \preccurlyeq)$$

We say that a po-groupoid has the property of *increasing idempotency* if every element is equal to or less than the result of operating this element with itself. Next, we show that this property is equivalent to the contraction property.

**Proposition 8.** *Let $\mathcal{A}$ be a po-groupoid. The following conditions are equivalent:*

i) $\mathcal{A}$ *satisfies the quasi-inequation* $(c \preccurlyeq)$.
ii) $\mathcal{A} \models x \preccurlyeq x * x$.

*Proof.* Due to Lemma 4. □

Now we define the algebra classes $\mathring{\mathbb{M}}^{s\ell}_\sigma$ and $\mathring{\mathbb{M}}^{\ell}_\sigma$. Let $\lambda \in \{s\ell, \ell\}$. We will use the following terms:

- $\mathring{\mathbb{M}}^{\lambda}_e$, which is the class of pointed $\lambda$-monoids that satisfies $(e \preccurlyeq)$; that is, of the commutative pointed $\lambda$-monoids.
- $\mathring{\mathbb{M}}^{\lambda}_{w_l}$, which is the class of pointed $\lambda$-monoids that satisfies $(w \preccurlyeq)$; that is, of the integral pointed $\lambda$-monoids.
- $\mathring{\mathbb{M}}^{\lambda}_{w_r}$, which is the class of pointed $\lambda$-monoids that satisfies $(\preccurlyeq w)$; that is, of the pointed $\lambda$-monoids with lower bound 0.
- $\mathring{\mathbb{M}}^{\lambda}_{w}$, which is the class of pointed $\lambda$-monoids that satisfies $(w \preccurlyeq)$ and $(\preccurlyeq w)$.
- $\mathring{\mathbb{M}}^{\lambda}_{c}$, which is the class of $\lambda$-pointed monoids that satisfies $(c \preccurlyeq)$; or, in a similar way, that satisfies increasing idempotency, and which we will refer to as *contractives.*

Let $\sigma$ be a subsequence (possibly empty) of the sequence $ew_lw_rc$. If in $\sigma$ there is the sequence $w_lw_r$, we will replace it by $w$ for short. We will denote by $\mathring{\mathbb{M}}^{\lambda}_\sigma$ the class of pointed $\lambda$-monoids that satisfies the properties codified by the letters appearing in $\sigma$; and if $\sigma$ is the empty sequence, then $\mathring{\mathbb{M}}^{\lambda}_\sigma$ is the class $\mathring{\mathbb{M}}^{\lambda}$. So, for example, $\mathring{\mathbb{M}}^{\lambda}_{w_rc}$ is the class of the pointed $\lambda$-monoids that satisfies $(\preccurlyeq w)$ and $(c \preccurlyeq)$.

The $\mathring{\mathbb{M}}^{\lambda}_\sigma$ classes are subvarieties of $\mathring{\mathbb{M}}^{\lambda}$; since as we have seen, the quasi-inequations $(e \preccurlyeq)$, $(w \preccurlyeq)$, $(\preccurlyeq w)$ and $(c \preccurlyeq)$ are equivalent to inequations and, as in the classes considered the order is definable by the equation $x \vee y \approx y$, the inequations are equivalent to equations. So:

- $\mathring{\mathbb{M}}^{\lambda}_e$ is the subvariety of $\mathring{\mathbb{M}}^{\lambda}$ defined by the equation $x * y \approx y * x$;
- $\mathring{\mathbb{M}}^{\lambda}_{w_l}$ is the subvariety of $\mathring{\mathbb{M}}^{\lambda}$ defined by the equation $x \vee 1 \approx 1$ or, equivalently, by the equation $(x * y) \vee x \approx x$;
- $\mathring{\mathbb{M}}^{\lambda}_{w_r}$ is the subvariety of $\mathring{\mathbb{M}}^{\lambda}$ defined by the equation $0 \vee x \approx x$;
- $\mathring{\mathbb{M}}^{\lambda}_{c}$ is the subvariety of $\mathring{\mathbb{M}}^{\lambda}$ defined by the equation $x \vee (x * x) \approx x * x$.

By combining these equations we obtain all the subvarieties $\mathring{\mathbb{M}}^{\lambda}_\sigma$.

**Corollary 10.** *For every subsequence, $\sigma$, of the sequence $ew_lw_rc$, the $\mathring{\mathbb{M}}^{s\ell}_\sigma$ and $\mathring{\mathbb{M}}^{\ell}_\sigma$ classes are varieties.*

**Remark 5.** *From now on, in the context of semilatticed and latticed algebras, given two terms $t_1$ and $t_2$, we will use the expression $t_1 \preccurlyeq t_2$ as an abbreviation for the equation $t_1 \vee t_2 \approx t_2$. Note that the expression $t_1 \preccurlyeq t_2$ may also be seen as an atomic formula of the language of the order-algebra associated with every semilatticed monoid.*

**Proposition 9.** *In a 0-bounded integral po-monoid, the minimum element is the zero of the monoid.*[10]

*Proof.* Let $\mathcal{A}$ be a 0-bounded integral po-monoid. If $a \in A$, due to integrality, we have $a \leq 1$ and, due to monotonicity, $a * 0 \leq 1 * 0 = 0$ and $0 * a \leq 0 * 1 = 0$. However, as 0 is the minimum, $a * 0 = 0 * a = 0$. $\square$

[10] Remember that the *zero* element or the *absorbent* element of a groupoid $\mathbf{G} = \langle G, * \rangle$ is an element $0 \in G$ such that, for every $a \in G$, $a * 0 = 0 * a = 0$ is satisfied. If there is such an element, then it is unique.

**Proposition 10.** *Let* $\mathbf{A} \in \mathring{\mathbb{M}}^{\ell}$. *The following conditions are equivalent:*

i) $\mathbf{A} \models x \preccurlyeq 1$,
ii) $\mathbf{A} \models x * y \preccurlyeq x \wedge y$.

*Proof.* By Proposition 6 we have that the equation $x \preccurlyeq 1$ is equivalent to the equations $x * y \preccurlyeq x$ and $x * y \preccurlyeq y$; but in a latticed structure these two equations are equivalent to the equation $x * y \preccurlyeq x \wedge y$. □

**Proposition 11.** *Let* $\mathbf{A} \in \mathring{\mathbb{M}}^{s\ell}_{w_l}$. *Then, the following are equivalent:*

i) $\mathbf{A} \models x * x \approx x$ *(**A** is idempotent),*
ii) $\mathbf{A} \models x \preccurlyeq x * x$ *(**A** is contractive),*
iii) *For every* $a, b \in A$, *the infimum of* $a$ *and* $b$ *exists and is equal to* $a * b$,
iv) $\mathbf{A} \models x \preccurlyeq y$ *iff* $\mathbf{A} \models x * y \approx x$.

*Proof.* We only prove *ii)* $\Rightarrow$ *iii)* since the other implications are trivial. Let $a, b, c \in A$. Due to integrality we have that $a * b$ is a common lower bound to $a$ and $b$. Suppose now that $c$ is a common lower bound to $a$ and $b$. Due to monotonicity we have: $c * c \leq a * b$. However, by *ii)* we have $c \leq c * c$ and, therefore, $c \leq a * b$. □

Similarly, for the latticed varieties we have the following proposition.

**Proposition 12.** *Let* $\mathbf{A} \in \mathring{\mathbb{M}}^{\ell}_{w_l}$. *The following conditions are equivalent:*

i) $\mathbf{A} \models x * x \approx x$ *(**A** is idempotent),*
ii) $\mathbf{A} \models x \preccurlyeq x * x$ *(**A** is contractive),*
iii) $\mathbf{A} \models x * y \approx x \wedge y$,
iv) $\mathbf{A} \models x \preccurlyeq y$ *iff* $\mathbf{A} \models x * y \approx x$.

**Proposition 13.** *Let* **A** *be in* $\mathring{\mathbb{M}}^{s\ell}_{\sigma}$ *or in* $\mathring{\mathbb{M}}^{\ell}_{\sigma}$, *with* $w_l c \leq \sigma$. *Then* $\mathbf{A} \models x * y \approx y * x$. *Therefore,* $\mathring{\mathbb{M}}^{s\ell}_{w_lc} = \mathring{\mathbb{M}}^{s\ell}_{ew_lc}$, $\mathring{\mathbb{M}}^{s\ell}_{wc} = \mathring{\mathbb{M}}^{s\ell}_{ewc}$, $\mathring{\mathbb{M}}^{\ell}_{w_lc} = \mathring{\mathbb{M}}^{\ell}_{ew_lc}$ *and* $\mathring{\mathbb{M}}^{\ell}_{wc} = \mathring{\mathbb{M}}^{\ell}_{ewc}$.

*Proof.* In these varieties it is clear that the operation $*$ is commutative since $a * b$ is the infimum of $\{a, b\}$. □

**Remark 6.** *Note that* $\mathring{\mathbb{M}}^{s\ell}_{ewc}$ *is the variety of the bounded distributive lattices (cf.(Balbes & Dwinger, 1974, Chapter II), where this class is denoted by* $\mathbf{D_{01}}$*) and* $\mathring{\mathbb{M}}^{s\ell}_{ew_lc}$ *is the variety of the upper bounded distributive lattices.*

**Remark 7.** *The varieties* $\mathring{\mathbb{M}}^{s\ell}_{ew_lc}$ *and* $\mathring{\mathbb{M}}^{\ell}_{ew_lc}$ *are definitionally equivalent and the operations* $*$ *and* $\wedge$ *are the same. This also applies to the varieties* $\mathring{\mathbb{M}}^{s\ell}_{ewc}$ *and* $\mathring{\mathbb{M}}^{\ell}_{ewc}$.

## 7. Relative Pseudocomplements, Residuated Structures.

In this section we present the residuated structures related to the systems introduced in Section 3 Some of the content is based on Dubreil-Jacotin et al. (1953), Birkhoff (1973), Jipsen & Tsinakis (2002) and Galatos et al. (2007a). A novel element we should emphasize is that we use the notion of *relative pseudocomplement.* This notion is a generalization of the same notion used traditionally within the framework of lattices.

7.1. **Relative pseudocomplement.** Recall that, given a lattice $\mathbf{L}$ and two elements $a, b \in L$, if the largest element $x \in L$ such that $a \wedge x \leq b$ exists, this element is denoted by $a \to b$ and is referred to as the *relative pseudocomplement of* $a$ *with respect to* $b$. We should bear in mind that a Heyting algebra is a lattice $\mathbf{L}$ with a minimum element such that the relative pseudocomplement of $a$ with respect to $b$ exists for every $a, b \in L$ (Cf. (Balbes & Dwinger, 1974, Chapter IX)). In order to introduce the notion of the relative pseudocomplement, we start from an ordered set, $\langle A, \leq \rangle$, in which we have defined a binary operation, $*$.

**Definition 16** (Relative Pseudocomplement)**.** *Let* $*$ *be a binary operation defined in an ordered set* $\langle A, \leq \rangle$*. Given* $a, b \in A$*, if the largest element* $x \in A$ *such that* $a*x \leq b$ $(x*a \leq b)$ *exists, we say that this element is the* $*$-right (left) relative pseudocomplement of $a$ with respect to $b$. *If the operation* $*$ *is commutative, then the notions of left and right relative pseudocomplement coincide and we call this element simply the* relative pseudocomplement of $a$ with respect to $b$.

Note that if the order defined in the set is a $\wedge$-semilattice or a lattice, then the notion of $\wedge$-relative pseudocomplement and the traditional notion of relative pseudocomplement coincide.

**Remark 8.** *Given a binary operation* $*$ *defined in an ordered set, when the context does not lend itself to any confusion, we will use the name (right or left) relative pseudocomplement instead of relative* $*$*-pseudocomplement.*

*The operation that assigns to every pair of elements* $a, b \in A$ *the* $*$*-right (left) relative pseudocomplement of* $a$ *with respect to* $b$ *is called* $*$-right (left) relative pseudocomplementation.

7.2. **Residuation, residuated structures.** We now recap the notion of *residuation* and establish the link between this notion and that of the relative pseudocomplement.

**Definition 17** (Residuated operation)**.** *A binary operation* $*$ *defined in a partially-ordered set* $\langle A, \leq \rangle$ *is called* residuated *if there are two binary operations* $\backslash$ *and* $/$ *defined in* $A$ *such that, for every* $a, b, c \in A$*:*

(LR) $$a * b \leq c \quad \textit{iff} \quad b \leq a \backslash c \quad \textit{iff} \quad a \leq c/b.$$

*This condition is called the* Law of Residuation*, and the operations* $\backslash$ *and* $/$ *are called the* right residual *and* left residual *of the operation* $*$*, respectively. We also call them* right *and* left residuation *of* $*$*, respectively. Given* $a, b \in A$*, the element* $a \backslash b$ *will be referred to as the* right residual of $a$ relative to $b$ *and the element* $b/a$ *will be referred to as the* left residual of $a$ relative to $b$*. In the event that the operation is commutative, the two residuals coincide and the corresponding operation is simply called the* residual *(or* residuation*) of the operation* $*$*. In this case, we will use the symbol* $\to$ *and annotate* $a \to b$ *instead of* $a \backslash b$ *or* $b/a$*.*

A *basic residuated structure* is an order-algebra $\mathcal{A} = \langle A, *, \backslash, /, \leq \rangle$ with the algebraic type $\langle 2, 2, 2 \rangle$ such that $*$ is a residuated operation with respect to the order and such that the operations $\backslash$ and $/$ are the right and left residuals of the operation $*$, respectively. We will call *residuated structure* any structure having a basic residuated structure as a reduct. Residuated structures are characterized by the fact that the operation $*$ is monotonous with respect to the order and by the fact that, for every $a, b$ in the universe, there exist the right and left pseudocomplements of $a$ relative to $b$. This last point implies the uniqueness of the residual.

**Proposition 14.** *Let $\langle A, \leq \rangle$ be a partially-ordered set and let $*$ be a binary operation defined in $A$. For every $a, b \in A$, we define the sets:*

$$\mathrm{R}^b(a) = \{x \in A : a * x \leq b\} \qquad \textit{and} \qquad \mathrm{L}^b(a) = \{x \in A : x * a \leq b\}.$$

*Then, the following conditions are equivalent:*

i) *The operation $*$ is residuated.*
ii) *The operation $*$ is monotonous and, for every $a, b \in A$, there exist the right and the left pseudocomplements of $a$ with respect to $b$, i.e., the sets $\mathrm{R}^b(a)$ and $\mathrm{L}^b(a)$ have a maximum element.*

*Under these conditions, for every $a, b \in A$, $a \backslash b = max\, \mathrm{R}^b(a)$ and $b/a = max\, \mathrm{L}^b(a)$.*

*Proof. i) $\Rightarrow$ ii)*: Suppose that $*$ is residuated and let $a, b \in A$. We will see that $a \backslash b$ is the maximum element of $\mathrm{R}^b(a)$. By applying (LR), from $a \backslash b \leq a \backslash b$ we obtain $a * (a \backslash b) \leq b$ and, therefore, $a \backslash b \in \mathrm{R}^b(a)$. Suppose now that $c \in \mathrm{R}^b(a)$: then we have $a * c \leq b$ which, by (LR), is equivalent to $c \leq a \backslash b$. Therefore, $a \backslash b$ is the maximum of $\mathrm{R}^b(a)$. Similarly we obtain that $b/a$ is the maximum of $\mathrm{L}^b(a)$. Now we will see that monotonicity is satisfied. Let $a, b, c, d \in A$ and suppose that $a \leq c$ and $b \leq d$. ¿From $c * d \leq c * d$, we obtain $b \leq d \leq c \backslash (c * d)$ and therefore $c * b \leq c * d$. Based on this, we obtain $a \leq c \leq (c * d)/b$ and consequently $a * b \leq c * d$.

*ii) $\Rightarrow$ i)*: For every $a, b \in A$ we define:

$$a \backslash b := \max \mathrm{R}^b(a) \textit{and} \qquad a/b := \max \mathrm{L}^a(b).$$

Given that for every $a, b \in A$ the sets $\mathrm{R}^b(a)$ and $\mathrm{L}^a(b)$ have a maximum, the operations $\backslash$ and $/$ are well defined. Let $c \in A$ and suppose that $a * c \leq b$. Then $c \in \mathrm{R}^b(a)$ and, therefore, $c \leq a \backslash b$. Now suppose $c \leq a \backslash b$. According to left monotonicity, we obtain $a * c \leq a * (a \backslash b)$. However, as $a \backslash b = \max \mathrm{R}^b(a)$, in particular $a \backslash b \in \mathrm{R}^b(a)$ and hence, $a * (a \backslash b) \leq b$. So, $a * c \leq b$. Similarly, from $c * a \leq b$ we obtain $c \leq b/a$ and from $c \leq b/a$, by applying right monotonicity, we obtain: $c * a \leq b$. □

**Corollary 11** (Uniqueness of residuals)**.** *If a binary operation defined in an ordered set $\langle A, \leq \rangle$ is residuated, then there are precisely two binary operations $\backslash$ and $/$ that satisfy (LR).*

We observe that, as a consequence of Proposition 14, we have that in the framework of po-groupoids, if the operation is residuated, the notions of residual and relative pseudocomplement coincide.

**Corollary 12.** *Let $\langle A, *, \leq \rangle$ be a po-groupoid. Then the following are equivalent:*

i) *The operation $*$ is residuated.*

ii) *For every $a, b \in A$, there are right and left relative $*$-pseudocomplements of $a$ with respect to $b$.*

*Under these conditions, for every $a, b \in A$, the right residual of $a$ relative to $b$ is the right relative pseudocomplement of $a$ with respect to $b$ and the left residual of $a$ relative to $b$ is the left relative pseudocomplement of $a$ with respect to $b$.*

**Definition 18** (Residuated po-groupoid)**.** *A* residuated po-groupoid *is an order-algebra $\mathcal{A} = \langle A, *, \backslash, /, \leq \rangle$, where $\langle A, *, \leq \rangle$ is a po-groupoid, the operation $*$ is residuated and the operations $\backslash$ and $/$ are its residuals.*

**Remark 9** (Convention)**.** *In accordance with the preceding results, we have that the basic residuated structures are precisely the residuated po-groupoids. In the nomenclature we use here, for practical reasons we will dispense with the term* partially ordered *or the prefix* po*, because the notion of* residuation *implicitly entails the presence of a partial order in the structure in such a way that the operation $*$ is monotonous with respect to that order. Residuals constitute a generalization of the division operation in the groups. In concordance with this idea, $a\backslash b$ is read as "a under b" and $b/a$ is read as "b above a". In both cases we may say that $b$ is the* numerator *and $a$ the* denominator*.*

As the operation $*$ of a po-groupoid is compatible with the order, the residuals of a residuated groupoid are also connected with the order in the following sense: the right (left) residual is antimonotonous in the first (second) argument and monotonous in the second (first) argument.

**Proposition 15.** *In every residuated groupoid $\mathcal{A}$ the following conditions, for any $a, b, c \in A$, are satisfied:*

*i) if $a \leq b$, then $c\backslash a \leq c\backslash b$ and $b\backslash c \leq a\backslash c$,*
*ii) if $a \leq b$, then $a/c \leq b/c$ and $c/b \leq c/a$.*

*Proof. i)*: Suppose $a \leq b$. According to reflexivity, we have that $c\backslash a \leq c\backslash a$. Hence, by applying $(LR)$ we obtain $c * (c\backslash a) \leq a$ and, therefore, $c * (c\backslash a) \leq b$ and again by $(LR)$ we obtain: $c\backslash a \leq c\backslash b$. In contrast, by applying monotonicity, from $a \leq b$ we obtain $a * (b\backslash c) \leq b * (b\backslash c)$ and from $b\backslash c \leq b\backslash c$, by $(LR)$, we obtain $b * (b\backslash c) \leq c$ and thus, $a * (b\backslash c) \leq c$ which is equivalent to $b\backslash c \leq a\backslash c$.
*ii)*: Proved in the same way. □

In a residuated groupoid, the operation $*$ preserves the existing suprema in each argument and the residuals preserve all the existing infima in the numerator and turn the existing suprema into infima in the denominator, as is shown in the following propositions.

**Remark 10.** *If $\{a_i : i \in I\}$ is a family of elements of a partially-ordered set $\langle A, \leq \rangle$, then the supremum and the infimum (if they exist) in $A$ of the family will be denoted by $\bigvee_{i \in I} a_i$ and $\bigwedge_{i \in I} a_i$, respectively.*

**Proposition 16** (Generalized distributivity)**.** *Let $\mathcal{A} = \langle A, *, \backslash, /, \leq \rangle$ be a residuated groupoid and $\{a_i : i \in I\}$ and $\{b_j : j \in J\}$ two families of elements of $A$. If $\bigvee_{i\in I} a_i$ and $\bigvee_{j\in J} b_j$ exist, then there exists $\bigvee_{\langle i,j\rangle \in I\times J} a_i * b_j$ and the following holds:*

$$\bigvee_{i\in I} a_i * \bigvee_{j\in J} b_j = \bigvee_{\langle i,j\rangle \in I\times J} a_i * b_j .$$

*Proof.* According to monotonicity it is clear that, for every $\langle i, j\rangle \in I \times J$, $\bigvee_{i\in I} a_i * \bigvee_{j\in J} b_j$ is an upper bound of $a_i * b_j$. Suppose that $a_i * b_j \leq c$. Then, by applying (LR) we have $b_j \leq a_i\backslash c$ and therefore $\bigvee_{j\in J} b_j \leq a_i\backslash c$. Once again by (LR) we obtain $a_i * \bigvee_{j\in J} b_j \leq c$ and, hence, $a_i \leq c/\bigvee_{j\in J} b_j$ and, thus, $\bigvee_{i\in I} a_i \leq c/\bigvee_{j\in J} b_j$ which, again by (LR), allows us to conclude $\bigvee_{i\in I} a_i * \bigvee_{j\in J} b_j \leq c$. □

**Proposition 17.** *Let $\mathcal{A} = \langle A, *, \backslash, /, \leq \rangle$ be a residuated groupoid and $\{a_i : i \in I\}$ and $\{b_j : j \in J\}$ two families of elements of $A$. If $\bigvee_{i\in I} a_i$ and $\bigwedge_{j\in J} b_j$ exist, then, for every $c \in A$, there exist $\bigwedge_{i\in I} a_i\backslash c$, $\bigwedge_{j\in J} c\backslash b_j$, $\bigwedge_{i\in I} c/a_i$ and $\bigwedge_{j\in J} b_j/c$ and the following are satisfied:*

$$(\bigvee_{i\in I} a_i)\backslash c = \bigwedge_{i\in I} a_i\backslash c\,;\ c\backslash(\bigwedge_{j\in J} b_j) = \bigwedge_{j\in J} c\backslash b_j\,;$$

$$c/(\bigvee_{i\in I} a_i) = \bigwedge_{i\in I} c/a_i\,;\ (\bigwedge_{j\in J} b_j)/c = \bigwedge_{j\in J} b_j/c.$$

*Proof.* From $a_i \leq \bigvee_{i\in I} a_i$, due to the antimonotonicity in the first argument of the right residual, we obtain $(\bigvee_{i\in I} a_i)\backslash c \leq a_i\backslash c$. Therefore, $(\bigvee_{i\in I} a_i)\backslash c$ is a lower bound of $a_i\backslash c$. Let $d \in A$ and suppose that, for every $i \in I$, $d \leq a_i\backslash c$. This is equivalent to $a_i * d \leq c$ which in turn is equivalent to $a_i \leq c/d$. Hence, we obtain $(\bigvee_{i\in I} a_i) \leq c/d$, which is equivalent to $(\bigvee_{i\in I} a_i) * d \leq c$ and, therefore, to $d \leq (\bigvee_{i\in I} a_i)\backslash c$. Consequently, $(\bigvee_{i\in I} a_i)\backslash c$ is the infimum of $a_i\backslash c$.

From $\bigwedge_{j\in J} b_j \leq b_j$, due to the monotonicity of the second argument of the right residual, we obtain $c\backslash(\bigwedge_{j\in J} b_j) \leq c\backslash b_j$. Suppose that, for every $j \in J$, $d \leq c\backslash b_j$. This is equivalent to $c*d \leq b_j$. Hence we obtain $c*d \leq (\bigwedge_{j\in J} b_j)$, which is equivalent to $d \leq c\backslash(\bigwedge_{j\in J} b_j)$. Therefore, $(\bigwedge_{j\in J} b_j)$ is the infimum of $c\backslash b_j$.

The other two equalities are proved similarly using (LR), the antimonotonicity of the second argument and the monotonicity of the first argument of the left residual. □

In Section 6 above, we saw that every $s\ell$-groupoid defines a *po*-groupoid (Proposition 1) and also that distributivity is a stronger condition than monotonicity (Proposition 2) and, therefore, it is not true in general that a *po*-groupoid that is a semilattice under its partial-order relation is an $s\ell$-groupoid. However, this will be true whenever the operation $*$ is residuated:

**Proposition 18.** *Let $\langle A, *, \leq\rangle$ be a residuated groupoid that is a semilattice under its partial-order relation. We define $x \vee y =: \bigvee\{x, y\}$. Then $\langle A, \vee, *\rangle$ is an $s\ell$-groupoid.*

*Proof.* It is a direct consequence of Proposition 16. □

**Definition 19** (Residuated monoid)**.** *A* residuated monoid *is an order-algebra*

$$\mathcal{A} = \langle A, *, \backslash, /, 1, \leq\rangle$$

*such that $\langle A, *, \backslash, /, \leq\rangle$ is a residuated groupoid and such that $\langle A, *, 1\rangle$ is a monoid.*

**Definition 20** (Residuated lattice)**.** *A* residuated lattice *is an algebra*

$$\mathbf{A} = \langle A, \vee, \wedge, *, \backslash, /, 1\rangle$$

*such that $\langle A, \vee, \wedge\rangle$ is a lattice and $\langle A, *, \backslash, /, 1, \leq\rangle$, where $\leq$ is the order of the lattice, is a residuated monoid. We will denote the class of residuated lattices by $\mathbb{RL}$.*

**Definition 21** (Pointed residuated lattice)**.** *A* pointed residuated lattice *is an algebra*

$$\mathbf{A} = \langle \vee, \wedge, *, \backslash, /, 0, 1\rangle$$

*such that its $\langle \vee, \wedge, *, \backslash, /, 1\rangle$-reduct is a residuated lattice and such that* 0 *is a fixed, but arbitrary, element of $A$.*

Pointed residuated lattices are called *full Lambek algebras* according to Ono (see for instance Ono (1993)) on account of their connection with the sequent calculus **FL**. We will denote by $\mathbb{FL}$ the class of pointed residuated lattices and will call its members $\mathbb{FL}$-*algebras*. Observe that (pointed) residuated lattices are (pointed) $\ell$-monoids such that their monoidal operation is residuated. Residuated lattices and $\mathbb{FL}$-algebras can be understood as order-algebras $\langle \mathbf{A}, \leq\rangle$, where $\leq$ is the order defined by the lattice.

**Definition 22** (Mirror image.)**.** *If $t$ is a term of an algebraic language $\mathcal{L}$ such that $\langle *, \backslash, /\rangle \leq \mathcal{L} \leq \langle \vee, \wedge, *, \backslash, /, 0, 1\rangle$, we define its* mirror image *$\mu(t)$ inductively on the complexity of $t$:*[11]

$$\mu(t) := \begin{cases} t, & \text{if } t \in Var \text{ or } t \in \{0,1\}, \\ \mu(u) \vee \mu(v), & \text{if } t = u \vee v, \\ \mu(u) \wedge \mu(v), & \text{if } t = u \wedge v, \\ \mu(v) * \mu(u), & \text{if } t = u * v, \\ \mu(u)\backslash\mu(v), & \text{if } t = v/u, \\ \mu(v)/\mu(u), & \text{if } t = u\backslash v. \end{cases}$$

*We define the* mirror image of a formula *of the first-order language with equality $\mathcal{L}^{\preccurlyeq} = \langle \mathcal{L}, \preccurlyeq\rangle$ as the formula obtained by replacing all the existing terms therein with their mirror images.*

[11] *The complexity of a term* means the number of functional occurrences of arity $k \geq 1$ in this term.

**Lemma 5.** *Let $\mathbb{K}$ be the class of residuated groupoids, residuated monoids, residuated lattices or $\mathbb{FL}$-algebras. We denote their algebraic language by $\mathcal{L}_{\mathbb{K}}$ and we denote by $\mathcal{L}_{\mathbb{K}}^{\preccurlyeq}$ the first-order language with equality $\langle \mathcal{L}_{\mathbb{K}}, \preccurlyeq \rangle$. If $\mathcal{A} \in \mathbb{K}$, consider the $\mathcal{L}_{\mathbb{K}}^{\preccurlyeq}$-structure $\mathcal{A}'$, with a universe equal to that of $\mathcal{A}$, with the same order as $\mathcal{A}$, and where the operations and constants of $\mathcal{A}'$ in $\{\vee, \wedge, 0, 1\}$ are the same as in $\mathcal{A}$ and the remaining operations are defined as follows: for every $a, b \in A$,*

$$a *^{\mathcal{A}'} b := b *^{\mathcal{A}} a, \quad a \backslash^{\mathcal{A}'} b := b /^{\mathcal{A}} a, \quad b /^{\mathcal{A}'} a := a \backslash^{\mathcal{A}} b.$$

*Then,*

*i) $\mathcal{A}'$, which we will name the* opposite *of $\mathcal{A}$, belongs to $\mathbb{K}$,*
*ii) for every term $t$ of $\mathcal{L}_{\mathbb{K}}$, $\mu(t)^{\mathcal{A}'} = t^{\mathcal{A}}$ is satisfied.*

*Proof. i)*: It is easy to see that the operation $*^{\mathcal{A}'}$ is residuated and that $\backslash^{\mathcal{A}'}$ and $/^{\mathcal{A}'}$ are its right and left residual, respectively.
*ii)*: By induction on the complexity of the term $t$. Suppose that the variables appearing in $t$ are in $\{x_1, \ldots, x_n\}$. We must prove that for every assignment $\bar{a}$ of the variables in $A$, if this assignment is such that $\bar{a}(x_i) = a_i$, for every $1 \leq i \leq n$, then $\mu(t)^{\mathcal{A}'}(a_1, \ldots, a_n) = t^{\mathcal{A}}(a_1, \ldots, a_n)$. If $t$ is a variable or $t \in \{u \vee v, u \wedge v, 0, 1\}$, then it is obvious. Suppose $t = u\backslash v$, where $u$ and $v$ are $\langle *, \backslash, / \rangle$-terms. Then we have:

$$\mu(u\backslash v)^{\mathcal{A}'}(a_1, \ldots, a_n) = (\mu(v)/\mu(u))^{\mathcal{A}'}(a_1, \ldots, a_n) =$$
$$= \mu(v)^{\mathcal{A}'}(a_1, \ldots, a_n)/^{\mathcal{A}'}\mu(u)^{\mathcal{A}'}(a_1, \ldots, a_n),$$

and by applying the induction hypothesis:

$$\mu(v)^{\mathcal{A}'}(a_1, \ldots, a_n)/^{\mathcal{A}'}\mu(u)^{\mathcal{A}'}(a_1, \ldots, a_n) = v^{\mathcal{A}}(a_1, \ldots, a_n)/^{\mathcal{A}'}u^{\mathcal{A}}(a_1, \ldots, a_n) =$$
$$= u^{\mathcal{A}}(a_1, \ldots, a_n)\backslash^{\mathcal{A}}v^{\mathcal{A}}(a_1, \ldots, a_n) = (u\backslash v)^{\mathcal{A}}(a_1, \ldots, a_n).$$

The cases $t = v/u$ and $t = u * v$ are similar. □

**Lemma 6.** *Let $\mathcal{A}$ be a residuated groupoid (residuated monoid, residuated lattice, $\mathbb{FL}$-algebra) and let $\mathcal{A}'$ be its opposite. Then the following is satisfied for every formula $\varphi$ of the first-order language of the residuated groupoids (residuated monoids, residuated lattices, $\mathbb{FL}$-algebras):*

$$\mathcal{A} \vDash \varphi \quad \textit{iff} \quad \mathcal{A}' \vDash \mu(\varphi).$$

*Proof.* Let $\mathbb{K}$ and $\mathcal{L}_{\mathbb{K}}^{\preccurlyeq}$ be as in Lemma 5 and let $\varphi$ be an $\mathcal{L}_{\mathbb{K}}^{\preccurlyeq}$-formula. We will see that if $\mathcal{A} \in \mathbb{K}$, then $\mathcal{A} \nvDash \varphi$ if and only if $\mathcal{A}' \nvDash \mu(\varphi)$. The proof follows by induction on the complexity of $\varphi$.[12] If $\varphi$ is an atomic formula, it will be an equation or an inequation. If it is in $t_1 \preccurlyeq t_2$ or $t_1 \approx t_2$ and variables in terms $t_1$ and $t_2$ are in $\{x_1, \ldots, x_m\}$, then we have that $\mathcal{A} \nvDash \varphi$ is equivalent to the fact that there are elements $a_1, \ldots, a_m$ such that:

$$t_1^{\mathcal{A}}(a_1, \ldots, a_m) > t_2^{\mathcal{A}}(a_1, \ldots, a_m) \text{ or } t_1^{\mathcal{A}}(a_1, \ldots, a_m) \neq t_2^{\mathcal{A}}(a_1, \ldots, a_m),$$

[12] *The complexity of a first-order formula* is the number of occurrences of the boolean operators and the quantifiers.

which, according to Lemma 5, is equivalent to:

$$\mu(t_1)^{\mathcal{A}'}(a_1, \ldots, a_m) > \mu(t_2)^{\mathcal{A}'}(a_1, \ldots, a_m) \text{ or } \mu(t_1)^{\mathcal{A}'}(a_1, \ldots, a_m) \neq \mu(t_2)^{\mathcal{A}'}(a_1, \ldots, a_m),$$

that is, $\mathcal{A}' \nvDash \varphi$. The remaining proof is a simple and routine task. □

**Theorem 9** (Law of Mirror Images). *A formula is valid in the class of the residuated groupoids (residuated monoids, residuated lattices, $\mathbb{FL}$-algebras) if and only if it is its own mirror image.*

*Proof.* Due to Lemma 6. □

**Corollary 13.** *A quasi-inequation (inequation, quasi-equation, equation) is valid in the class of the residuated groupoids (residuated monoids, residuated lattices, $\mathbb{FL}$-algebras) if and only if it is its own mirror image.*

**Remark 11.** *Note that every subclass of the classes considered that is defined by a set of formulas and their mirror images satisfies the Law of Mirror Images.*

7.3. **Properties of residuated monoids.** In the following proposition, we give some properties of residuated monoids which are easy to prove.

**Proposition 19.** *In all residuated monoids, the following inequations and equations (and their mirror images) are satisfied:*

a) $x * (x \backslash y) \preccurlyeq y$,
b) $1 \preccurlyeq x \backslash x$,
c) $(x \backslash y) * z \preccurlyeq x \backslash (y * z)$,
d) $x \backslash y \preccurlyeq (z * x) \backslash (z * y)$,
e) $(x \backslash y) * (y \backslash z) \preccurlyeq x \backslash z$,
f) $(x * y) \backslash z \approx y \backslash (x \backslash z)$,
g) $x \backslash (y / z) \approx (x \backslash y) / z$,
h) $(x \backslash 1) * y \preccurlyeq x \backslash y$.
i) $x * (x \backslash x) \approx x$,
j) $(x \backslash x) * (x \backslash x) \approx x \backslash x$.

In the following proposition we give some properties of residuated groupoids (and, therefore, of residuated monoids) with a minimum element.

**Proposition 20.** *If a residuated groupoid $\mathcal{A}$ has a minimum element $\bot$, then the element $\bot \backslash \bot$ ($\bot / \bot$) is the maximum element of $A$. Furthermore, for every $a \in A$, we have:*

$$i) \quad a * \bot = \bot = \bot * a\,, \qquad ii) \quad \bot \backslash a = \top = a \backslash \top, \qquad iii) \quad a / \bot = \top = \top / a\,,$$

*where we denote by $\top$ the element $\bot \backslash \bot = \bot / \bot$.*

*Proof.* Let $a \in A$. As $\bot$ is the minimum element, we have $\bot \leq a \backslash \bot$ and, by applying (LR), this is equivalent to $a * \bot \leq \bot$ which in turn is equal to $a \leq \bot / \bot$. Furthermore, $a \leq \bot \backslash \bot$ is obtained as the mirror image. So, $\top := \bot \backslash \bot = \bot / \bot$ is the maximum element of $A$.
*i*) Given that $a \leq \bot \backslash \bot$, we have: $a * \bot \leq \bot$ and, as $\bot$ is the minimum, $a * \bot = \bot$. According to the Law of Mirror Images, $\bot * a = \bot$ is obtained.

*ii*) As $\bot * (\bot\backslash\bot) \leq \bot$, we have: $\bot * \top \leq \bot$ and, therefore, for every $a \in A$, $\bot * \top \leq a$ which is equal to $\top \leq \bot\backslash a$. Thus, $\bot\backslash a = \top$. Based on the fact that $\top$ is the maximum, we have: $a * \top \leq \top$, which is equivalent to $\top \leq a\backslash\top$. Therefore, $\top = a\backslash\top$.
*iii*) The equalities referred to are mirror images of the equalities in *ii*). □

The class of the integral residuated monoids (that is, of the residuated monoids where the unit element of the monoid is the maximum element with respect to the order) is definitionally equivalent to the class formed by all its algebraic reducts, as is seen in the following proposition.

**Proposition 21.** *If $\mathcal{A}$ is an integral residuated monoid, then the following conditions are equivalent:*

i) $\mathcal{A} \models x \preccurlyeq y$,
ii) $\mathcal{A} \models x\backslash y \approx 1$,
iii) $\mathcal{A} \models y/x \approx 1$.

*Proof.* Pursuant to the Law of Mirror Images, the proof of equivalence of the first two items will be enough. Let $a, b \in A$ be such that $a \leq b$. We have $a * 1 \leq b$ and this, according to the law of residuation, is equivalent to $1 \leq a\backslash b$ but due to integrality, this is equal to $1 = a\backslash b$. □

So, the class of integral residuated monoids can be defined as a class of algebras $\mathbf{A} = \langle A, *, \backslash, /, 1\rangle$ of type $\langle 2, 2, 2, 0\rangle$. Clearly, it is a quasivariety. The class of commutative integral residuated monoids is known in the literature by the acronym $\mathbb{POCRIM}$ (*partially-ordered commutative residuated integral monoids*). It is a quasivariety that is not a variety (see Higgs (1984)), since it is not a class closed by homomorphic images. As a consequence of this fact, we have that the quasivariety corresponding to the non-commutative case will not be a variety either.

7.4. **Properties of residuated lattices.** In the following proposition we give some properties for the residuated lattices which are a consequence of (*i*) from Propositions 16 and 17.

**Proposition 22.** *In every residuated lattice* $\mathbf{A}$ *the following equations are satisfied:*

1) $(x \vee y)\backslash z \approx (x\backslash z) \wedge (y\backslash z)$,
2) $z/(x \vee y) \approx (z/x) \wedge (z/y)$,
3) $z\backslash(x \wedge y) \approx (z\backslash x) \wedge (z\backslash y)$,
4) $(x \wedge y)/z \approx (x/z) \wedge (y/z)$.

The class of residuated lattices is a variety. Below we give an equational base (see Jipsen & Tsinakis (2002)).

**Theorem 10** (Equational presentation of $\mathbb{RL}$). $\mathbb{RL}$ *is the equational class of algebras* $\mathbf{A} = \langle A, \vee, \wedge, *, \backslash, /, 1\rangle$ *of type* $\langle 2, 2, 2, 2, 2, 0\rangle$ *that satisfies:*[13]

1) *Any set of equations defining the class of lattices,*

[13] Recall that we use the inequation $t_1 \preccurlyeq t_2$ as an abbreviation for the equation $t_1 \vee t_2 \approx t_2$.

2) *Any set of equations defining the class of the monoids with an identity element 1,*
3) *(r)* $x * ((x\backslash z) \wedge y) \preccurlyeq z$*; (l)* $((z/x) \wedge y) * x \preccurlyeq z$*,*
4) *(r)* $y \preccurlyeq x\backslash((x * y) \vee z)$*; (l)* $y \preccurlyeq ((y * x) \vee z)/x$*.*

7.5. **The $\mathbb{FL}_\sigma$ varieties.** In this work, we use the presentation of the class of the $\mathbb{FL}$-algebras in the language $\langle \vee, \wedge, *, \backslash, /, ^{\backprime}, ', 0, 1\rangle$ including the right and left negation operations as primitives connectives. The reason for this is that in this paper we study some fragments without implication and with negation of the logic systems in question. Now, we give an equational presentation of the class $\mathbb{FL}$ in the language $\langle \vee, \wedge, *, \backslash, /, ^{\backprime}, ', 0, 1\rangle$.

**Theorem 11** (Equational presentation of $\mathbb{FL}$)**.** $\mathbb{FL}$ *is the equational class of the algebras* $\mathbf{A} = \langle A, \vee, \wedge, *, \backslash, /, ^{\backprime}, ', 0, 1\rangle$ *of type* $\langle 2, 2, 2, 2, 2, 1, 1, 0, 0\rangle$ *that satisfies:*

1) *Any set of equations defining the lattices,*
2) *Any set of equations defining the monoids with an identity element 1,*
3) *(r)* $x * ((x\backslash z) \wedge y) \preccurlyeq z$*; (l)* $((z/x) \wedge y) * x \preccurlyeq z$*,*
4) *(r)* $y \preccurlyeq x\backslash((x * y) \vee z)$*; (l)* $y \preccurlyeq ((y * x) \vee z)/x$*,*
5) *(r)* $x^{\backprime} \approx x\backslash 0$*; (l)* $'x \approx 0/x$*.*

Much as we did for pointed monoids, we will define $\mathbb{FL}_\sigma$ classes; that is, subclasses defined by the properties $(e \preccurlyeq)$, $(w \preccurlyeq)$, $(\preccurlyeq w)$ and $(c \preccurlyeq)$. We observe that the $\langle \vee, *, 0, 1\rangle$-reduct of an $\mathbb{FL}$-algebra is an $\mathring{\mathbb{M}}^{s\ell}$-algebra, since residuated lattices satisfy the distributivity of the operation $*$ with respect to the operation $\vee$.

**Definition 23** ($\mathbb{FL}_\sigma$-algebra)**.** *An* $\mathbb{FL}_\sigma$*-algebra is an* $\mathbb{FL}$*-algebra such that its* $\langle \vee, *, 0, 1\rangle$*-reduct is an* $\mathring{\mathbb{M}}^{s\ell}_\sigma$*-algebra.*

Note that in the $\mathbb{FL}_e$-algebras, the equations $x\backslash y \approx y/x$ and $x^{\backprime} \approx {}'x$ are satisfied on account of the commutativity of the monoidal operation. For this reason, the class $\mathbb{FL}_e$ is presented as having only one residual and one negation which, in accordance with its logical interpretation, are denoted by $\to$ and $\neg$, respectively. For future reference, we give an equational presentation of the class $\mathbb{FL}_e$ below.

**Theorem 12** (Equational presentation of $\mathbb{FL}_e$)**.** $\mathbb{FL}_e$ *is the equational class of algebras* $\mathbf{A} = \langle A, \vee, \wedge, *, \to, \neg, 0, 1\rangle$ *of type* $\langle 2, 2, 2, 2, 1, 0, 0\rangle$ *that satisfies:*

1) *Any set of equations defining the lattices,*
2) *Any set of equations defining the monoids with an identity element 1,*
3) $x * ((x \to z) \wedge y) \preccurlyeq z$*,*
4) $y \preccurlyeq x \to ((x * y) \vee z)$*,*
5) $\neg x \approx x \to 0$*.*

Clearly, all the $\mathbb{FL}_\sigma$ classes are varieties. For the equational presentation of a class $\mathbb{FL}_\sigma$, we will add the following to the equations defining $\mathbb{FL}$ or $\mathbb{FL}_e$:

- $x \preccurlyeq 1$ or $x * y \preccurlyeq x$ if $w_l \leq \sigma$,
- $0 \preccurlyeq x$ if $w_r \leq \sigma$,
- $x \preccurlyeq x * x$ if $c \leq \sigma$.

**Remark 12.** *Class $\mathbb{FL}_{ewc}$ is definitionally equivalent to the class of Heyting algebras, the semantic counterpart of intuitionistic logic. For more details on Heyting algebras, refer to Balbes & Dwinger (1974).*

## 8. Adding Negation Operators: Pseudocomplemented Structures.

In this section, we present the notion of *pseudocomplementation* in the framework of *pointed groupoids* and we define the $\mathbb{PM}^{\preccurlyeq}$ class of pseudocomplemented po-monoids and the classes $\mathbb{PM}^{s\ell}$ and $\mathbb{PM}^{\ell}$ of semilatticed and latticed pseudocomplemented monoids. The notion of pseudocomplement with respect to the monoidal operation can be seen as a generalization of the same notion defined in the framework of pseudocomplemented distributive lattices (see Balbes & Dwinger (1974); Lakser (1973)). We show that the class $\mathbb{PM}^{\preccurlyeq}$ can be defined by means a set of inequations and thus the $\mathbb{PM}^{s\ell}$ and $\mathbb{PM}^{\ell}$ classes are varieties (Sections 8.2 and 8.3). Section 8.4 analyzes the case in which the pseudocomplementation is with respect to the minimum element of the monoid. Pseudocomplementation constitute the algebraic counterpart of negation: in Part 5 we will state the connection between, on the one hand, the varieties $\mathbb{PM}^{s\ell}_{\sigma}$ and $\mathbb{PM}^{\ell}_{\sigma}$ (subvarieties of $\mathbb{PM}^{s\ell}$ and $\mathbb{PM}^{\ell}$ defined by the equations codified by $\sigma$) and on the other hand, the fragments of the Gentzen system $\mathcal{FL}_{\sigma}$ and the associated external system $\mathfrak{e}\mathcal{FL}_{\sigma}$ in the languages $\langle \vee, *, ^{\backprime}, ', 0, 1\rangle$ and $\langle \vee, \wedge, *, ^{\backprime}, ', 0, 1\rangle$.

### 8.1. Operations of pseudocomplementation.

In this subsection, we introduce the operations of left and right pseudocomplementation in the general context of pointed po-groupoids.

**Proposition 23.** *Let $\mathcal{A} = \langle A, *, 0, \leq\rangle$ be a pointed po-groupoid. The following conditions are equivalent:*

i) *For every $a \in A$, there exist right and left relative pseudocomplements of $a$ with respect to* 0*.*
ii) *There exist two unary operations $^{\backprime}$ and $'$ defined in $A$ such that, for every $a, b \in A$:*

$$a * b \leq 0 \quad \textit{iff} \quad b \leq a^{\backprime} \quad \textit{iff} \quad a \leq {}'b. \tag{LP}$$

*Given these conditions, for each $a \in A$, $a^{\backprime}$ is the right relative pseudocomplement of $a$ with respect to* 0 *and $'a$ is the left relative pseudocomplement of $a$ with respect to* 0*. Hence, there are exactly two operations $^{\backprime}$ and $'$, which we will call* right pseudocomplementation *and* left pseudocomplemention of the operation $*$*, satisfying condition* (LP)*. This condition will be called the* Law of Pseudocomplementation*.*

*Proof. i) ⇒ ii)*: For each $c \in A$, let us consider the sets:

$$\mathrm{R}^0(c) = \{x \in A : c * x \leq 0\} \quad \text{and} \quad \mathrm{L}^0(c) = \{x \in A : x * c \leq 0\}.$$

By *i)*, we have that these sets both have a minimum element. Then, for each $c \in A$, let us define $c^{\backprime} := \max \mathrm{R}^0(c)$ and $'c := \max \mathrm{L}^0(c)$. Let $a, b \in A$ and suppose $a * b \leq 0$. Then, clearly, $b \leq a^{\backprime}$ and $a \leq {}'b$. However, if $b \leq a^{\backprime}$, by left monotonicity we have $a * b \leq a * a^{\backprime}$ but, since $a^{\backprime} \in \mathrm{R}^0(a)$, we have $a * b \leq 0$; if $a \leq {}'b$, by right monotonicity, we have $a * b \leq {}'b * b$ and, since $'b \in \mathrm{L}^0(b)$, we obtain $a * b \leq 0$.

*ii*) $\Rightarrow$ *i*): Let $a \in A$. We want to demonstrate that $a^{\backprime}$ is the maximum element of $\mathrm{R}^0(a)$. Observe that by *ii*), $a^{\backprime} \leq a^{\backprime}$ is equivalent to $a * a^{\backprime} \leq 0$ and so $a \in \mathrm{R}^0(a)$. Suppose now that $c \in \mathrm{R}^0(a)$: then $a * c \leq 0$ which, by *ii*), is equivalent to $c \leq a^{\backprime}$. Analogously we can prove that ${}'a$ is the maximum element of $\mathrm{L}^0(a)$. □

**Remark 13.** *If the operation $*$ is commutative, then the two pseudocomplements coincide and the corresponding operation receives the name of pseudocomplementation of $*$. In this case, we use the symbol $\neg$ and we write $\neg a$ instead of $a^{\backprime}$ or ${}'a$.*

**Definition 24** (Pseudocomplemented po-groupoids)**.** *A* pseudocomplemented po-groupoid *is an order-algebra $\mathcal{A} = \langle A, *, {}^{\backprime}, {}', 0, \leq \rangle$ with algebraic type $\langle 2, 1, 1, 0 \rangle$ such that $\langle A, *, 0, \leq \rangle$ is a pointed po-groupoid and the operations ${}^{\backprime}$ and ${}'$ are the left and right pseudocomplements of the operation $*$, respectively. By* pseudocomplemented *we refer to every structure that has a pseudocomplemented po-groupoid as its reduct.*

**Remark 14.** *Note that every $\mathbb{FL}$-algebra is a pseudocomplemented structure since its $\langle \leq, *, {}^{\backprime}, {}', 0, 1 \rangle$-reduct is a pseudocomplemented po-groupoid ($\leq$ is the order associated with the lattice).*

The following result is a reformulation of Corollary 12

**Proposition 24.** *Let $\mathcal{A}$ be a po-groupoid. The following conditions are equivalent:*

*i*) *For each $b \in A$, there are two unary operations ${}^{\backprime\mathcal{A}_b}$ and ${}^{\prime\mathcal{A}_b}$ defined on $A$ such that the structure $\mathcal{A}_b = \langle A, *, {}^{\backprime\mathcal{A}_b}, {}^{\prime\mathcal{A}_b}, b, \leq \rangle$ is a pseudocomplemented po-groupoid.*

*ii*) *The operation $*$ is residuated.*

*Under these conditions, for each $a, b \in A$, the right residual of $a$ relative to $b$ is the right relative pseudocomplement of $a$ with respect to $b$ and the left residual of $a$ relative to $b$ is the left relative pseudocomplement of $a$ with respect to $b$.*

*Proof.* *i*) $\Rightarrow$ *ii*): Let us define on $A$ two binary operations, $\backslash$ and $/$, in the following way: for each $a, b \in A$, $a \backslash b := a^{\backprime\mathcal{A}_b}$ and $b/a := {}^{\prime\mathcal{A}_b}a$. Since, for each $b \in A$, the operations $*$, ${}^{\backprime\mathcal{A}_b}$ and ${}^{\prime\mathcal{A}_b}$ satisfy the condition (LP), we have that, for each $a, b \in A$ the operations $*$, $\backslash$ and $/$ satisfy the condition (LR). Hence, $*$ is residuated.

*ii*) $\Rightarrow$ *i*): Let $\backslash$ and $/$ be the operations of left and right residuation corresponding to the operation $*$ and let $b \in A$. We define on $A$ two binary operations, ${}^{\backprime}$ and ${}'$, in the following way: for each $a \in A$, $a^{\backprime} := a \backslash b$ and ${}'a := b/a$. Then, by (LR), we have that the operations $*$, ${}^{\backprime}$ and ${}'$ satisfy condition (LP). Thus, $\mathcal{A}_b = \langle A, *, {}^{\backprime}, {}', b, \leq \rangle$ is a pseudocomplemented po-groupoid. □

**Definition 25** (Pseudocomplemented po-monoid)**.** *A* pseudocomplemented po-monoid *is an order-algebra $\mathcal{A} = \langle A, *, {}^{\backprime}, {}', 0, 1, \leq \rangle$ such that $\langle A, *, 1 \rangle$ is a monoid and $\langle A, *, {}^{\backprime}, {}', 0, \leq \rangle$ is a pseudocomplemented po-groupoid.*

**Definition 26** (Pseudocomplemented $s\ell$-monoid)**.** *A* pseudocomplemented $s\ell$-monoid *is an algebra $\mathbf{A} = \langle A, \vee, *, {}^{\backprime}, {}', 0, 1 \rangle$ such that $\langle A, \vee, *, 1 \rangle$ is an $s\ell$-monoid and $\langle A, *, {}^{\backprime}, {}', 0, \leq \rangle$, where $\leq$ is the semilattice order, is a pseudocomplemented po-groupoid.*

**Definition 27** (Pseudocomplemented $\ell$-monoid)**.** *A* pseudocomplemented $\ell$-monoid *is an algebra* $\mathbf{A} = \langle A, \vee, \wedge, *, ^{\backprime}, {}', 0, 1\rangle$ *such that* $\langle A, \vee, \wedge, *, 1\rangle$ *is an* $\ell$*-monoid and* $\langle A, *, ^{\backprime}, {}', 0, \leq\rangle$*, where* $\leq$ *is the order of the lattice, is a pseudocomplemented po-groupoid.*

**Definition 28** (The classes $\mathbb{PM}^{\preccurlyeq}_{\sigma}$, $\mathbb{PM}^{s\ell}_{\sigma}$, $\mathbb{PM}^{\ell}_{\sigma}$)**.** *We will denote by* $\mathbb{PM}^{\preccurlyeq}$*,* $\mathbb{PM}^{s\ell}$ *and* $\mathbb{PM}^{\ell}$*, the classes of pseudocomplemented po-monoids, pseudocomplemented s*$\ell$*-monoids, and pseudocomplemented* $\ell$*-monoids, respectively. Let* $\lambda \in \{\preccurlyeq, s\ell, \ell\}$ *and let* $\sigma$ *be a subsequence, possibly empty, of* $ew_lw_rc$*. We define the classes* $\mathbb{PM}^{\lambda}_{\sigma}$ *as the substructures of* $\mathbb{PM}^{\lambda}$ *that satisfy the properties in the set*

$$\{(e \preccurlyeq), (w \preccurlyeq), (\preccurlyeq w), (c \preccurlyeq)\}$$

*codified by the sequence* $\sigma$ *and, if* $\sigma$ *is the empty sequence, then* $\mathbb{PM}^{\lambda}_{\sigma}$ *is* $\mathbb{PM}^{\lambda}$*. The members of* $\mathbb{PM}^{s\ell}_{\sigma}$ ($\mathbb{PM}^{\ell}_{\sigma}$) *are called* $\mathbb{PM}^{s\ell}_{\sigma}$*-algebras* ($\mathbb{PM}^{\ell}_{\sigma}$*-algebras).*

**Definition 29** (Mirror image)**.** *Let* $t$ *be a term of an algebraic language,* $\mathcal{L}$*, such that* $\langle *, ^{\backprime}, {}'\rangle \leq \mathcal{L} \leq \langle \vee, \wedge, *, ^{\backprime}, {}', 0, 1\rangle$*. We define its* mirror image*,* $\mu(t)$*, inductively on the complexity of* $t$*:*

$$\mu(t) := \begin{cases} t, & \text{if } t \in Var \text{ or } t \in \{0,1\}, \\ \mu(u) \vee \mu(v), & \text{if } t = u \vee v, \\ \mu(u) \wedge \mu(v), & \text{if } t = u \wedge v, \\ \mu(v) * \mu(u), & \text{if } t = u * v, \\ \mu(u)^{\backprime}, & \text{if } t = u^{\backprime}, \\ {}'\mu(u), & \text{if } t = {}'u. \end{cases}$$

*We define the* mirror image of a formula *of the first-order language with equality* $\mathcal{L}^{\preccurlyeq} = \langle \mathcal{L}, \preccurlyeq\rangle$ *as the formula obtained by substituting all the terms that occur in the formula by their mirror images.*

**Lemma 7.** *We will denote by* $\mathbb{K}$ *the class of the pseudocomplemented po-groupoids (po-monoids, s*$\ell$*-monoids,* $\ell$*-monoids). Let us denote by* $\mathcal{L}_{\mathbb{K}}$ *its algebraic language and by* $\mathcal{L}^{\preccurlyeq}_{\mathbb{K}}$ *the first-order language with equality* $\langle \mathcal{L}_{\mathbb{K}}, \preccurlyeq\rangle$*. Let* $\mathcal{A} \in \mathbb{K}$*, and let* $\mathcal{A}^{o}$ *be the* $\mathcal{L}^{\preccurlyeq}_{\mathbb{K}}$*-structure defined in the following way:*

i) *the universe and the order of* $\mathcal{A}^{o}$ *are as in* $\mathcal{A}$*,*
ii) *the operations and distinguished elemnts of* $\mathcal{A}^{o}$ *belonging to* $\{\vee, \wedge, 0, 1\}$ *are the same as in* $\mathcal{A}$*, and*
iii) *the rest of the operations are defined in the following way: for each* $a, b \in A$*,*

$$a *^{\mathcal{A}^{o}} b := b *^{\mathcal{A}} a, \quad a^{\backprime^{\mathcal{A}^{o}}} := {}'^{\mathcal{A}}a, \quad {}'^{\mathcal{A}^{o}}a := a^{\backprime^{\mathcal{A}}}.$$

*Then,*

*i)* $\mathcal{A}^{o}$*, which will be called the* opposite *of* $\mathcal{A}$*, belongs to* $\mathbb{K}$*,*
*ii)* *for every term* $t$ *of* $\mathcal{L}_{\mathbb{K}}$*, it holds that* $\mu(t)^{\mathcal{A}^{o}} = t^{\mathcal{A}}$*.*

*Proof. i)*: It is easy to see that the operations $*^{\mathcal{A}^{o}}$, $^{\backprime^{\mathcal{A}^{o}}}$ and ${}'^{\mathcal{A}^{o}}$ satisfy (LP).
*ii)*: By induction on the complexity of $t$. □

**Lemma 8.** *Let $\mathcal{A}$ be a pseudocomplemented po-groupoid (po-monoid, s$\ell$-monoid, $\ell$-monoid) and let $\mathcal{A}^\partial$ be its opposite. Then, for every first-order formula $\varphi$ of the language of the pseudocomplemented po-groupoids (po-monoids, s$\ell$-monoids, $\ell$-monoids), the following holds:*

$$\mathcal{A} \models \varphi \quad \text{if and only if} \quad \mathcal{A}^\partial \models \mu(\varphi).$$

*Proof.* Let $\mathbb{K}$ and $\mathcal{L}_{\mathbb{K}}^{\preccurlyeq}$ be as in Lemma 7. Given $\mathcal{A} \in \mathbb{K}$ and an $\mathcal{L}_{\mathbb{K}}^{\preccurlyeq}$-formula, $\varphi$, it is easy to see by induction on the complexity of $\varphi$ that $\mathcal{A} \not\models \varphi$ if and only if $\mathcal{A}^\partial \not\models \mu(\varphi)$. □

**Theorem 13** (Law of Mirror Images)**.** *A formula is valid in the class of the pseudocomplemented po-groupoids (po-monoids, s$\ell$-monoids, $\ell$-monoids) if and only if its mirror image is also valid in that class.*

*Proof.* By Lemma 8. □

**Corollary 14.** *A quasi-inequation is valid in the class of the pseudocomplemented po-groupoids (po-monoids, s$\ell$-monoids, $\ell$-monoids) if and only if its mirror image is valid in that class.*

**Remark 15.** *Observe that every subclass of the classes considered defined by a set of formulas and their mirror images satisfies the Law of Mirror Images.*

**Proposition 25.** *In every pseudocomplemented po-groupoid $\mathcal{A}$ the following inequations and quasi-inequations are satisfied:*

$$\begin{array}{rlrl}
i) & x * x^{\backprime} \preccurlyeq 0, & i') & {}'x * x \preccurlyeq 0, \\
ii) & x \preccurlyeq y \supset x * y^{\backprime} \preccurlyeq 0, & ii') & x \preccurlyeq y \supset {}'y * x \preccurlyeq 0, \\
iii) & x \preccurlyeq y \supset y^{\backprime} \preccurlyeq x^{\backprime}, & iii') & x \preccurlyeq y \supset {}'y \preccurlyeq {}'x, \\
iv) & x \preccurlyeq {}'(x^{\backprime}), & iv') & x \preccurlyeq ({}'x)^{\backprime}, \\
v) & x^{\backprime} \approx ({}'(x^{\backprime}))^{\backprime}, & v') & {}'x \approx {}'(({}'x)^{\backprime}).
\end{array}$$

*Proof.* By the Law of Mirror Images, it will be sufficient to prove one of the two inequations or quasi-inequations in each row. Let $a, b \in A$.
*i)* By reflexivity we have $a^{\backprime} \leq a^{\backprime}$ and, by (LP), this is equivalent to $a * a^{\backprime} \leq 0$.
*ii)* Suppose $a \leq b$. By monotonicity and by *i)* we have $a * b^{\backprime} \leq b * b^{\backprime} \leq 0$.
*iii)* By *ii*), if $a \leq b$, then $a * b^{\backprime} \leq 0$ and by applying (LP), we obtain $b^{\backprime} \leq a^{\backprime}$.
*iv)* From $a * a^{\backprime} \leq 0$ we obtain $a \leq {}'(a^{\backprime})$ by applying (LP).
*v)* By *iv*) we have $a \leq {}'(a^{\backprime})$; from this by *iii*) we obtain $({}'(a^{\backprime}))^{\backprime} \leq a^{\backprime}$, which by *iv'*) gives $a^{\backprime} \leq ({}'(a^{\backprime}))^{\backprime}$. Hence, $a^{\backprime} = ({}'(a^{\backprime}))^{\backprime}$. □

Thus, by *iii*) and *iii'*) of the last proposition, we have that the left and right pseudocomplements are antimonotonous operations.

**Proposition 26.** *Let $\mathcal{A} = \langle A, *, {}^{\backprime}, {}', 0, \leq \rangle$ be a pseudocomplemented po-groupoid and let $\{a_i : i \in I\}$ be a family of elements in $A$. If $\bigvee_{i \in I} a_i$ exists, then there also exist $\bigwedge_{i \in I} a_i{}^{\backprime}$ and $\bigwedge_{i \in I} {}'a_i$ and the following holds:*

$$(\bigvee_{i \in I} a_i)^{\backprime} = \bigwedge_{i \in I} a_i{}^{\backprime} \;\; ; \; {}'(\bigvee_{i \in I} a_i) = \bigwedge_{i \in I} {}'a_i.$$

*Proof.* Since $a_i \leq \bigvee_{i\in I} a_i$, by the antimonotonicity of the right pseudocomplementation, we obtain $(\bigvee_{i\in I} a_i)^{\backslash} \leq a_i^{\backslash}$. Therefore, $(\bigvee_{i\in I} a_i)^{\backslash}$ is a lower bound of the set $\{a_i^{\backslash} : i \in I\}$. Let $b \in A$ and suppose that $b$ is also a lower bound of that set. By (LP), $b \leq a_i^{\backslash}$ is equivalent to $a_i * b \leq 0$, which is equivalent to $a_i \leq {}'b$. Therefore, $(\bigvee_{i\in I} a_i) \leq {}'b$, which is equivalent to $b \leq (\bigvee_{i\in I} a_i)^{\backslash}$. Consequently, $(\bigvee_{i\in I} a_i)^{\backslash}$ is the infimum of the set $\{a_i^{\backslash} : i \in I\}$. The other identity is analogously proved by using the antimonotonicity of the left pseudocomplementation and (LP). □

8.2. **Characterization of the $\mathbb{PM}^{\preccurlyeq}$ class.** In this section we present a set of inequalities which, together with the condition of antimonotonicity of the operations $^{\backslash}$ and $'$, characterize the law of pseudocomplementation (LP) in the class $\mathbb{PM}^{\preccurlyeq}$.

**Theorem 14.** *Let $\mathcal{A} = \langle A, *, ^{\backslash}, ', 0, 1, \leq\rangle$ be an order-algebra with algebraic type $\langle 2, 1, 1, 0, 0\rangle$ such that $\langle A, *, 1, \leq\rangle$ is a po-monoid. The following conditions are equivalent:*

*i) $\mathcal{A}$ is a pseudocomplemented po-monoid.*

*ii) The operations $^{\backslash}$ and $'$ are antimonotonous and $\mathcal{A}$ is a model of the inequations:*

$$r_1)\ \ 1^{\backslash} \preccurlyeq 0 \qquad r_2)\ \ 1 \preccurlyeq 0^{\backslash} \qquad r_3)\ \ x * (y * x)^{\backslash} \preccurlyeq y^{\backslash}$$
$$l_1)\ \ {}'1 \preccurlyeq 0 \qquad l_2)\ \ 1 \preccurlyeq {}'0 \qquad l_3)\ \ {}'(x * y) * x \preccurlyeq {}'y$$

*Proof.* *i) ⇒ ii)*: By Proposition 25, the operations $^{\backslash}$ and $'$ are antimonotonous. Now we will show that $(r_1)$, $(r_2)$ and $(r_3)$ are satisfied. This will be sufficient since $(l_1)$, $(l_2)$ and $(l_3)$ are their respective mirror images. We will use the fact that in every pseudocomplemented monoid $x * x^{\backslash} \preccurlyeq 0$ holds. We have that $1^{\backslash} \leq 1 * 1^{\backslash} \leq 0$. Thus, since $0 * 1 = 0 \leq 0$, by (LP) we obtain $1 \leq 0^{\backslash}$. Let $a, b \in A$. Then $b * (a * (b * a)^{\backslash}) \leq (b * a) * (b * a)^{\backslash} \leq 0$. Thus, by (LP), we can conclude $a * (b * a)^{\backslash} \leq b^{\backslash}$.
*ii) ⇒ i)*: Let $a, b \in A$. If $a * b \leq 0$, then, by the antimonotonicity of $^{\backslash}$, we have $0^{\backslash} \leq (a * b)^{\backslash}$ and therefore $1 \leq (a * b)^{\backslash}$; now we have $b = b * 1 \leq b * (a * b)^{\backslash} \leq a^{\backslash}$. Suppose now that $b \leq a^{\backslash}$. Then, $a * b \leq a * a^{\backslash} = a * (1 * a)^{\backslash} \leq 1^{\backslash} \leq 0$. To prove the equivalence between $a * b \leq 0$ and $a \leq {}'b$, we can proceed analogously. Therefore, $\mathcal{A}$ satisfies (LP). □

In the commutative case this result takes the form given in the following theorem.

**Theorem 15.** *Let $\mathcal{A} = \langle A, *, \neg, 0, 1, \leq\rangle$ be an order-algebra with algebraic type $\langle 2, 1, 1, 0, 0\rangle$ such that $\langle A, *, 1, \leq\rangle$ is a commutative* po*-monoid. The following conditions are equivalent:*

i) *$\mathcal{A}$ is a pseudocomplemented commutative po-monoid.*

ii) *The operation $\neg$ is antimonotonous and the following inequations are valid in $\mathcal{A}$:*

$$p_1)\ \ \neg 1 \preccurlyeq 0 \qquad p_2)\ \ 1 \preccurlyeq \neg 0 \qquad p_3)\ \ x * \neg(y * x) \preccurlyeq \neg y$$

As corollaries of the previous results we obtain the following characterizations for the classes $\mathbb{PM}^{\preccurlyeq}$ i $\mathbb{PM}^{\preccurlyeq}_e$.

**Corollary 15.** *An order-algebra $\mathcal{A} = \langle A, *, ^{\backslash}, ', 0, 1, \leq\rangle$ with algebraic type $\langle 2, 1, 1, 0, 0\rangle$ is a pseudocomplemented monoid if and only if:*

1) *$\langle A, *, 1, \leq\rangle$ is a po-monoid.*

2) *The operations* $^{\backprime}$ *and* $'$ *are antimonotonous with respect to the partial order.*
3) $\mathcal{A}$ *satisfies the inequations:*

$$r_1)\ \ 1^{\backprime} \preccurlyeq 0 \qquad r_2)\ \ 1 \preccurlyeq 0^{\backprime} \qquad r_3)\ \ x * (y * x)^{\backprime} \preccurlyeq y^{\backprime}$$
$$l_1)\ \ {}'1 \preccurlyeq 0 \qquad l_2)\ \ 1 \preccurlyeq {}'0 \qquad l_3)\ \ {}'(x * y) * x \preccurlyeq {}'y$$

**Corollary 16.** *An order-algebra* $\mathcal{A} = \langle A, *, \neg, 0, 1, \leq \rangle$ *with algebraic type* $\langle 2, 1, 0, 0 \rangle$ *is a commutative pseudocomplemented monoid if and only if:*

1) $\langle A, *, 1, \leq \rangle$ *is a commutative po-monoid.*
2) *The operation* $\neg$ *is antimonotonous with respect to the partial order.*
3) $\mathcal{A}$ *satisfies the inequations:*

$$p_1)\ \ \neg 1 \preccurlyeq 0 \qquad p_2)\ \ 1 \preccurlyeq \neg 0 \qquad p_3)\ \ x * \neg(y * x) \preccurlyeq \neg y$$

The results we present in the following two propositions give some equations and inequations that are satisfied in the classes $\mathbb{PM}^{\preccurlyeq}$ and $\mathbb{PM}_e^{\preccurlyeq}$.

**Proposition 27.** *In every pseudocomplemented po-monoid the following equations and inequations are satisfied:*

$$i)\ \ 0 \approx 1^{\backprime} \qquad\qquad i')\ \ 0 \approx {}'1$$
$$ii)\ \ x \preccurlyeq {}'(y * (x * y)^{\backprime}) \qquad ii')\ \ x \preccurlyeq ({}'(y * x) * y)^{\backprime}$$

*Proof.* It is sufficient to prove *i)* and *ii)*, since *i')* and *ii')*, are their respective mirror images.
*i)* We have $1 * 0 = 0 \leq 0$. Thus, by applying (LP), $0 \leq 1^{\backprime}$; this and $(r_1)$ imply $0 \approx 1^{\backprime}$.
*ii)* If $a, b \in A$, since $(a * b) * (a * b)^{\backprime} \leq 0$, by associativity, $a * (b * (a * b)^{\backprime}) \leq 0$ and from this, by applying (LP), $a \leq {}'(b * (a * b)^{\backprime})$. □

In the commutative case this result takes the form given in the following proposition.

**Proposition 28.** *In every commutative pseudocomplemented po-monoid the following conditions are satisfied:*

$$i)\ \ 0 \approx \neg 1 \qquad ii)\ \ x \preccurlyeq \neg(y * \neg(x * y))$$

8.3. **The $\mathbb{PM}^{s\ell}$ and $\mathbb{PM}^{\ell}$ classes are varieties.** In this section, we characterize the classes $\mathbb{PM}^{s\ell}$ and $\mathbb{PM}^{\ell}$ as equational classes.

**Lemma 9.** *Let* $\mathbf{A} = \langle A, \vee \rangle$ *be a semilattice and let* $\iota$ *be a unary operation defined in* $A$. *The following conditions are equivalent:*

i) *The operation* $\iota$ *is antimonotonous with respect to the order of the semilattice.*
ii) $\mathbf{A} \models \iota(x \vee y) \vee \iota(x) \approx \iota(x)$.

*Proof.* Suppose $a, b \in A$ and let $\leq$ be the semilattice order.
*i)* $\Rightarrow$ *ii)*: Since $a \leq a \vee b$, by the antimonotonicity of the operation $\iota$ we obtain $\iota(a \vee b) \leq \iota(a)$; i.e., $\iota(a \vee b) \vee \iota(a) = \iota(a)$.
*ii)* $\Rightarrow$ *i)*: Suppose $a \leq b$. Then $a \vee b = b$ and so, using *ii)*, we have $\iota(b) = \iota(a \vee b) \leq \iota(a)$. □

**Theorem 16.** $\mathbb{PM}^{s\ell}$ *is the equational class of algebras* $\mathbf{A} = \langle A, \vee, *, ^{\backprime}, ', 0, 1 \rangle$ *of type* $\langle 2, 2, 1, 1, 0, 0 \rangle$ *that satisfy:*

1) *Any set of equations defining the class of $s\ell$-monoids.*
2) *The equations (which characterize the Law of Pseudocomplementation):*
   $r_1$) $1^{\backprime} \approx 0$  $r_2$) $1 \vee 0^{\backprime} \approx 0^{\backprime}$  $r_3$) $(x * (y * x)^{\backprime}) \vee y^{\backprime} \approx y^{\backprime}$
   $l_1$) ${}'1 \approx 0$  $l_2$) $1 \vee {}'0 \approx {}'0$  $l_3$) $({}'(x * y) * x) \vee {}'y \approx {}'y$
3) *The equations (which characterize the antimonotonicity of pseudocomplementations):*
   $r_a$) $(x \vee y)^{\backprime} \vee x^{\backprime} \approx x^{\backprime}$
   $l_a$) ${}'(x \vee y) \vee {}'x \approx {}'x$

*Proof.* From the definition, the fact that $\mathbf{A}$ is a pseudocomplemented $s\ell$-monoid is equivalent to the facts that $\langle A, \vee, *, 1\rangle$ is an $s\ell$-monoid and that, if $\leq$ is the order of the semilattice, then $\mathcal{A} = \langle A, *, {}^{\backprime}, {}', 0, 1, \leq\rangle$ is a pseudocomplemented po-monoid by the characterization of Theorem 14. This is equivalent to saying that the operations ${}^{\backprime}$ and ${}'$ are antimonotonous and that in $\mathcal{A}$ the inequations $(r_1), \ldots, (l_3)$ of Theorem 14 are satisfied. Observe that by substituting, given two terms $t_1$ and $t_2$, the inequalities $t_1 \preccurlyeq t_2$ by the equations $t_1 \vee t_2 \approx t_2$ (bearing in mind that $(r_1)$ and $(l_1)$ can be substituted by the equations $1^{\backprime} \approx 0$ and ${}'1 \approx 0$) we obtain the equations $(r_1), \ldots, (l_3)$ of the present theorem. Moreover, by Lemma 9, the antimonotonicity of the operations ${}^{\backprime}$ and ${}'$ is equivalent to the fact that the equations $(r_a)$ and $(l_a)$ are satisfied in $\mathbf{A}$. □

Analogously, we have the following characterization of $\mathbb{PM}^{\ell}$ as an equational class.

**Theorem 17.** *$\mathbb{PM}^{\ell}$ is the equational class of algebras $\mathbf{A} = \langle A, \vee, \wedge, *, {}^{\backprime}, {}', 0, 1\rangle$ of type $\langle 2, 2, 2, 1, 1, 0, 0\rangle$ that satisfy any set of equations defining the class of the $\ell$-monoids and the equations $(r_1), \ldots, (l_3)$, $(r_a)$ and $(l_a)$ of Theorem 16.*

**Corollary 17.** *The classes $\mathbb{PM}^{s\ell}_{\sigma}$ and $\mathbb{PM}^{\ell}_{\sigma}$, with $\sigma \leq ew_lw_rc$, are varieties.*

*Proof.* Any class $\mathbb{PM}^{s\ell}_{\sigma}$ ($\mathbb{PM}^{\ell}_{\sigma}$) is obtained by adding some of the equations

$$x * y \approx y * x,\ x \vee 1 \approx 1,\ 0 \vee x \approx x,\ x \vee (x * x) \approx x * x,$$

to the set of equations characterizing the class $\mathbb{PM}^{s\ell}$ ($\mathbb{PM}^{\ell}$). □

Next, by adapting the notation, we give an equational characterization of the commutative classes $\mathbb{PM}^{s\ell}_e$ and $\mathbb{PM}^{\ell}_e$.

**Theorem 18.** *$\mathbb{PM}^{s\ell}_e$ is the equational class of the algebras $\mathbf{A} = \langle A, \vee, *, \neg, 0, 1\rangle$ of type $\langle 2, 2, 1, 0, 0\rangle$ that satisfy:*

1) *Any set of equations defining the class of the commutative $s\ell$-monoids,*
2) $p_1$) $\neg 1 \approx 0$, $p_2$) $1 \vee \neg 0 \approx \neg 0$, $p_3$) $(x * \neg(y * x)) \vee \neg y \approx \neg y$,
   $a$) $\neg(x \vee y) \vee \neg x \approx \neg x$.

**Theorem 19.** *$\mathbb{PM}^{\ell}_e$ is the equational class of the algebras $\mathbf{A} = \langle A, \vee, \wedge, *, \neg, 0, 1\rangle$ of type $\langle 2, 2, 2, 1, 0, 0\rangle$ that satisfy any set of equations defining the commutative $\ell$-monoids and the equations $(p_1)$, $(p_2)$, $(p_3)$ and $(a)$ of Theorem 18.*

8.4. **Pseudocomplementation with respect to the minimum.** In this section we show that in the framework of pseudocomplemented po-monoids, when the distinguished element 0 is the minimum with respect to the partial order, condition (LP) is equivalent to the inequalities $(r_1)$, $(r_2)$, $(r_3)$, $(l_1)$, $(l_2)$ and $(l_3)$ without the need to add to these inequalities the condition of antimonotonicity of the pseudocomplementations. We also give an alternative set of inequalities which, together with the antimonotonicity of pseudocomplementations, characterizes (LP) in this class of pseudocomplemented po-monoids. We then also analyze the case in which the structures of this kind are integral.

**Theorem 20.** *Let $\mathcal{A} = \langle A, *, ^{\backprime}, ', 0, 1, \leq\rangle$ be an order-algebra with algebraic type $\langle 2, 1, 1, 0, 0\rangle$ such that $\langle A, *, 1, \leq\rangle$ is a* po*-monoid. Let us suppose that* 0 *is the minimum element with respect to the partial order. Then, the following conditions are equivalent:*

*i) $\mathcal{A}$ is a pseudocomplemented po-monoid.*

*ii) $\mathcal{A}$ satisfies the inequations:* [14]

$$r_1)\quad 1^{\backprime} \preccurlyeq 0 \qquad r_2)\quad 1 \preccurlyeq 0^{\backprime} \qquad r_3)\quad x * (y * x)^{\backprime} \preccurlyeq y^{\backprime}$$
$$l_1)\quad {}'1 \preccurlyeq 0 \qquad l_2)\quad 1 \preccurlyeq {}'0 \qquad l_3)\quad {}'(x * y) * x \preccurlyeq {}'y$$

*iii) The operations $^{\backprime}$ and $'$ are antimonotonous and $\mathcal{A}$ satisfies the inequations:*

$$r_4)\quad x * x^{\backprime} \preccurlyeq 0 \qquad r_5)\quad x \preccurlyeq {}'(y * (x * y)^{\backprime})$$
$$l_4)\quad {}'x * x \preccurlyeq 0 \qquad l_5)\quad x \preccurlyeq ({}'(y * x) * y)^{\backprime}$$

*Proof.* *i)* $\Leftrightarrow$ *ii)*: As we have seen in Theorem 14, in every pseudocomplemented po-monoid the inequalities $(r_1), \ldots, (l_3)$ are satisfied and, reciprocally, their validity and the antimonotonicity of the operations $^{\backprime}$ and $'$ allow us to prove (LP). Therefore, to prove this equivalence, it is sufficient to show that when $(r_1), \ldots, (l_3)$ hold in $\mathcal{A}$, the operations $^{\backprime}$ and $'$ are antimonotonous.

First, observe that the validity of $(r_1)$ and $(r_3)$ in $\mathcal{A}$ allows us to prove that $x * x^{\backprime} \preccurlyeq 0$ is also valid in $\mathcal{A}$: if $a \in A$, then $a * a^{\backprime} = a * (1 * a)^{\backprime} \leq 1^{\backprime} \leq 0$. Suppose $a, b \in A$ and $a \leq b$. We have $a * b^{\backprime} \leq b * b^{\backprime} \leq 0$ and, since 0 is the minimum, $a * b^{\backprime} = 0$. Now, by using $(r_2)$ and $(r_3)$, we obtain $b^{\backprime} = b^{\backprime} * 1 \leq b^{\backprime} * 0^{\backprime} = b^{\backprime} * (a * b^{\backprime})^{\backprime} \leq a^{\backprime}$. We can now proceed analogously to show that $a \leq b$ implies ${}'b \leq {}'a$, by using $(l_1)$, $(l_2)$ and $(l_3)$.

*i)* $\Rightarrow$ *iii)*: In every pseudocomplemented po-monoid the operations $^{\backprime}$ and $'$ are antimonotonous and the inequations $(r_4)$, $(l_4)$ (see Proposition 25), $(r_5)$ and $(l_5)$ (see Proposition 27) hold.

*iii)* $\Rightarrow$ *ii)*: By using $(r_4)$ we obtain: $1^{\backprime} = 1 * 1^{\backprime} \leq 0$. Since 0 is the minimum, we have $0 \leq {}'0 * 0 = {}'(1 * 0) * 0$ and, by monotonicity, $({}'(1 * 0) * 0)^{\backprime} \leq 0^{\backprime}$. Now, by $(l_5)$, we obtain $1 \leq ({}'(1 * 0) * 0)^{\backprime}$ and thus $1 \leq 0^{\backprime}$. By using $(l_5)$, we prove that $x \preccurlyeq ({}'x)^{\backprime}$ holds in $\mathcal{A}$: if $a \in A$, then $a \leq ({}'(1 * a) * 1)^{\backprime} = ({}'a * 1)^{\backprime} = ({}'a)^{\backprime}$. Finally, let $a, b \in A$. By $(l_4)$ we obtain $b \leq {}'(a * (b * a)^{\backprime})$, and by applying monotonicity we have $({}'(a * (b * a)^{\backprime}))^{\backprime} \leq b^{\backprime}$. But $a * (b * a)^{\backprime} \leq ({}'(a * (b * a)^{\backprime}))^{\backprime}$. Therefore, $a * (b * a)^{\backprime} \leq b^{\backprime}$. The inequalities $(l_1)$, $(l_2)$ and $(l_3)$ are proved in an analogous way. □

As immediate consequences of this result, we obtain the following characterizations for the class $\mathbb{PM}^{\preccurlyeq}_{w_r}$ and the varieties $\mathbb{PM}^{s\ell}_{w_r}$ and $\mathbb{PM}^{\ell}_{w_r}$.

[14] The inequations $r_1$) and $r_2$) can be substituted by the equations $1^{\backprime} \approx 0$ and ${}'1 \approx 0$.

**Corollary 18.** *An order-algebra* $\mathcal{A} = \langle A, *, ^{\backprime}, ', 0, 1, \leq\rangle$ *with algebraic type* $\langle 2, 1, 1, 0, 0\rangle$ *belongs to the class* $\mathbb{PM}^{\preccurlyeq}_{w_r}$ *if and only if the following conditions are satisfied:*

a) $\langle A, \leq, *, 1\rangle$ *is a po-monoid,*
b) $\mathcal{A} \models 0 \preccurlyeq x$,
c) $(r_1), \ldots, (l_3)$ *hold in* $\mathcal{A}$.

**Corollary 19.** $\mathbb{PM}^{s\ell}_{w_r}$ $(\mathbb{PM}^{\ell}_{w_r})$ *is the equational class of algebras*

$$\mathbf{A} = \langle A, \vee, *, ^{\backprime}, ', 0, 1\rangle \quad (\mathbf{A} = \langle A, \vee, \wedge, *, ^{\backprime}, ', 0, 1\rangle)$$

*of type* $\langle 2, 2, 1, 1, 0, 0\rangle$ $(\langle 2, 2, 2, 1, 1, 0, 0\rangle)$ *that satisfies:* [15]

1) *Any set of equations defining the class of the sℓ-monoids (ℓ-monoids).*
2) *The equation* $0 \preccurlyeq x$.
3) *The equations:*

$$\begin{array}{llllll} r_1) & 1^{\backprime} \preccurlyeq 0 & r_2) & 1 \preccurlyeq 0^{\backprime} & r_3) & x * (y * x)^{\backprime} \preccurlyeq y^{\backprime} \\ l_1) & '1 \preccurlyeq 0 & l_2) & 1 \preccurlyeq '0 & l_3) & '(x * y) * x \preccurlyeq 'y \end{array}$$

*The set of equations (3) can be substituted in the axiomatization of* $\mathbb{PM}^{s\ell}_{w_r}$ $(\mathbb{PM}^{\ell}_{w_r})$ *by the set of equations:*

$$\begin{array}{llllll} r_4) & x * x^{\backprime} \preccurlyeq 0 & r_5) & x \preccurlyeq '(y * (x * y)^{\backprime}) & r_6) & (x \vee y)^{\backprime} \preccurlyeq x^{\backprime} \\ l_4) & 'x * x \preccurlyeq 0 & l_5) & x \preccurlyeq ('(y * x) * y)^{\backprime} & l_6) & '(x \vee y) \preccurlyeq 'x \end{array}$$

*The equations* $(r_1)$ *and* $(l_1)$ *can be substituted by the equations* $1^{\backprime} \approx 0$ *and* $'1 \approx 0$, *respectively. The equations* $(r_4)$ *and* $(l_4)$ *can be substituted by* $x * x^{\backprime} \approx 0$ *and* $'x * x \approx 0$, *respectively.*

*Proof.* Immediate by using the characterization of Corollary 18 and the fact that the quasi-equations expressing the monotonicity can be substituted by the equations $(r_6)$ and $(l_6)$. □

**Proposition 29.** *Let* $\mathcal{A} = \langle A, *, ^{\backprime}, ', 0, 1, \leq\rangle$ *be a pseudocomplemented po-monoid with* 0 *as the minimum element. Then A has a maximum element, say* $\top$, *and* $'0 = 0^{\backprime} = \top$ *holds.*

*Proof.* Since 0 is the minimum, we have that, for every $a \in A$, $0 \leq 'a$ and $0 \leq a^{\backprime}$, that is, $0 * a \leq 0$ and $a * 0 \leq 0$. From this, we obtain $a \leq 0^{\backprime}$ and $a \leq '0$. In particular, $'0 \leq 0^{\backprime}$ and $0^{\backprime} \leq '0$. Therefore, $'0$ and $0^{\backprime}$ are the same element and this element is the maximum. □

In the following proposition we prove some properties of the pseudocomplemented monoids with 0 as the minimum element which are integral (i.e., 1 is the maximum element).

**Proposition 30.** *In every* $\mathcal{A} \in \mathbb{PM}^{\preccurlyeq}_{w}$, *the following are satisfied:*

*i)* $x * y \preccurlyeq x$ *ii)* $x * y \preccurlyeq y$ *iii)* $x * 0 \approx 0 * x \approx 0$ *iv)* $0^{\backprime} \approx 1 \approx '0$

*Proof.* The reduct $\langle A, *, 0, 1, \leq\rangle$ is a po-monoid with minimum element 0 and integral and thus (see Propositions 6 and 9) *i), ii)* and *iii)* are valid in it. Property *iv)* is a consequence of $(r_2)$, $(l_2)$ and the fact that 1 is the maximum. □

**Corollary 20.** *An order-algebra* $\mathcal{A} = \langle A, *, ^{\backprime}, ', 0, 1, \leq\rangle$ *of algebraic type* $\langle 2, 1, 1, 0, 0\rangle$ *is of the class* $\mathbb{PM}^{\preccurlyeq}_{w}$ *if and only if:*

---

[15] In the context of semilatticed algebras we use $t_1 \preccurlyeq t_2$ as an abbreviation for the equation $t_1 \vee t_2 \approx t_2$.

a) $\langle A, \leq, *, 1\rangle$ *is a po-monoid,*
b) *in* $\mathcal{A}$ *the inequations* $0 \preccurlyeq x$ *and* $x \preccurlyeq 1$ *are satisfied,*
c) *in* $\mathcal{A}$ *the following are satisfied:*

$$r_1)\ \ 1^{\backprime} \approx 0 \qquad r_2)\ \ 1 \approx 0^{\backprime} \qquad r_3)\ \ x * (y * x)^{\backprime} \preccurlyeq y^{\backprime}$$
$$l_1)\ \ {}'1 \approx 0 \qquad l_2)\ \ 1 \approx {}'0 \qquad l_3)\ \ {}'(x * y) * x \preccurlyeq {}'y$$

A characterization of the $\mathbb{PM}^{\preccurlyeq}_{ew}$ class is clearly obtained by adding the condition of the commutativity of the monoidal operation and substituting in Corollary 20 conditions *c)* with the following ones:

$$p_1)\ \ \neg 1 \approx 0 \qquad p_2)\ \ 1 \approx \neg 0 \qquad p_3)\ \ x * \neg(y * x) \preccurlyeq \neg y$$

The following result is an equational characterization of the $\mathbb{PM}^{s\ell}_{w}$ and $\mathbb{PM}^{\ell}_{w}$ classes.

**Corollary 21.** $\mathbb{PM}^{s\ell}_{w}$ ($\mathbb{PM}^{\ell}_{w}$) *is the equational class of the algebras*

$$\mathbf{A} = \langle A, \vee, *, {}^{\backprime}, {}', 0, 1\rangle \quad (\mathbf{A} = \langle A, \vee, \wedge, *, {}^{\backprime}, {}', 0, 1\rangle)$$

*of type* $\langle 2, 2, 1, 1, 0, 0\rangle$ ($\langle 2, 2, 2, 1, 1, 0, 0\rangle$) *that satisfy:*

1) *Any set of equations defining the class of the* $s\ell$*-monoids (*$\ell$*-monoids),*
2) $0 \preccurlyeq x,\ x \preccurlyeq 1$
3) $r_1)\ \ 1^{\backprime} \approx 0 \qquad r_2)\ \ 1 \approx 0^{\backprime} \qquad r_3)\ \ (x * (y * x)^{\backprime}) \preccurlyeq y^{\backprime}.$
   $l_1)\ \ {}'1 \approx 0 \qquad l_2)\ \ 1 \approx {}'0 \qquad l_3)\ \ ({}'(x * y) * x) \preccurlyeq {}'y.$

*Proof.* Immediate by using the characterization of Corollary 20. □

In the following statements we adapt the notation and give equational characterizations for the commutative classes $\mathbb{PM}^{s\ell}_{ew_r}$, $\mathbb{PM}^{\ell}_{ew_r}$, $\mathbb{PM}^{s\ell}_{ew}$ and $\mathbb{PM}^{\ell}_{ew}$.

**Corollary 22.** $\mathbb{PM}^{s\ell}_{ew_r}$ ($\mathbb{PM}^{\ell}_{ew_r}$) *is the equational class of the algebras*

$$\mathbf{A} = \langle A, \vee, *, \neg, 0, 1\rangle \quad (\mathbf{A} = \langle A, \vee, \wedge, *, \neg, 0, 1\rangle)$$

*of type* $\langle 2, 2, 1, 0, 0\rangle$ ($\langle 2, 2, 2, 1, 0, 0\rangle$) *that satisfy:*

1) *Any set of equations defining the class of the commutative* $s\ell$*-monoids (*$\ell$*-monoids),*
2) $0 \preccurlyeq x,$
3) $p_1)\ \ \neg 1 \approx 0 \qquad p_2)\ \ 1 \preccurlyeq \neg 0 \qquad p_3)\ \ (x * \neg(y * x)) \preccurlyeq \neg y.$

**Corollary 23.** (Cf. (Bou et al., 2006, Theorem 4.6)) $\mathbb{PM}^{s\ell}_{ew}$ ($\mathbb{PM}^{\ell}_{ew}$) *is the equational class of the algebras*

$$\mathbf{A} = \langle A, \vee, *, \neg, 0, 1\rangle \quad (\mathbf{A} = \langle A, \vee, \wedge, *, \neg, 0, 1\rangle)$$

*of type* $\langle 2, 2, 1, 0, 0\rangle$ ($\langle 2, 2, 2, 1, 0, 0\rangle$) *satisfying:*

1) *Any set of equations defining the class of the commutative* $s\ell$*-monoids (*$\ell$*-monoids),*
2) $0 \preccurlyeq x,\ \ x \preccurlyeq 1,$
3) $p_1)\ \ \neg 1 \approx 0 \qquad p_2)\ \ 1 \approx \neg 0 \qquad p_3)\ \ (x * \neg(y * x)) \preccurlyeq \neg y.$

**Remark 16.** *Let us stress that the variety* $\mathbb{PM}^{s\ell}_{ewc}$ *is precisely the variety of the* pseudo-complemented distributive lattices. *For an equational presentation of this class, see Balbes & Dwinger (1974). In the following statement we give an equational characterization of this variety which is an alternative to the one in Balbes & Dwinger (1974).*

**Corollary 24.** $\mathbb{PM}^{s\ell}_{ewc}$ *is the equational class of the algebras*

$$\mathbf{A} = \langle A, \vee, \wedge, \neg, 0, 1\rangle$$

*of type* $\langle 2, 2, 1, 0, 0\rangle$ *that satisfy:*

1) *Any set of equations defining the class of the distributive lattices,*
2) $0 \preccurlyeq x,\ x \preccurlyeq 1,$
3) $p_1)\ \ \neg 1 \approx 0 \qquad p_2)\ \ 1 \approx \neg 0 \qquad p_3)\ \ (x \wedge \neg(y \wedge x)) \preccurlyeq \neg y.$

*Proof.* Equations (1) and (2) in Theorem 23 define the $\mathring{\mathbb{M}}^{s\ell}_{ew}$ class. If we add the equation $x \preccurlyeq x * x$, then: the operation $*$ is equal to the operation $\wedge$ (see Proposition 12); the class $\mathring{\mathbb{M}}^{s\ell}_{ewc}$ is that of the distributive lattices; and the pseudocomplementation is with respect to the operation $\wedge$. □

## Part 3. CONNECTING GENTZEN SYSTEMS AND ALGEBRAS

This part has just one single section. In it we prove that the subsystems $\mathcal{FL}_\sigma[\vee, *, 0, 1]$, $\mathcal{FL}_\sigma[\vee, \wedge, *, 0, 1]$, $\mathcal{FL}_\sigma[\vee, *, `, ', 0, 1]$, and $\mathcal{FL}_\sigma[\vee, \wedge, *, `, ', 0, 1]$ are algebraizable, having as their respective equivalent algebraic semantics the varieties: $\mathring{\mathbb{M}}^{s\ell}_\sigma$, $\mathring{\mathbb{M}}^{\ell}_\sigma$, $\mathbb{PM}^{s\ell}_\sigma$ and $\mathbb{PM}^{\ell}_\sigma$. We also prove that the $\mathcal{FL}_\sigma$ system is algebraizable with its equivalent algebraic semantics being the $\mathbb{FL}_\sigma$ variety.

### 9. Algebraic Analysis of some Implication-Free Subsystems.

We will use the letter $\Psi$ as a generic denotation of the languages of the classes $\mathring{\mathbb{M}}^{s\ell}_\sigma$, $\mathring{\mathbb{M}}^{\ell}_\sigma$, $\mathbb{PM}^{s\ell}_\sigma$, $\mathbb{PM}^{\ell}_\sigma$ and $\mathbb{FL}_\sigma$ and we will use $\mathbb{K}_\sigma[\Psi]$ as a generic denotation for all these classes of algebras. In what follows, we will show that every subsystem $\mathcal{FL}_\sigma[\Psi]$ is algebraizable and that the class $\mathbb{K}_\sigma[\Psi]$ is its equivalent quasivariety semantics. To prove the algebraization results, we use Lemma 1. To this end, we need a translation $\tau$ from sequents to equations and a translation $\rho$ from equations to sequents. So we start by stating the definitions of these translations.

**Definition 30.** *We define the translations* $\tau$ *from* $\Psi$*-sequents to* $\Psi$*-equations and* $\rho$ *from* $\Psi$*-equations to* $\Psi$*-sequents in the following way:*

$$\tau(\varphi_0, ..., \varphi_{m-1} \Rightarrow \varphi) := \begin{cases} \{(\varphi_0 * ... * \varphi_{m-1}) \vee \varphi \approx \varphi\}, & \text{if } m \geq 1, \\ \{1 \vee \varphi \approx \varphi\}, & \text{if } m = 0, \end{cases}$$

$$\tau(\varphi_0, ..., \varphi_{m-1} \Rightarrow \emptyset) := \begin{cases} \{\varphi_0 * ... * \varphi_{m-1} \vee 0 \approx 0\}, & \text{if } m \geq 1, \\ \{1 \vee 0 \approx 0\}, & \text{if } m = 0, \end{cases}$$

$$\rho(\varphi \approx \psi) := \{\varphi \Rightarrow \psi,\ \psi \Rightarrow \varphi\}.$$

Note that the translation $\tau$ is well defined since the languages $\Psi$ contain all the connectives in $\{\vee, *, 0, 1\}$.

Now we are going to prove condition 1) of Lemma 1. We need some previous results and some notation.

**Lemma 10.** *For every $\mathcal{FL}_\sigma[\Psi]$-theory, $\Phi$,*

$$\rho(\varphi \preccurlyeq \psi) \subseteq \Phi \quad \textit{iff} \quad \varphi \Rightarrow \psi \in \Phi.$$

*Thus, in particular, the derivability of the sequents in $\rho(\varphi \preccurlyeq \psi)$ is equivalent to the derivability of the sequent $\varphi \Rightarrow \psi$.*

*Proof.* We have that $\rho(\varphi \preccurlyeq \psi) = \rho(\varphi \vee \psi \approx \psi) = \{\varphi \vee \psi \Rightarrow \psi, \psi \Rightarrow \varphi \vee \psi\}$. The sequent $\psi \Rightarrow \varphi \vee \psi$ is derivable from $\psi \Rightarrow \psi$ using $(\Rightarrow \vee_2)$. Thus, it is sufficient to prove that the sequents $\varphi \vee \psi \Rightarrow \psi$ and $\varphi \Rightarrow \psi$ are interderivable in $\mathbf{FL}_\sigma[\Psi]$. Let us consider the following formal proofs:

$$\dfrac{\dfrac{\varphi \Rightarrow \varphi}{\varphi \Rightarrow \varphi \vee \psi}\,(\Rightarrow \vee_1) \qquad \varphi \vee \psi \Rightarrow \psi}{\varphi \Rightarrow \psi}\,(Cut) \qquad \dfrac{\varphi \Rightarrow \varphi \qquad \psi \Rightarrow \psi}{\varphi \vee \psi \Rightarrow \psi}\,(\vee \Rightarrow)$$

□

**Notation 1.** *If $\bar{x} = x_0, \ldots, x_{m-1}$ is a sequence of elements in a $\Psi$-algebra, $\mathbf{A}$, then we define $\prod \bar{x} := 1$ if $\bar{x}$ is the empty sequence, $\prod \bar{x} := x_0$ if $m = 1$, and $\prod \bar{x} := x_0 * \ldots * x_{m-1}$ if $m \geq 1$. In particular, if $\mathbf{A}$ is the algebra of $\Psi$-formulas, for each sequence $\Gamma = \varphi_0, \ldots, \varphi_{m-1}$, we will have $\prod \Gamma := 1$ if $\Gamma$ is the empty sequence, $\prod \Gamma := \varphi_0$ if $m = 1$ and $\prod \Gamma := \varphi_0 * \ldots * \varphi_{m-1}$ if $m \geq 1$.*

**Lemma 11.** *If $\Gamma$ is a sequence of $\Psi$-formulas, then the sequent $\Gamma \Rightarrow \prod \Gamma$ is derivable in $\mathbf{FL}_\sigma[\Psi]$.*

*Proof.* By induction on the length of the sequence $\Gamma$.
- If $m = 0$, then $\Gamma \Rightarrow \prod \Gamma$ is the sequent $\emptyset \Rightarrow 1$, that is, the axiom $(\Rightarrow 1)$.
- If $m > 0$ and $\Gamma = \varphi_0, \ldots, \varphi_{m-1}$, by the induction hypothesis we have that $\varphi_0, \ldots, \varphi_{m-2} \Rightarrow \varphi_0 * \ldots * \varphi_{m-2}$ is derivable. From this sequent and $\varphi_{m-1} \Rightarrow \varphi_{m-1}$, applying $(\Rightarrow *)$, we obtain $\Gamma \Rightarrow \prod \Gamma$. □

In the next lemma we already prove that condition 1) of Lemma 1 is satisfied.

**Lemma 12.** *For every $\varsigma \in Seq_\Psi^{\omega \times \{0,1\}}$, $\varsigma \dashv\vdash_{\mathbf{FL}_\sigma[\Psi]} \rho\tau(\varsigma)$.*

*Proof.* We consider two cases: *a)* $\varsigma = \Gamma \Rightarrow \varphi$ and *b)* $\varsigma = \Gamma \Rightarrow \emptyset$.
*a)*: By the definition of $\tau$ we have $\rho\tau(\varsigma) = \rho(\prod \Gamma \preccurlyeq \varphi)$. But, by Lemma 11, we have $\rho(\prod \Gamma \preccurlyeq \varphi) \dashv\vdash \prod \Gamma \Rightarrow \varphi$. Thus, it will be sufficient to prove $\Gamma \Rightarrow \varphi \dashv\vdash \prod \Gamma \Rightarrow \varphi$. The sequent $\Gamma \Rightarrow \varphi$ is obtained from the derivable sequent $\Gamma \Rightarrow \prod \Gamma$ (see Lemma 11) and the sequent $\prod \Gamma \Rightarrow \varphi$ by applying $(Cut)$.
We will show that $\prod \Gamma \Rightarrow \varphi$ is obtained from $\Gamma \Rightarrow \varphi$ using induction on the length of the sequence $\Gamma$. If $n = 0$, we must see that $\emptyset \Rightarrow \varphi \vdash 1 \Rightarrow \varphi$ and this is clear by applying $(1 \Rightarrow)$. If $n > 0$, we have the following derivation:

$$\dfrac{\dfrac{\Gamma \Rightarrow \varphi}{\varphi_0, \ldots, \varphi_{m-2} * \varphi_{m-1} \Rightarrow \varphi}\,(* \Rightarrow)}{\prod \Gamma \Rightarrow \varphi}\,(\text{Induction hypothesis})$$

$b)$: In this case we have $\rho\tau(\varsigma) = \rho(\prod\Gamma \preccurlyeq 0)$. Thus, we must prove

$$\Gamma \Rightarrow \emptyset \dashv\vdash \prod\Gamma \Rightarrow 0.$$

The sequent $\Gamma \Rightarrow \emptyset$ can be obtained from $\prod\Gamma \Rightarrow 0$ in the following way:

$$\dfrac{\dfrac{\Gamma \Rightarrow \prod\Gamma \qquad \prod\Gamma \Rightarrow 0}{\Gamma \Rightarrow 0}\,(\text{Tall}) \qquad 0 \Rightarrow \emptyset}{\Gamma \Rightarrow \emptyset}\,(Cut)$$

To see that $\prod\Gamma \Rightarrow 0$ can be obtained from $\Gamma \Rightarrow \emptyset$ we use induction on the length of the sequence $\Gamma$. If $n = 0$, we must see that $\emptyset \Rightarrow \emptyset \vdash 1 \Rightarrow 0$ and this is immediate using $(1 \Rightarrow)$ and $(\Rightarrow 0)$. If $n > 0$ we have the following derivation:

$$\dfrac{\dfrac{\Gamma \Rightarrow \emptyset}{\varphi_0, \dots, \varphi_{m-2} * \varphi_{m-1} \Rightarrow \emptyset}\,(* \Rightarrow)}{\prod\Gamma \Rightarrow 0}\,(\text{Induction hypothesis})$$

□

In the following lemma we prove that condition 2) of Lemma 1 is satisfied.

**Lemma 13.** *For every* $\varphi \approx \psi \in Eq_\Psi$, $\varphi \approx \psi \dashv\models_{\mathbb{K}_\sigma[\Psi]} \tau\rho(\varphi \approx \psi)$.

*Proof.* We must prove that $\varphi \approx \psi \dashv\models_{\mathbb{K}_\sigma[\Psi]} \{\varphi \vee \psi \approx \psi,\ \psi \vee \varphi \approx \varphi\}$; which is trivial. □

In the following lemma we prove that condition 3) of Lemma 1 is satisfied.

**Lemma 14.** *For every* $\mathbf{A} \in \mathbb{K}_\sigma[\Psi]$ *we define* $R$ *as the set*

$$\{\langle \bar{x}, \bar{y}\rangle \in A^m \times A^n : \langle m, n\rangle \in \omega \times \{0,1\},\ \mathbf{A} \models \tau(p_0, \dots, p_{m-1} \Rightarrow q_0, \dots, q_{n-1})[\![\bar{x}, \bar{y}]\!]\}.$$

*Then* $R$ *is an* $\mathcal{FL}_\sigma[\Psi]$*-filter.*

*Proof.* Given one of the languages $\Psi$, we must prove that, for every $\sigma$ and every $\mathbf{A} \in \mathbb{K}_\sigma[\Psi]$, the set $R$ contains the interpretations of the axioms of $\mathbf{FL}_\sigma[\Psi]$ and is closed under the interpretation of the rules of $\mathbf{FL}_\sigma[\Psi]$. First, observe that the set $R$ is equal to:

$$\{\langle \bar{x}, a\rangle \in A^m \times A : m \in \omega, \prod\bar{x} \leq a\} \cup \{\langle \bar{x}, \emptyset\rangle \in A^m \times \{\emptyset\} : m \in \omega, \prod\bar{x} \leq 0\}.$$

For every algebra $\mathbf{A}$ in the classes considered, $R$ contains all the pairs of the form $\langle a, a\rangle$ (where $a \in A$), $\langle 0, \emptyset\rangle$ and $\langle \emptyset, 1\rangle$; that is, $R$ contains the interpretations of the axioms. Next we will see that in each case $R$ is closed under the interpretation of the rules. We begin with the rules common to all the calculi under consideration.

From now on we will use the symbol $\delta$ to denote the empty set or an arbitrary element of $A$. Then $c_\delta \in A$ is defined as 0 if $\delta = \emptyset$ and as $\delta$ if $\delta \in A$.

- *Cut rule*:

$$\dfrac{\Gamma \Rightarrow \varphi \qquad \Sigma, \varphi, \Pi \Rightarrow \Delta}{\Sigma, \Gamma, \Pi \Rightarrow \Delta}\,(\text{Cut})$$

Suppose $\langle \bar{x}, a\rangle \in R$ and $\langle\langle \bar{y}, a, \bar{z}\rangle, \delta\rangle \in R$. Then $\prod\bar{x} \leq a$ and $\prod\bar{y} * a * \prod\bar{z} \leq c_\delta$. By monotonicity we have $\prod\bar{y} * \prod\bar{x} * \prod\bar{z} \leq \prod\bar{y} * a * \prod\bar{z}$ and hence $\prod\bar{y} * \prod\bar{x} * \prod\bar{z} \leq c_\delta$. Therefore, $\langle\langle \bar{y}, \bar{x}, \bar{z}\rangle, \delta\rangle \in R$.

• *Rules for* $\vee$:

$$\frac{\Sigma,\varphi,\Gamma\Rightarrow\Delta\qquad\Sigma,\psi,\Gamma\Rightarrow\Delta}{\Sigma,\varphi\vee\psi,\Gamma\Rightarrow\Delta}\;(\vee\Rightarrow)\qquad\frac{\Gamma\Rightarrow\varphi}{\Gamma\Rightarrow\varphi\vee\psi}\;(\Rightarrow\vee_1)\qquad\frac{\Gamma\Rightarrow\psi}{\Gamma\Rightarrow\varphi\vee\psi}\;(\Rightarrow\vee_2)$$

$(\vee\Rightarrow)$: If $\langle\langle\bar{x},a,\bar{y}\rangle,\delta\rangle\in R$ and $\langle\langle\bar{x},b,\bar{y}\rangle,\delta\rangle\in R$, then we have $\prod\bar{x}*a*\prod\bar{y}\leq c_\delta$ and $\prod\bar{x}*b*\prod\bar{y}\leq c_\delta$. Thus, $(\prod\bar{x}*a*\prod\bar{y})\vee(\prod\bar{x}*b*\prod\bar{y})\leq c_\delta$ and by distributivity we have $\prod\bar{x}*(a\vee b)*\prod\bar{y}\leq c_\delta$ and, as a consequence, $\langle\langle\bar{x},a\vee b,\bar{y}\rangle,\delta\rangle\in R$.
$(\Rightarrow\vee_1)$: If $\langle\bar{x},a\rangle\in R$, then $\prod\bar{x}\leq a\leq a\vee b$. Therefore, $\langle\bar{x},a\vee b\rangle\in R$.
$(\Rightarrow\vee_2)$: Analogous to the case above.

• *Rules for* $*$:

$$\frac{\Sigma,\varphi,\psi,\Gamma\Rightarrow\Delta}{\Sigma,\varphi*\psi,\Gamma\Rightarrow\Delta}\;(*\Rightarrow)\qquad\qquad\frac{\Gamma\Rightarrow\varphi\qquad\Pi\Rightarrow\psi}{\Gamma,\Pi\Rightarrow\varphi*\psi}\;(\Rightarrow*)$$

$(*\Rightarrow)$: If $\langle\langle\bar{x},a,b,\bar{y}\rangle,\delta\rangle\in R$, then $\prod\bar{x}*a*b*\prod\bar{y}\leq c_\delta$. Therefore, $\langle\langle\bar{x},a*b,\bar{y}\rangle,\delta\rangle\in R$.
$(\Rightarrow*)$: If $\langle\bar{x},a\rangle\in R$ and $\langle\bar{y},b\rangle\in R$, then $\prod\bar{x}\leq a$ and $\prod\bar{y}\leq b$ and thus, by monotonicity, we have that $\prod\bar{x}*\prod\bar{y}\leq a*b$ and, therefore, $\langle\langle\bar{x},\bar{y}),a*b\rangle\in R$.

• *Rules* $(1\Rightarrow)$ *and* $(\Rightarrow 0)$:

$$\frac{\Sigma,\Gamma\Rightarrow\Delta}{\Sigma,1,\Gamma\Rightarrow\Delta}\;(1\Rightarrow)\qquad\frac{\Gamma\Rightarrow\emptyset}{\Gamma\Rightarrow 0}\;(\Rightarrow 0)$$

$(1\Rightarrow)$: If $\langle\bar{x},\bar{y}\rangle\in R$, then $\prod\bar{x}*\prod\bar{y}\leq c_\delta$ and hence $\prod\bar{x}*1*\prod\bar{y}\leq c_\delta$. Therefore, $\langle\langle\bar{x},1,\bar{y}\rangle,\delta\rangle\in R$.
$(\Rightarrow 0)$: If $\langle\bar{x},\emptyset\rangle\in R$, then $\prod\bar{x}\leq 0$; that is, $\langle\bar{x},0\rangle\in R$.

If $\mathbf{A}$ is a $\mathbb{PM}^{s\ell}$-algebra, we must prove, moreover, that the set $R$ is closed under the interpretation of the rules for $^{\backprime}$ and $'$:

$$\frac{\Gamma\Rightarrow\varphi}{\Gamma,\varphi^{\backprime}\Rightarrow\emptyset}\;(^{\backprime}\Rightarrow)\qquad\qquad\frac{\varphi,\Gamma\Rightarrow\emptyset}{\Gamma\Rightarrow\varphi^{\backprime}}\;(\Rightarrow{}^{\backprime})$$

$$\frac{\Gamma\Rightarrow\varphi}{{}'\varphi,\Gamma\Rightarrow\emptyset}\;({}'\Rightarrow)\qquad\qquad\frac{\Gamma,\varphi\Rightarrow\emptyset}{\Gamma\Rightarrow{}'\varphi}\;(\Rightarrow{}')$$

$(^{\backprime}\Rightarrow)$: If $\langle\bar{x},a\rangle\in R$, then $\prod\bar{x}\leq a$. By monotonicity and the fact that $a^{\backprime}$ is the right pseudocomplement of $a$, we have: $\prod\bar{x}*a^{\backprime}\leq a*a^{\backprime}\leq 0$; that is, $\langle\langle\bar{x},a^{\backprime}\rangle,\emptyset\rangle\in R$.
$(\Rightarrow{}^{\backprime})$: If $\langle\langle a,\bar{x}\rangle,\emptyset\rangle\in R$, then $a*\prod\bar{x}\leq 0$ and, by the pseudocomplementation law, $\prod\bar{x}\leq a^{\backprime}$; that is, $\langle\bar{x},a^{\backprime}\rangle\in R$.
$({}'\Rightarrow),(\Rightarrow{}')$: We can proceed analogously by using the properties of the left pseudocomplement.

If $\mathbf{A}$ is one of the classes $\mathring{\mathbb{M}}^{\ell}$, $\mathbb{PM}^{\ell}$ or $\mathbb{FL}$, that is, if $\Psi$ contains the connective $\wedge$, we must see, moreover, that the set $R$ is closed under the interpretation of the introduction rules for this connective:

$$\frac{\Sigma,\varphi,\Gamma\Rightarrow\Delta}{\Sigma,\varphi\wedge\psi,\Gamma\Rightarrow\Delta}\;(\wedge_1\Rightarrow)\qquad\frac{\Sigma,\psi,\Gamma\Rightarrow\Delta}{\Sigma,\varphi\wedge\psi,\Gamma\Rightarrow\Delta}\;(\wedge_2\Rightarrow)\qquad\frac{\Gamma\Rightarrow\varphi\quad\Gamma\Rightarrow\psi}{\Gamma\Rightarrow\varphi\wedge\psi}\;(\Rightarrow\wedge)$$

$(\wedge_1\Rightarrow)$: If $\langle\langle\bar{x},a,\bar{y}\rangle,\delta\rangle\in R$, then $\prod\bar{x}*a*\prod\bar{y}\leq c_\delta$. From $a\wedge b\leq a$ by monotonicity we obtain $\prod\bar{x}*(a\wedge b)*\prod\bar{y}\leq\prod\bar{x}*a*\prod\bar{y}$. Thus $\prod\bar{x}*(a\wedge b)*\prod\bar{y}\leq c_\delta$. Therefore, $\langle\langle\bar{x},a\wedge b,\bar{y}\rangle,\delta\rangle\in R$.

$(\wedge_2 \Rightarrow)$: Analogously.
$(\Rightarrow \wedge)$: If $\langle \bar{x}, a\rangle \in R$ and $\langle \bar{x}, b\rangle \in R$, then $\prod \bar{x} \leq a$ and $\prod \bar{x} \leq b$. Thus $\prod \bar{x} \leq a \wedge b$. Therefore, $\langle \bar{x}, a \wedge b\rangle \in R$.

If $\mathbf{A} \in \mathbb{FL}$, we must show, moreover, that $R$ is closed under the introduction rules to the connectives $\backslash$ and $/$:

$$\frac{\Gamma \Rightarrow \varphi \qquad \Sigma, \psi, \Pi \Rightarrow \Delta}{\Sigma, \Gamma, \varphi\backslash\psi, \Pi \Rightarrow \Delta}\ (\backslash \Rightarrow) \qquad \frac{\varphi, \Gamma \Rightarrow \psi}{\Gamma \Rightarrow \varphi\backslash\psi}\ (\Rightarrow \backslash)$$

$$\frac{\Gamma \Rightarrow \varphi \quad \Sigma, \psi, \Pi \Rightarrow \Delta}{\Sigma, \psi/\varphi, \Gamma, \Pi \Rightarrow \Delta}\ (/ \Rightarrow) \qquad \frac{\Gamma, \varphi \Rightarrow \psi}{\Gamma \Rightarrow \psi/\varphi}\ (\Rightarrow /)$$

$(\backslash \Rightarrow)$: If $\langle \bar{x}, a\rangle \in R$ and $\langle\langle \bar{y}, b, \bar{z}\rangle, \delta\rangle \in R$, then $\prod \bar{x} \leq a$ and $\prod \bar{y} * b * \prod \bar{z} \leq c_\delta$. By monotonicity and the properties of the right residuum we have $\prod \bar{x} * (a\backslash b) \leq a * (a\backslash b) \leq b$ and hence: $\prod \bar{y} * \prod \bar{x} * (a\backslash b) * \prod \bar{z} \leq \prod \bar{y} * b * \prod \bar{z} \leq c_\delta$. As a consequence, $\langle\langle \bar{y}, \bar{x}, a\backslash b, \bar{z}\rangle, \delta\rangle \in R$.
$(\Rightarrow \backslash)$: If $\langle\langle a, \bar{x}\rangle, b\rangle \in R$, then $a * \prod \bar{x} \leq b$ and, by the law of residuation, $\prod \bar{x} \leq a\backslash b$. Therefore, $\langle \bar{x}, a\backslash b\rangle \in R$.
$(/ \Rightarrow)$, $(\Rightarrow /)$: We can proceed analogously using the properties of the left residuum.

Finally, if $\mathbf{A} \in \mathbb{K}_\sigma[\Psi]$, with $\sigma$ non-empty, we want to see that $R$ is closed under all the structural rules codified by $\sigma$.

• *Exchange*:

$$\frac{\Gamma, \varphi, \psi, \Pi \Rightarrow \Delta}{\Gamma, \psi, \varphi, \Pi \Rightarrow \Delta}\ (e \Rightarrow)$$

If $\langle\langle \bar{x}, a, b, \bar{y}\rangle, \delta\rangle \in R$, then $\prod \bar{x} * a * b * \prod \bar{y} \leq c_\delta$. Thus, if $\mathbf{A}$ is commutative, we have $\prod \bar{x} * b * a * \prod \bar{y} \leq c_\delta$. Consequently, $\langle\langle \bar{x}, b, a, \bar{y}\rangle, \delta\rangle \in R$.

• *Weakening*:

$$\frac{\Sigma, \Gamma \Rightarrow \Delta}{\Sigma, \varphi, \Gamma \Rightarrow \Delta}\ (w \Rightarrow) \qquad \frac{\Gamma \Rightarrow \emptyset}{\Gamma \Rightarrow \varphi}\ (\Rightarrow w)$$

$(w \Rightarrow)$: If $\langle\langle \bar{x}, \bar{y}\rangle\delta\rangle \in R$, then $\prod \bar{x} * \prod \bar{y} \leq c_\delta$; thus, if $\mathbf{A}$ is integral, by monotonicity we have: $\prod \bar{x} * a * \prod \bar{y} \leq \prod \bar{x} * 1 * \prod \bar{y} \leq \prod \bar{x} * \prod \bar{y} \leq c_\delta$. That is, $\langle\langle \bar{x}, a, \bar{y}\rangle, \delta\rangle \in R$.
$(\Rightarrow w)$: If $\langle \bar{x}, \emptyset\rangle \in R$, then $\prod \bar{x} \leq 0$. Thus, if 0 is the minimum, then $\prod \bar{x} \leq a$. Therefore, $\langle \bar{x}, a\rangle \in R$.

• *Contraction*:

$$\frac{\Sigma, \varphi, \varphi, \Gamma \Rightarrow \Delta}{\Sigma, \varphi, \Gamma \Rightarrow \Delta}\ (c \Rightarrow)$$

If $\langle\langle \bar{x}, a, a, \bar{y}\rangle, \delta\rangle \in R$, then we have: $\prod \bar{x} * a * a * \prod \bar{y} \leq c_\delta$. From this, when $\mathbf{A}$ has the property of increasing idempotency, we obtain $\prod \bar{x} * a * \prod \bar{y} \leq c_\delta$. Therefore, $\langle\langle \bar{x}, a, \bar{y}\rangle, \delta\rangle \in R$. □

In the following lemma we prove that condition 4) of Lemma 1 is satisfied.

**Lemma 15.** *For every theory $\Phi \in Th\,\mathcal{FL}_\sigma[\Psi]$, the set:*

$$\theta_\Phi := \{\langle \varphi, \psi\rangle \in Fm^2_\Psi : \rho(\varphi \approx \psi) \subseteq \Phi\}$$

*is a congruence relative to the quasivariety $\mathbb{K}_\sigma[\Psi]$.*

*Proof.* Let $\Phi$ be an $\mathcal{FL}_\sigma[\Psi]$-theory. By the definition of the translation $\rho$ we have

$$\theta_\Phi := \{\langle\varphi,\psi\rangle \in Fm_\Psi^2 : \{\varphi\Rightarrow\psi, \psi\Rightarrow\varphi\}\subseteq\Phi\}$$

but, by Corollary 6, this is the Leibniz congruence $\Omega\Phi$ of the theory $\Phi$.

Let us denote $\mathbf{Q}[\Psi] := Fm_\Psi/\Omega\Phi$. We want to see that $\mathbf{Q}[\Psi]$ is a $\mathbb{K}_\sigma[\Psi]$-algebra. To see this, we will prove the following conditions:

(1) If $\Psi = \langle\vee, *, 0, 1\rangle$, then $\mathbf{Q}[\Psi] \in \mathring{\mathbb{M}}^{s\ell}$.
(2) If $\Psi = \langle\vee, *, \backslash, /, 0, 1\rangle$, then $\mathbf{Q}[\Psi] \in \mathbb{PM}^{s\ell}$.
(3) If $\langle\wedge\rangle \leq \Psi$, then $\mathbf{Q}[\Psi]$ satisfies a set of equations that defines the class of lattices.
(4) If $\Psi = \mathfrak{L}$, then $\mathbf{Q}[\Psi] \in \mathbb{FL}$.
(5) If $e \leq \sigma$, then $\mathbf{Q}[\Psi] \in \mathbb{K}_e[\Psi]$.
(6) If $w_l \leq \sigma$, then $\mathbf{Q}[\Psi] \in \mathbb{K}_{w_l}[\Psi]$.
(7) If $w_r \leq \sigma$, then $\mathbf{Q}[\Psi] \in \mathbb{K}_{w_r}[\Psi]$.
(8) If $c \leq \sigma$ then, $\mathbf{Q}[\Psi] \in \mathbb{K}_c[\Psi]$.

1) If $\varphi \approx \psi$ belongs to a set of equations that defines the $\mathring{\mathbb{M}}^{s\ell}$-algebras (Theorem 3), we will see that the sequents $\varphi\Rightarrow\psi$ and $\psi\Rightarrow\varphi$ are derivable and thus $\langle\varphi,\psi\rangle\in\Omega\Phi$.

• *Equations defining the* $\langle\vee\rangle$*-semilattices:*
$(x\vee y)\vee z \approx x\vee(y\vee z)$; $x\vee y\approx y\vee x$; $x\vee x\approx x$. The corresponding sequents are easily obtained by using $(Axiom)$ and the introduction rules for $\vee$.

• *Equations defining the* $\langle *, 1\rangle$*-monoids:*
$(x*y)*z \approx x*(y*z)$. The corresponding sequents are easily obtained by using $(Axiom)$ and the rules $(\Rightarrow *)$ and $(*\Rightarrow)$.
$1*x\approx x$. The sequent $1*\varphi\Rightarrow\varphi$ is obtained by applying $(*\Rightarrow)$ to the sequent $1,\varphi\Rightarrow\varphi$ which it is obtained from $\varphi\Rightarrow\varphi$ using $(1\Rightarrow)$. The sequent $\varphi\Rightarrow 1*\varphi$ is obtained from the axioms $\emptyset\Rightarrow 1$ and $\varphi\Rightarrow\varphi$ using the rule $(\Rightarrow *)$.

• *Equations of distributivity of* $*$ *with respect to* $\vee$*:*
$(x\vee y)*z\approx(x*z)\vee(y*z)$. Let us consider the following derivations:

$$\dfrac{\dfrac{\dfrac{\dfrac{\varphi\Rightarrow\varphi \qquad \gamma\Rightarrow\gamma}{\varphi,\gamma\Rightarrow\varphi*\gamma}\,(\Rightarrow *)}{\varphi,\gamma\Rightarrow(\varphi*\gamma)\vee(\psi*\gamma)}\,(\Rightarrow\vee_1) \qquad \dfrac{\dfrac{\psi\Rightarrow\psi \qquad \gamma\Rightarrow\gamma}{\psi,\gamma\Rightarrow\psi*\gamma}\,(\Rightarrow *)}{\psi,\gamma\Rightarrow(\varphi*\gamma)\vee(\psi*\gamma)}\,(\Rightarrow\vee_2)}{\varphi\vee\psi,\gamma\Rightarrow(\varphi*\gamma)\vee(\psi*\gamma)}\,(\vee\Rightarrow)}{(\varphi\vee\psi)*\gamma\Rightarrow(\varphi*\gamma)\vee(\psi*\gamma)}\,(*\Rightarrow)$$

$$\dfrac{\dfrac{\dfrac{\dfrac{\dfrac{\varphi\Rightarrow\varphi}{\varphi\Rightarrow\varphi\vee\psi}\,(\Rightarrow\vee_1) \qquad \gamma\Rightarrow\gamma}{\varphi,\gamma\Rightarrow(\varphi\vee\psi)*\gamma}\,(\Rightarrow *)}{\varphi*\gamma\Rightarrow(\varphi\vee\psi)*\gamma}\,(*\Rightarrow) \qquad \dfrac{\dfrac{\dfrac{\psi\Rightarrow\psi}{\psi\Rightarrow\varphi\vee\psi}\,(\Rightarrow\vee_2) \qquad \gamma\Rightarrow\gamma}{\psi,\gamma\Rightarrow(\varphi\vee\psi)*\gamma}\,(\Rightarrow *)}{\psi*\gamma\Rightarrow(\varphi\vee\psi)*\gamma}\,(*\Rightarrow)}{(\varphi*\gamma)\vee(\psi*\gamma)\Rightarrow(\varphi\vee\psi)*\gamma}\,(\vee\Rightarrow)$$

$z*(x\vee y)\approx(z*x)\vee(z*x)$. The corresponding sequents are mirror images of the previous ones and therefore are also derivable (Theorem 3).

2) If $\Psi = \langle \vee, *, ^{\backprime}, ', 0, 1\rangle$, in addition we must prove that $\mathbf{Q}[\Psi]$ satisfies the equations in Theorem 16 relative to the pseudocomplements. Indeed:

• *Equations* $1^{\backprime} \approx 0$; $1 \preccurlyeq 0^{\backprime}$; $x * (y * x)^{\backprime} \preccurlyeq y^{\backprime}$.

$1^{\backprime} \approx 0$. The sequents $^{\backprime}1 \Rightarrow 0$ and $0 \Rightarrow 1^{\backprime}$ are $\mathcal{FL}_\sigma[\Psi]$-derivable:

$$\cfrac{\cfrac{\emptyset \Rightarrow 1}{1^{\backprime} \Rightarrow \emptyset}\,(^{\backprime} \Rightarrow)}{1^{\backprime} \Rightarrow 0}\,(\Rightarrow 0) \qquad \cfrac{\cfrac{0 \Rightarrow \emptyset}{1, 0 \Rightarrow \emptyset}\,(1 \Rightarrow)}{0 \Rightarrow 1^{\backprime}}\,(\Rightarrow {}^{\backprime})$$

$1 \preccurlyeq 0^{\backprime}$. The sequent $1 \Rightarrow 0^{\backprime}$ is obtained by means of the following derivation:

$$\cfrac{\cfrac{0 \Rightarrow \emptyset}{\emptyset \Rightarrow 0^{\backprime}}\,(\Rightarrow {}^{\backprime})}{1 \Rightarrow 0^{\backprime}}\,(1 \Rightarrow)$$

$x * (y * x)^{\backprime} \preccurlyeq y^{\backprime}$. The sequent $\varphi * (\psi * \varphi)^{\backprime} \Rightarrow \psi^{\backprime}$ is obtained by means of the following derivation:

$$\cfrac{\cfrac{\cfrac{\cfrac{\psi \Rightarrow \psi \qquad \varphi \Rightarrow \varphi}{\psi, \varphi \Rightarrow \psi * \varphi}\,(\Rightarrow *)}{\psi, \varphi, (\psi * \varphi)^{\backprime} \Rightarrow \emptyset}\,(^{\backprime} \Rightarrow)}{\varphi, (\psi * \varphi)^{\backprime} \Rightarrow \psi^{\backprime}}\,(\Rightarrow {}^{\backprime})}{\varphi * (\psi * \varphi)^{\backprime} \Rightarrow \psi^{\backprime}}\,(* \Rightarrow)$$

• *Equations* $'1 \approx 0$; $1 \preccurlyeq {}'0$; $(x * {}'(y * x)) \preccurlyeq {}'y$. The sequents corresponding to these equations are the mirror images of the sequents corresponding to the equations concerning the right pseudocomplement and, by Theorem 3, they are derivable.

• *Equations* $(x \vee y)^{\backprime} \preccurlyeq x^{\backprime}$; $'(x \vee y) \preccurlyeq {}'x$. By the law of the mirror images, it is sufficient to see that the sequent $(\varphi \vee \psi)^{\backprime} \Rightarrow \varphi^{\backprime}$ is derivable. Let us consider the following derivation:

$$\cfrac{\cfrac{\cfrac{\varphi \Rightarrow \varphi}{\varphi \Rightarrow \varphi \vee \psi}\,(\Rightarrow \vee_1)}{\varphi, (\varphi \vee \psi)^{\backprime} \Rightarrow \emptyset}\,(^{\backprime} \Rightarrow)}{(\varphi \vee \psi)^{\backprime} \Rightarrow \varphi^{\backprime}}\,(\Rightarrow {}^{\backprime})$$

3) If $\wedge \in \Psi$, we must see that $\mathbf{Q}[\Psi]$ satisfies a set of equations defining the lattices. We have already seen that using the introduction rules for $\vee$ we can prove that $\mathbf{Q}[\Psi]$ satisfies the commutativity and the associativity of the operation $\vee$. Therefore, it will be sufficient to prove that it also satisfies the commutativity and the associativity of the operation $\wedge$ and the absorption laws:

• *Equations* $x \wedge y \approx y \wedge x$; $x \wedge (y \wedge z) \approx (x \wedge y) \wedge z$; $x \wedge y \preccurlyeq x$; $x \preccurlyeq x \vee y$.

The sequents corresponding to the first three equations are easily obtained using $(Axiom)$ and the introduction rules for the connective $\wedge$. The sequent corresponding to the last equation is obtained from $(Axiom)$ using $(\Rightarrow \vee_1)$.

4) If $\Psi = \langle \vee, \wedge, *, \backslash, /, {}^{\backprime}, {}', 0, 1\rangle$, to prove that $\mathbf{Q}[\Psi]$ is an $\mathbb{FL}$-algebra, it only remains to see that it satisfies the equations involving the residuals and the pseudocomplements in the equational characterization of the class $\mathbb{FL}$ in Theorem 11

• *Equations* $x * ((x\backslash z) \wedge y) \preccurlyeq z$, $y \preccurlyeq x\backslash((x * y) \vee z)$. The corresponding formal proofs are the following:

$$\dfrac{\dfrac{\dfrac{\varphi \Rightarrow \varphi \qquad \gamma \Rightarrow \gamma}{\varphi, \varphi\backslash\gamma \Rightarrow \gamma}\,(\backslash \Rightarrow)}{\varphi, (\varphi\backslash\gamma) \wedge \psi \Rightarrow \gamma}\,(\wedge_1 \Rightarrow)}{\varphi * ((\varphi\backslash\gamma) \wedge \psi) \Rightarrow \gamma}\,(* \Rightarrow) \qquad \dfrac{\dfrac{\dfrac{\varphi \Rightarrow \varphi \qquad \psi \Rightarrow \psi}{\varphi, \psi \Rightarrow \varphi * \psi}\,(\Rightarrow *)}{\varphi, \psi \Rightarrow (\varphi * \psi) \vee \gamma}\,(\Rightarrow \vee_1)}{\psi \Rightarrow \varphi\backslash((\varphi * \psi) \vee \gamma)}\,(\Rightarrow \backslash)$$

• *Equations* $((z/x) \wedge y) * x \preceq z$, $y \preceq ((y * x) \vee z)/x$. The sequents corresponding to these equations are the mirror images of the sequents corresponding to the equations involving the right residuum and, therefore, they are also derivable.

• *Equations defining the pseudocomplements:* $x^{\backprime} \approx x\backslash 0$; ${}'x \approx 0/x$. Let us consider the following derivations corresponding to the equation defining the right pseudocomplement:

$$\dfrac{\dfrac{\dfrac{\varphi \Rightarrow \varphi}{\varphi, \varphi^{\backprime} \Rightarrow \emptyset}\,({}^{\backprime} \Rightarrow)}{\varphi, \varphi^{\backprime} \Rightarrow 0}\,(\Rightarrow 0)}{\varphi^{\backprime} \Rightarrow \varphi\backslash 0}\,(\Rightarrow \backslash) \qquad \dfrac{\dfrac{\varphi \Rightarrow \varphi \qquad 0 \Rightarrow \emptyset}{\varphi, \varphi\backslash 0 \Rightarrow \emptyset}\,(\backslash \Rightarrow)}{\varphi\backslash 0 \Rightarrow \varphi^{\backprime}}\,(\Rightarrow {}^{\backprime})$$

The sequents corresponding to the left pseudocomplement are the mirror images of these and thus are also derivable.

5) If $e \leq \sigma$ we must prove that $\mathbf{Q}[\Psi]$ satisfies $x * y \approx y * x$. By symmetry, it is sufficient to derive the sequent $\varphi * \psi \Rightarrow \psi * \varphi$, which is easily obtained using (*Axiom*), $(\Rightarrow *)$, $(e \Rightarrow)$ and $(* \Rightarrow)$.

6) If $w_l \leq \sigma$, we must prove that $\mathbf{Q}[\Psi]$ satisfies $x \preccurlyeq 1$. It will be sufficient to prove that the sequent $\varphi \Rightarrow 1$ is derivable. Indeed, from $\emptyset \Rightarrow 1$ we obtain $\varphi \Rightarrow 1$ using $(w \Rightarrow)$.

7) If $w_r \leq \sigma$, we must prove that $\mathbf{Q}[\Psi]$ satisfies $0 \preccurlyeq x$. It will be sufficient to prove that $0 \Rightarrow \varphi$ is derivable. Indeed, from the sequent $0 \Rightarrow \emptyset$ we obtain $0 \Rightarrow \varphi$ using $(\Rightarrow w)$.

8) If $w_r \leq \sigma$, we must prove that $\mathbf{Q}[\Psi]$ satisfies $x \preccurlyeq x * x$. It will be sufficient to prove that $\varphi \Rightarrow \varphi * \varphi$ is derivable. Let us consider the following derivation:

$$\dfrac{\dfrac{\varphi \Rightarrow \varphi \qquad \varphi \Rightarrow \varphi}{\varphi, \varphi \Rightarrow \varphi * \varphi}\,(\Rightarrow *)}{\varphi \Rightarrow \varphi * \varphi}\,(c \Rightarrow)$$

□

As a consequence of Lemma 12, Lemma 13, Lemma 14, and Lemma 15, we already have the algebraization result:

**Theorem 21** (Algebraization)**.** *Every Gentzen system* $\mathcal{FL}_\sigma[\Psi]$ *is algebraizable, with the variety* $\mathbb{K}_\sigma[\Psi]$ *as its equivalent algebraic semantics.*

TABLE 2. Systems $\mathcal{FL}_\sigma[\Psi]$ and their equivalent quasivariety semantics

| **Gentzen system** | **e.q.s.** |
|---:|---|
| $\mathcal{FL}_\sigma[\vee, *, 0, 1]$ | $\mathring{\mathbb{M}}^{s\ell}_\sigma$ |
| $\mathcal{FL}_\sigma[\vee, \wedge, *, 0, 1]$ | $\mathring{\mathbb{M}}^{\ell}_\sigma$ |
| $\mathcal{FL}_\sigma[\vee, *, \neg, 0, 1]$ | $\mathbb{PM}^{s\ell}_\sigma$ |
| $\mathcal{FL}_\sigma[\vee, \wedge, *, \neg, 0, 1]$ | $\mathbb{PM}^{\ell}_\sigma$ |
| $\mathcal{FL}_\sigma$ | $\mathbb{FL}_\sigma$ |

*Proof.* We use the translations $\tau$ and $\rho$ in Definition 30. With these translations, the four conditions of Lemma 1 are satisfied: condition 1) by Lemma 12, condition 2) by Lemma 13, condition 3) by Lemma 14, and condition 4) by Lemma 15. □

Table 9 shows the subsystems $\mathcal{FL}_\sigma[\Psi]$ and the corresponding classes of algebras that are their equivalent quasivariety semantics (e.q.s.).

The results of algebraization allow us to obtain the following corollaries as consequences.

**Corollary 25.** *If $\mathbf{A} \in \mathbb{K}_\sigma[\Psi]$, then the sequential Leibniz operator $\Omega_{\mathbf{A}}$ is an isomorphism between the lattice of $\mathcal{FL}_\sigma[\Psi]$-filters and the lattice of $\mathbb{K}_\sigma[\Psi]$-congruences.*

*Proof.* It is an immediate consequence of Theorem 1 and Theorem 21. □

**Corollary 26.** *The subdirectly irreducible algebras of one of the classes $\mathbb{K}_\sigma[\Psi]$ are exactly the algebras of the class with the smallest non-trivial $\mathcal{FL}_\sigma[\Psi]$-filter.*

The algebraization result for $\sigma = ewc$ and the sublanguage $\langle \vee, \wedge, \neg \rangle$ appears in Rebagliato & Verdú (1995). The results for the case $\sigma = ew$ and the sublanguages $\langle \vee, *, \neg, 0, 1 \rangle$ and $\langle \vee, \wedge, *, \neg, 0, 1 \rangle$ appear in Bou et al. (2006). The results for the case $\sigma = ew$ and the sublanguages $\langle \vee, *, 0, 1 \rangle$ and $\langle \vee, \wedge, *, 0, 1 \rangle$ are in Adillon & Verdú (2000); Adillon et al. (2007).

## Part 4. COMPLETIONS AND SUBREDUCTS

In this part we show that the method called *ideal completion* allows us to embed every algebra in $\mathring{\mathbb{M}}^{s\ell}_\sigma$, $\mathring{\mathbb{M}}^{\ell}_\sigma$, $\mathbb{PM}^{s\ell}_\sigma$ or $\mathbb{PM}^{\ell}_\sigma$ in a complete $\mathbb{FL}_\sigma$-algebra. Such embeddings have as a consequence that the classes $\mathring{\mathbb{M}}^{s\ell}_\sigma$, $\mathring{\mathbb{M}}^{\ell}_\sigma$, $\mathbb{PM}^{s\ell}_\sigma$ and $\mathbb{PM}^{\ell}_\sigma$ are the classes of all the subreducts of the algebras in the $\mathbb{FL}_\sigma$ class. These results will be used in Part 5 to prove that the subsystems $\mathcal{FL}_\sigma[\vee, *, 0, 1]$, $\mathcal{FL}_\sigma[\vee, \wedge, *, 0, 1]$, $\mathcal{FL}_\sigma[\vee, *, \grave{}, ', 0, 1]$, and $\mathcal{FL}_\sigma[\vee, \wedge, *, \grave{}, ', 0, 1]$ are fragments of $\mathcal{FL}_\sigma$ in the languages $\langle \vee, *, 0, 1 \rangle$, $\langle \vee, \wedge, *, 0, 1 \rangle$, $\langle \vee, *, \grave{}, ', 0, 1 \rangle$, and $\langle \vee, \wedge, *, \grave{}, ', 0, 1 \rangle$. These last results, jointly with the algebraization, will allow us to prove that the external systems $\mathfrak{e}\mathcal{FL}_\sigma[\vee, *, 0, 1]$, $\mathfrak{e}\mathcal{FL}_\sigma[\vee, \wedge, *, 0, 1]$, $\mathfrak{e}\mathcal{FL}_\sigma[\vee, *, \grave{}, ', 0, 1]$, and $\mathfrak{e}\mathcal{FL}_\sigma[\vee, \wedge, *, \grave{}, ', 0, 1]$ are also the fragments of $\mathfrak{e}\mathcal{FL}_\sigma$ in the languages considered. In Section 10, we recap some basic concepts and results for complete algebras; and in Section 11, we recap some of the results regarding completions that we will need. Finally, in Section 12, we present the construction of ideal completion and we obtain the results regarding subreducts.

## 10. Basic Concepts.

In this section, we recap some basic concepts concerning complete lattices. The notion of a complete algebra for the classes $\mathring{\mathbb{M}}^{s\ell}$, $\mathring{\mathbb{M}}^{\ell}$, $\mathbb{PM}^{s\ell}$, $\mathbb{PM}^{\ell}$ and $\mathbb{FL}$ is defined; it is stated that every complete $\mathring{\mathbb{M}}^{s\ell}$-algebra is the reduct of a complete $\mathring{\mathbb{M}}^{\ell}$-algebra and that every complete $\mathbb{PM}^{s\ell}$-algebra is the reduct of a complete $\mathbb{PM}^{\ell}$-algebra.

**Remark 17.** *If $\mathcal{A} = \langle A, \leq \rangle$ is an ordered set and $X \subseteq A$, we denote by $X^{\rightarrow}$ the set of all the upper bounds of $X$ in $A$ and by $X^{\leftarrow}$ the set of all the lower bounds of $X$ in $A$. Note that the supremum of $X$ in $A$, if it exists, is the minimum of the set $X^{\rightarrow}$; and that the infimum of $X$ in $A$, if it exists, is the maximum of the set $X^{\leftarrow}$.*

**Proposition 31.** (Cf. Davey & Priestley (1990) ) *Let $\mathcal{A} = \langle A, \leq \rangle$ be an ordered set. Then:*

a) *the supremum of $A$ in $A$ exists if and only if $A$ has a greatest element $\top$ and, in this case, $\bigvee_{\mathcal{A}} A = \top$;*
b) *the supremum of $\emptyset$ in $A$ exists if and only if $A$ has a smallest element $\perp$ and, in this case, $\bigvee_{\mathcal{A}} \emptyset = \perp$.*

*Proof.* *a*): If $A$ has a *top* $\top$, then $A^{\rightarrow} = \{\top\}$ and thus, $\bigvee_{\mathcal{A}} A = \top$. If $A$ has no *top*, then $A^{\rightarrow} = \emptyset$ and so $\bigvee_{\mathcal{A}} A$ does not exist.
*b*): Each element $a \in A$ satisfies (emptily) $x \leq a$ for every $x \in \emptyset$. Hence, $\emptyset^{\rightarrow} = A$ and consequently $\bigvee_{\mathcal{A}} \emptyset$ exists if and only if $A$ has a *bottom* $\perp$. In this case, $\bigvee_{\mathcal{A}} \emptyset = \perp$. □

**Proposition 32.** *Every complete lattice has a minimum and a maximum element.*

*Proof.* It is an immediate consequence of Proposition 31. □

**Proposition 33.** (Cf. Davey & Priestley (1990) ) *Let $\mathcal{A} = \langle A, \leq \rangle$ be a non-empty ordered set. Then, the following conditions are equivalent:*

a) *$\mathcal{A}$ is a complete lattice.*
b) *$\bigvee_{\mathcal{A}} X$ exists for every subset $X \subseteq A$.*

*Proof.* *a*) $\Rightarrow$ *b*) is trivial.
*b*) $\Rightarrow$ *a*): Suppose that every subset of $A$ has a supremum. Let $X \subseteq A$ and $a = \bigvee_{\mathcal{A}} X^{\leftarrow}$. Since $X \subseteq (X^{\leftarrow})^{\rightarrow}$ and $a$ is the minimum of $(X^{\leftarrow})^{\rightarrow}$, we have $a \leq x$ for each $x \in X$ and so $a \in X^{\leftarrow}$. On the other hand, since $a \in (X^{\leftarrow})^{\rightarrow}$, we have that if $z \in X^{\leftarrow}$, then $z \leq a$ and hence $a$ is the maximum of $X^{\leftarrow}$. Therefore, $\bigwedge_{\mathcal{A}} X = \bigvee_{\mathcal{A}} X^{\leftarrow}$. □

**Definition 31.** *If $\mathbf{A}$ is an $\mathring{\mathbb{M}}^{s\ell}$-algebra or an $\mathring{\mathbb{M}}^{\ell}$-algebra, we will say that it is* complete *if every subset of its universe, $A$, has a supremum.*

**Proposition 34.** *Let $\mathbf{A}$ be an $\mathring{\mathbb{M}}^{s\ell}$-algebra or an $\mathring{\mathbb{M}}^{\ell}$-algebra. $\mathbf{A}$ is complete if and only if the ordered set associated with its semilatticed reduct is a complete lattice.*

*Proof.* By Proposition 33. □

The notion of a complete algebra for the algebras in $\mathbb{PM}^{s\ell}$, $\mathbb{PM}^{\ell}$ and $\mathbb{FL}$ is defined in an analogous way.

**Proposition 35.** *Every complete algebra in* $\mathring{\mathbb{M}}^{s\ell}$, $\mathring{\mathbb{M}}^{\ell}$, $\mathbb{PM}^{s\ell}$, $\mathbb{PM}^{\ell}$ *and* $\mathbb{FL}$ *has a minimum element and a maximum element.*

*Proof.* By Proposition 32. □

**Proposition 36.** *Let* $\sigma \leq ew_l w_r c$. *Every complete* $\mathring{\mathbb{M}}^{s\ell}_{\sigma}$*-algebra is the* $\langle \vee, *, 0, 1\rangle$*-reduct of a complete* $\mathring{\mathbb{M}}^{\ell}_{\sigma}$*-algebra.*

*Proof.* If $\mathbf{A}$ is a complete $\mathring{\mathbb{M}}^{s\ell}_{\sigma}$-algebra, then the ordered set associated with its semilatticed reduct is a complete lattice. It is clear that $\mathbf{A}$ is the $\langle \vee, *, 0, 1\rangle$-reduct of the complete $\mathring{\mathbb{M}}^{\ell}_{\sigma}$-algebra of universe $A$ where the operation $\wedge$ is defined in the following way: for each $a, b \in A$,

$$a \wedge b =: \bigvee_{\mathbf{A}} \{x \in A : x \leq a \text{ i } x \leq b\}.$$

The rest of the operations are those in $\mathbf{A}$. □

The following result is obtained in an analogous way.

**Proposition 37.** *Every complete* $\mathbb{PM}^{s\ell}_{\sigma}$*-algebra is the* $\langle \vee, *, ^{\backprime}, ^{\prime}, 0, 1\rangle$*-reduct of a complete* $\mathbb{PM}^{\ell}_{\sigma}$*-algebra.*

## 11. Completions.

Among $\mathbb{FL}_{\sigma}$-algebras, the complete ones are particularly interesting, because of the following theorem.

**Theorem 22.** *Every* $\mathbb{FL}_{\sigma}$*-algebra is embeddable in a complete* $\mathbb{FL}_{\sigma}$*-algebra.*

There are at least two methods in the literature that are well known as means to obtaining these completions: the *Dedekind-MacNeille completion* and the *ideal completion.* Before explaining how these two methods work, we will recall a characterization obtained by Ono (see Ono (1993, 2003a)) in which a complete $\mathbb{FL}$-algebra is constructed from a monoid $\mathbf{M} = \langle M, *, 1\rangle$ and a closure operator on $\mathcal{P}(M)$ that satisfies a certain additional condition with respect to the monoidal operation. Here we present this result in the most general case, that is, in a language containing two residuals (left and right) and two pseudocomplements. The inclusion of the two pseudocomplements in the language allows us to see the $\mathbb{PM}^{s\ell}$-algebras and the $\mathbb{PM}^{\ell}$-algebras as subreducts of the $\mathbb{FL}$-algebras, as will be shown in Section 12

**Remark 18** (Notation)**.** *Given a monoid* $\mathbf{M} = \langle M, *, 1\rangle$ *and two subsets* $X, Y \subseteq M$, *we will denote by* $X * Y$ *the set* $\{a * b : a \in X,\ b \in Y\}$.

**Proposition 38.** (Cf. Ono (1993, 2003b)) *Suppose that* $\mathbf{M} = \langle M, *, 1\rangle$ *is a monoid and that the mapping* $C : \mathcal{P}(M) \longrightarrow \mathcal{P}(M)$ *is a closure operator that satisfies*[16]*:*

$$(*) \qquad C(X) * C(Y) \subseteq C(X * Y),$$

[16]For notions concerning closure operators, see Section 2.4

*for each $X, Y \subseteq M$. Then, for every $C$-closed $D$, the structure*

$$\mathbf{C}_{\mathbf{M}}^{D} = \langle \mathcal{C}_M, \vee_C, \cap, *_C, \backslash, /, ^{\backprime}, {}', D, C(1) \rangle$$

*is a complete $\mathbb{FL}$-algebra, where*

$\mathcal{C}_M = \{X \subseteq M : C(X) = X\}$ *(the closure system associated with $C$),*
$X \vee_C Y = C(X \cup Y)$,
$X *_C Y = C(X * Y)$,
$X \backslash Y = \{z \in M : x * z \in Y$ *for every* $x \in X\}$,
$Y/X = \{z \in M : z * x \in Y$ *for every* $x \in X\}$,
$X^{\backprime} = \{z \in M : x * z \in D$ *for every* $x \in X\}$,
${}'X = \{z \in M : z * x \in D$ *for every* $x \in X\}$,

*and for every family $\{X_i\}_{i \in I} \subseteq \mathcal{P}(M)$,*

$\bigvee_{\mathbf{C_M}} \{X_i\}_{i \in I} = C(\bigcup_{i \in I} X_i)$,
$\bigwedge_{\mathbf{C_M}} \{X_i\}_{i \in I} = \bigcap_{i \in I} X_i$.

Observe that from the definitions of the residuals and the pseudocomplements we have that $X^{\backprime} = X \backslash D$ and ${}'X = D/X$.

In the next proposition we recap two well-known identities which are satisfied by all closure operators.

**Proposition 39.** *Let $M$ be a set, $\{X_i\}_{i \in I}$ be a family of subsets of $M$ and $C$ be a closure operator on $\mathcal{P}(M)$. Then:*

(1) $C(\bigcup_{i \in I} X_i) = C(\bigcup_{i \in I} C(X_i))$,
(2) $\bigcap_{i \in I} C(X_i) = C(\bigcap_{i \in I} C(X_i))$.

In the following proposition we state a condition equivalent to condition $(*)$.

**Proposition 40.** *Let $\mathbf{M} = \langle M, *, 1 \rangle$ be a monoid and $C$ be a closure operator on $\mathcal{P}(M)$. For each $X, Y \subseteq M$, the following conditions are equivalent:*

(a) $C(X) * C(Y) \subseteq C(X * Y)$,
(b) $C(C(X) * C(Y)) = C(X * Y)$.

*Proof.* Let $X, Y \subseteq M$.
$(a) \Rightarrow (b)$: On the one hand, by (a) we have $C(X) * C(Y) \subseteq C(X * Y)$. From this, using the properties of closure operators, we obtain $C(C(X) * C(Y)) \subseteq CC(X * Y) = C(X * Y)$. On the other hand, from $X \subseteq C(X)$ and $Y \subseteq C(Y)$, it is immediately apparent that $X * Y \subseteq C(X) * C(Y)$ and, therefore, $C(X * Y) \subseteq C(C(X) * C(Y))$.
$(b) \Rightarrow (a)$: By (b) we have $C(C(X) * C(Y)) = C(X * Y)$ but $C(X) * C(Y) \subseteq C(C(X) * C(Y))$. Therefore, $C(X) * C(Y) \subseteq C(X * Y)$. □

## 12. Ideal-Completion of $\mathring{\mathbb{M}}^{s\ell}$-Algebras and $\mathbb{FL}$-Algebras.

In this section we present the construction called *ideal completion* for $\mathring{\mathbb{M}}^{s\ell}$-algebras. The main result is Theorem 23, where it is shown that every $\mathring{\mathbb{M}}^{s\ell}_{\sigma}$-algebra is embeddable in its ideal completion in such a way that all the existing residua and all the existing meets are preserved, although it can be shown that this embedding in general does not preserve arbitrary joins (see for instance (Bou et al., 2006, Proposition 4.14)). Let us stress that this embedding result was already known and was obtained by Ono (cf. (Ono, 2003a, Theorem 7)) using a method different from the one that we will use (see Note 1). Theorem 23 has as a consequence that every algebra in $\mathring{\mathbb{M}}^{\ell}_{\sigma}$, $\mathbb{PM}^{s\ell}_{\sigma}$, $\mathbb{PM}^{\ell}_{\sigma}$ and $\mathbb{FL}_{\sigma}$ is embeddable in its ideal completion.

**Definition 32** (Ideal of a $\vee$-semilattice)**.** *Given a $\vee$-semilattice $\mathbf{A} = \langle A, \vee\rangle$, an* ideal *of $\mathbf{A}$ is a subset $I \subseteq A$ such that:*

i) $I \neq \emptyset$,
ii) *if $y \leq x$ and $x \in I$, then $y \in I$,*
iii) *if $x, y \in I$, then $x \vee y \in I$.*

If $a \in A$, the set $\{x \in A : x \leq a\}$, which is denoted by $(a]$, is an ideal and is called the *principal ideal generated by a*. Given $X \subseteq A$, the smallest ideal containing $X$ is denoted by $(X]$ and it is called the *ideal generated by X*. Clearly, $(a] = (\{a\}]$. Note that $A$ is the greatest ideal and note also that when $A$ has a minimum element $\bot$, then $\bot$ belongs to every ideal and $\{\bot\}$ is the smallest ideal. An *ideal of a lattice* $\mathbf{A} = \langle A, \vee, \wedge\rangle$ is an ideal of the semilattice $\langle A, \vee\rangle$.

Let $\mathbf{A} = \langle A, \vee, *, 0, 1\rangle$ be an $\mathring{\mathbb{M}}^{s\ell}$-algebra and $C^{Id}$ be an operator on $\mathcal{P}(A)$ defined by $C^{Id}(X) = (X]$ for every $X \subseteq A$. It is easy to see that $C^{Id}$ is a closure operator. Clearly the sets that are $C^{Id}$-closed are exactly the ideals of the semilattice reduct of $\mathbf{A}$. In what follows, we will see that $C^{Id}$ satisfies the condition $(*)$ of Proposition 38.

**Lemma 16.** (Cf. Grätzer (1979)) *Let $\mathbf{A}$ be a $\vee$-semilattice and let $X \subseteq A$. Then the following holds:*

$$(X] = \{c \in A : c \leq a_1 \vee \cdots \vee a_n \text{ for some } a_i \in X\}$$

**Lemma 17.** *Let $\mathbf{A}$ be an $\mathring{\mathbb{M}}^{s\ell}$-algebra. For each $X, Y \subseteq A$, $(X] * (Y] \subseteq (X * Y]$.*

*Proof.* Let $c \in (X]$, $d \in (Y]$. By the previous lemma, we have that $c \leq a_1 \vee \cdots \vee a_n$ for some $a_i \in X$, and $d \leq b_1 \vee \cdots \vee b_m$ for some $b_j \in Y$. From this, applying monotonicity we have $c * d \leq (a_1 \vee \cdots \vee a_n) * (b_1 \vee \cdots \vee b_m)$ and so, applying the distributivity of $*$ with respect to $\vee$, we have that $c * d$ is smaller than or equal to the union of $n \times m$ elements of $X * Y$ and, therefore, $c * d \in (X * Y]$. □

**Proposition 41.** *Let $\mathbf{A}$ be an $\mathring{\mathbb{M}}^{s\ell}$-algebra. Then the structure $\mathbf{C}^{D}_{\langle A,*,1\rangle}$, built using the monoidal reduct of $\mathbf{A}$, the closure operator $C = C^{Id}$ and a $C^{Id}$-closed $D$, is a complete $\mathbb{FL}$-algebra.*

*Proof.* If $\mathbf{A}$ is an $\mathring{\mathbb{M}}^{s\ell}$-algebra, by Lemma 17, we have $C^{Id}(X) * C^{Id}(Y) \subseteq C^{Id}(X * Y)$. Thus, $C^{DM}$ satisfies condition $(*)$ of Proposition 38 and consequently, $\mathbf{C}^{D}_{\langle A,*,1\rangle}$ is a complete $\mathbb{FL}$-algebra. □

**Definition 33** (Ideal completion of an $\mathring{\mathbb{M}}^{s\ell}$-algebra)**.** *Let* $\mathbf{A}$ *be an* $\mathring{\mathbb{M}}^{s\ell}$*-algebra. The complete* $\mathbb{FL}$*-algebra* $\mathbf{C}^{(0]}_{\langle A,*,1\rangle}$ *built using the monoidal reduct of* $\mathbf{A}$*, the operator* $C^{Id}$ *and the* $C^{Id}$*-closed* $(0]$ *by the method of Proposition 38 will be called* ideal completion *of* $\mathbf{A}$ *and will be denoted by* $\mathbf{A^{Id}}$*. The set of all the* $C^{Id}$*-closed sets, that is, the set of the ideals of* $\mathbf{A}$*, will be denoted by* $A^{Id}$*.*

In the following result we show that if the initial $\mathring{\mathbb{M}}^{s\ell}$-algebra, $\mathbf{A}$, is an $\mathring{\mathbb{M}}^{s\ell}_{\sigma}$-algebra, then its ideal completion also satisfies the properties codified by $\sigma$; that is, $\mathbf{A^{Id}}$ is an $\mathbb{FL}_{\sigma}$-algebra.

**Proposition 42.** *When* $\mathbf{A}$ *is an* $\mathring{\mathbb{M}}^{s\ell}_{\sigma}$*-algebra, with* $\sigma \leq e w_l w_r c$*, then* $\mathbf{A^{Id}}$ *is a complete* $\mathbb{FL}_{\sigma}$*-algebra.*

*Proof.* If the monoidal operation of $\mathbf{A}$ is commutative, then clearly $\mathbf{A^{Id}}$ is an $\mathbb{FL}_e$-algebra. If $\mathbf{A}$ is a $\mathring{\mathbb{M}}^{s\ell}_{w_r}$-algebra, then 0 is the minimum of $A$. In this case, $(0\,]$ is the smallest ideal; that is, $(0] = (\emptyset]$ and therefore $\mathbf{A^{Id}}$ is an $\mathbb{FL}_{w_r}$-algebra. If $\mathbf{A}$ is an $\mathring{\mathbb{M}}^{s\ell}_{w_l}$-algebra, then 1 is the maximum of $A$. In this case, $(1]$ is the greatest ideal; that is, $(1] = A$ and, therefore, $\mathbf{A^{Id}}$ is an $\mathbb{FL}_{w_l}$-algebra. Thus, clearly, if $\mathbf{A}$ is an $\mathring{\mathbb{M}}^{s\ell}_{w}$-algebra, $\mathbf{A^{Id}}$ is an $\mathbb{FL}_w$-algebra. Finally, we suppose that $\mathbf{A}$ is an $\mathring{\mathbb{M}}^{s\ell}_{c}$-algebra. To see that $\mathbf{A^{Id}}$ is an $\mathbb{FL}_c$-algebra, we have to show that, for every ideal $I$ of $\mathbf{A}$, $I \subseteq (I * I]$. Let $a \in I$. Clearly, $a^2 \in I * I$ and, therefore, $a^2 \in (I * I]$ but, as $\mathbf{A}$ is an $\mathring{\mathbb{M}}^{s\ell}_{c}$-algebra, $a \leq a^2$ and thus, as $(I * I]$ is an ideal, we obtain $a \in (I * I]$. □

The following lemmas will be used to state the embedding theorem.

**Lemma 18.** *Let* $\mathbf{A}$ *be an* $\mathring{\mathbb{M}}^{s\ell}$*-algebra and let* $I, I_1, I_2 \in A^{Id}$*. The following conditions are satisfied:*[17]

$$
\begin{array}{lrcl}
 & \textit{1)} \quad I_1 \vee_C I_2 & = & \{a \in A : a \leq i_1 \vee i_2 \textit{ for some } i_1 \in I_1, i_2 \in I_2\} \\
 & \textit{2)} \quad I_1 \cap I_2 & = & \{a \in A : a \leq i_1 \wedge i_2 \textit{ for some } i_1 \in I_1, i_2 \in I_2\} \\
 & \textit{3)} \quad I_1 *_C I_2 & = & \{a \in A : a \leq i_1 * i_2 \textit{ for some } i_1 \in I_1, i_2 \in I_2\} \\
(3) & \textit{4)} \quad I_1 \backslash I_2 & = & \{a \in A : i_1 * a \in I_2 \textit{ for every } i_1 \in I_1\} \\
 & \textit{5)} \quad I_2 / I_1 & = & \{a \in A : a * i_1 \in I_2 \textit{ for every } i_1 \in I_1\} \\
 & \textit{6)} \quad I^{\backprime} & = & \{a \in A : a \leq i^{\backprime} \textit{ for every } i \in I\} \\
 & \textit{7)} \quad {}'I & = & \{a \in A : a \leq {}'i \textit{ for every } i \in I\}
\end{array}
$$

*Proof.* 1): $I_1 \vee_C I_2$ is the smallest ideal containing $I_1 \cup I_2$ and every such ideal must contain the set:

$$J = \{a \in A : a \leq i_1 \vee i_2 \text{ for some } i_1 \in I_1, i_2 \in I_2\}.$$

[17] To simplify the notation, we denote by $C$ the operator $C^{Id}$.

But the set $J$ is an ideal. So, if $b \leq a \in J$ and $a \leq i_1 \vee i_2$ then also $b \leq i_1 \vee i_2$ and, therefore, $b \in J$. If $a, b \in J$ and $a \leq i_1 \vee i_2$, $b \leq i'_1 \vee i'_2$, then $a \vee b \leq (i_1 \vee i_2) \vee (i'_1 \vee i'_2) = (i_1 \vee i'_1) \vee (i_2 \vee i'_2)$ and since $i_1 \vee i'_1 \in I_1$ and $i_2 \vee i'_2 \in I_2$, we have $a \vee b \in J$. As a consequence, $I_1 \vee_C I_2 = J$.
2): Let us consider the set:

$$K = \{a \in A : a \leq i_1 \wedge i_2 \text{ for some } i_1 \in I_1, i_2 \in I_2\}.$$

If $a \in K$ and $a \leq i_1 \wedge i_2$, then also $a \leq i_1$ i $a \leq i_2$ and hence $a \in I_1$ and $a \in I_2$; that is, $a \in I_1 \cap I_2$. Reciprocally, if $a \in I_1 \cap I_2$, then from $a \leq a \wedge a$ we can deduce $a \in K$.
3): $I_1 *_C I_2$ is the smallest ideal containing $I_1 * I_2$ and every such ideal must contain the set:

$$L = \{a \in A : a \leq i_1 * i_2 \text{ for some } i_1 \in I_1, i_2 \in I_2\}.$$

But $L$ is an ideal. So, if $b \leq a \in L$ and $a \leq i_1 * i_2$, then also $b \leq i_1 * i_2$ and thus, $b \in L$. If $a, b \in L$ and $a \leq i_1 * i_2$, $b \leq i'_1 * i'_2$, then $i_1 * i_2, i'_1 * i'_2 \leq (i_1 \vee i'_1) * (i_2 \vee i'_2)$. Therefore, $a \vee b \leq (i_1 \vee i'_1) * (i_2 \vee i'_2)$ and, since $i_1 \vee i'_1 \in I_1$ i $i_2 \vee i'_2 \in I_2$, we have that $a \vee b \in L$.
4), 5): These are exactly the definitions of the residuals in $\mathbf{A}^{\mathbf{Id}}$.
6), 7): By the definition we have $I^{\backprime} = \{a \in A : i * a \in (0] \text{ for each } i \in I\}$. But $i * a \in (0]$ is equivalent to $i * a \leq 0$ which by (LP) is equivalent to $a \leq i^{\backprime}$. The proof of (7) is analogous. □

**Lemma 19.** *Let* $\mathbf{A}$ *be an* $\mathring{\mathbb{M}}^{s\ell}$*-algebra.*

i) *For each* $a, b \in A$, *if* $a \backslash b$ *and* $b/a$ *exist, then* $(a \backslash b] = (a] \backslash (b]$ *and* $(a/b] = (a]/(b]$.
ii) *For each* $a \in A$, *if* $a^{\backprime}$ *and* $'a$ *exist, then* $(a^{\backprime}] = (a]^{\backprime}$ *and* $('a] = '(a]$.

*Proof. i*) Let $a, b \in A$ and suppose that the residual $a \backslash b$ exists. By (4) of Lemma 18 and (LR), we have:

$$(a] \backslash (b] = \{z \in A : x * z \leq b \text{ for every } x \leq a\} = \{z \in A : z \leq x \backslash b \text{ for every } x \leq a\}.$$

If $z \in (a] \backslash (b]$, then $z \leq x \backslash b$ for every $x \leq a$ and, in particular, $z \leq a \backslash b$; i.e., $z \in (a \backslash b]$. Now suppose that $z \leq a \backslash b$, which is equivalent to $a * z \leq b$. If $x \leq a$, by monotonicity we have that $x * z \leq a * z$ and thus $x * z \leq b$. Therefore, $z \in (a] \backslash (b]$. It is proved analogously that if $a/b$ exists, then $(a]/(b] = (a/b]$.
*ii*) This is a particular case of *i*) because if $a^{\backprime}$ and $'a$ exist, then $a^{\backprime} = a \backslash 0$ and $'a = 0/a$ and so $(a]^{\backprime} = (a] \backslash (0] = (a \backslash 0] = (a^{\backprime}]$ and, analogously, $'(a] = ('a]$. □

**Lemma 20.** *Let* $\mathbf{A}$ *be a* $\vee$*-semilattice and suppose* $X \subseteq A$. *Then, if there exists the infimum* $\bigwedge_{\mathbf{A}} X$ *of* $X$ *in* $\mathbf{A}$, *then the principal ideal generated by* $\bigwedge_{\mathbf{A}} X$ *is equal to the meet of the principal ideals generated by the elements of* $X$; *i.e.,* $(\bigwedge_{\mathbf{A}} X] = \bigcap_{a \in X} (a]$.

*Proof.* $b \in \bigcap_{a \in X} (a]$ *iff* $b \leq a$ for all $a \in X$ *iff* $b \leq \bigwedge_{\mathbf{A}} X$ *iff* $b \in (\bigwedge_{\mathbf{A}} X]$. □

**Theorem 23.** *For every* $\mathring{\mathbb{M}}^{s\ell}_{\sigma}$*-algebra* $\mathbf{A} = \langle A, \vee, *, 0, 1 \rangle$, *the mapping defined by* $i_{\mathbf{A}}(a) = (a]$, *for each* $a \in A$, *is an embedding (i.e., an* $\langle \vee, *, 0, 1 \rangle$*-monomorphism) from* $\mathbf{A}$ *into the* $\langle \vee, *, 0, 1 \rangle$*-reduct of the complete* $\mathbb{FL}_{\sigma}$*-algebra* $\mathbf{A}^{\mathbf{Id}}$, *which preserves all the residuals, pseudocomplements and existing meets.*

*Proof.* Clearly $i_{\mathbf{A}}$ is injective, because for each $a, b \in A$, if $a \neq b$ then $(a] \neq (b]$. Let $a, b \in A$ and $\odot \in \{\vee, *\}$. By the characterizations in Lemma 18 we have:

$$(a] \odot_C (b] = \{z \in A : z = x \odot y \text{ for some } x \leq a, y \leq b\}.$$

It is clear that $(a \odot b] \subseteq (a] \odot_C (b]$. If $z \in (a] \odot_C (b]$, then $z \leq x \odot y$ with $x \leq a$ and $y \leq b$ and, by the monotonicity of $\odot$, $z \leq a \odot b$; i.e., $z \in (a \odot b]$. Therefore, $(a] \odot_C (b] = (a \odot b]$. Moreover, the distinguished element and the unit element of $\mathbf{A^{Id}}$ are, respectively, the principal ideals $(0] = i_{\mathbf{A}}(0)$ and $(1] = i_{\mathbf{A}}(1)$. Therefore, $i_{\mathbf{A}}$ is an homomorphism. By Lemma 19 we have that $i_{\mathbf{A}}$ preserves the existing residuals and pseudocomplements. Finally, by Lemma 20, we have that $i_{\mathbf{A}}$ preserves the existing meets. □

The previous result was already given in the context of $\mathring{\mathbb{M}}^{s\ell}_{ew}$ in (Bou et al., 2006, Theorem 4.15). Theorem 23 allows us to obtain the following consequences.

**Corollary 27.** *For every algebra* $\mathbf{A}$ *in* $\mathbb{PM}^{s\ell}_{\sigma}$*, the mapping* $i_{\mathbf{A}}$ *is an embedding from* $\mathbf{A}$ *into the* $\langle \vee, *, ^{\backprime}, ^{\prime}, 0, 1\rangle$*-reduct of the complete* $\mathbb{FL}_{\sigma}$*-algebra* $\mathbf{A^{Id}}$*. This embedding preserves the existing residuals and meets.*

**Corollary 28.** *For every algebra* $\mathbf{A}$ *in* $\mathring{\mathbb{M}}^{\ell}_{\sigma}$*, the mapping* $i_{\mathbf{A}}$ *is an embedding from* $\mathbf{A}$ *into the* $\langle \vee, \wedge, *, 0, 1\rangle$*-reduct of the complete* $\mathbb{FL}_{\sigma}$*-algebra* $\mathbf{A^{Id}}$*. This embedding preserves the existing pseudocomplements, residuals and meets.*

**Corollary 29.** *For every algebra* $\mathbf{A}$ *in* $\mathbb{PM}^{\ell}_{\sigma}$*, the mapping* $i_{\mathbf{A}}$ *is an embedding from* $\mathbf{A}$ *into the* $\langle \vee, *, ^{\backprime}, ^{\prime}, 0, 1\rangle$*-reduct of the complete* $\mathbb{FL}_{\sigma}$*-algebra* $\mathbf{A^{Id}}$*. This embedding preserves the existing residuals and meets.*

**Corollary 30.** *For every algebra* $\mathbf{A}$ *in* $\mathbb{FL}_{\sigma}$*, the mapping* $i_{\mathbf{A}}$ *is an embedding from* $\mathbf{A}$ *into the* $\langle \vee, \wedge, *, \backslash, /, ^{\backprime}, ^{\prime}, 0, 1\rangle$*-reduct of the complete* $\mathbb{FL}_{\sigma}$*-algebra* $\mathbf{A^{Id}}$*. This embedding preserves the existing meets.*

**Definition 34.** *If* $\mathbf{A}$ *is an algebra in* $\mathring{\mathbb{M}}^{\ell}_{\sigma}$*,* $\mathbb{PM}^{s\ell}_{\sigma}$*,* $\mathbb{PM}^{\ell}_{\sigma}$ *or* $\mathbb{FL}_{\sigma}$*, the complete* $\mathbb{FL}_{\sigma}$*-algebra obtained by ideal completion of its* $\langle \vee, *, 0, 1\rangle$*-reduct it will be denoted by* $\mathbf{A^{Id}}$ *and we will call it the* ideal completion *of* $\mathbf{A}$*.*

**Remark 1.** *The results in Theorem 23 and Corollary 30 are obtained by Ono in (Ono, 2003a, Theorem 7) for* $\mathring{\mathbb{M}}^{s\ell}_{ew}$*-algebras and* $\mathbb{FL}_{ew}$*-algebras by using a different method. As Ono notes in (Ono, 2003a, p. 435), his method can easily be adapted to the non-*0*-bounded, non-integral or non-commutative cases.*

As an immediate consequence of the previous embedding, we have the following results:

**Theorem 24.** *Let* $\sigma \leq ew_l w_r c$*.*

i) $\mathring{\mathbb{M}}^{s\ell}_{\sigma}$ *is the class of the* $\langle \vee, *, 0, 1\rangle$*-subreducts of the algebras of the classes* $\mathbb{PM}^{s\ell}_{\sigma}$ *and* $\mathbb{FL}_{\sigma}$*.*
ii) $\mathring{\mathbb{M}}^{\ell}_{\sigma}$ *is the class of the* $\langle \vee, \wedge, *, 0, 1\rangle$*-subreducts of the algebras of the classes* $\mathbb{PM}^{\ell}_{\sigma}$ *and* $\mathbb{FL}_{\sigma}$*.*
iii) $\mathbb{PM}^{s\ell}_{\sigma}$ *is the class of the* $\langle \vee, *, ^{\backprime}, ^{\prime}, 0, 1\rangle$*-subreducts of the algebras of the class* $\mathbb{FL}_{\sigma}$*.*

iv) $\mathbb{PM}^{\ell}_{\sigma}$ *is the class of the* $\langle \vee, \wedge, *, ^{\backslash}, ', 0, 1\rangle$*-subreducts of the algebras of the class* $\mathbb{FL}_{\sigma}$.

*Proof.* *i*) is a consequence of Theorem 23, *ii*) of Corollary 28, *iii*) of Corollary 27 and *iv*) of Corollary 29. □

This last theorem was already proved in (Bou et al., 2006, Theorem 4.16) for the $\sigma = ew$ case.

# Part 5. ANALYSIS OF THE FRAGMENTS WITHOUT IMPLICATIONS OF $\mathcal{FL}_{\sigma}[\Psi]$ AND $\mathfrak{e}\mathcal{FL}_{\sigma}[\Psi]$.

In this part we study the fragments in the languages

$$\langle \vee, *, 0, 1\rangle, \langle \vee, \wedge, *, 0, 1\rangle, \langle \vee, \wedge, *, ^{\backslash}, ', 0, 1\rangle, \text{ and } \langle \vee, *, ^{\backslash}, ', 0, 1\rangle$$

of the systems $\mathcal{FL}_{\sigma}$ and their associated external systems $\mathfrak{e}\mathcal{FL}_{\sigma}$. Using the algebraization results of Part 3 and those in Part 4 where it was established that these classes are subreducts of $\mathbb{FL}_{\sigma}$, we obtain that the mentioned subsystems are fragments of $\mathcal{FL}_{\sigma}$ and that the corresponding external systems are fragments of $\mathfrak{e}\mathcal{FL}_{\sigma}$. It is also shown that each system, $\mathcal{FL}_{\sigma}$, is equivalent to its associated external system; but it is shown that the fragments considered are not equivalent to any Hilbert system. We also show that $\mathfrak{e}\mathcal{FL}_{\sigma}[\vee, *, 0, 1]$, $\mathfrak{e}\mathcal{FL}_{\sigma}[\vee, \wedge, *, 0, 1]$, $\mathfrak{e}\mathcal{FL}_{\sigma}[\vee, *, ^{\backslash}, ', 0, 1]$ and $\mathfrak{e}\mathcal{FL}_{\sigma}[\vee, \wedge, *, ^{\backslash}, ', 0, 1]$ are not protoalgebraic, but have respectively the varieties $\mathring{\mathbb{M}}^{s\ell}_{\sigma}$, $\mathring{\mathbb{M}}^{\ell}_{\sigma}$, $\mathbb{PM}^{s\ell}_{\sigma}$ and $\mathbb{PM}^{\ell}_{\sigma}$ as algebraic semantics with defining equation $1 \vee p \approx p$. In Section 16, we give some decidability results for some of the fragments considered.

## 13. Fragments.

First of all, we characterize the fragments of $\mathbb{FL}_{\sigma}$ in the languages considered.

**Theorem 25.** *Let* $\Psi$ *be one of the languages* $\langle \vee, *, 0, 1\rangle$, $\langle \vee, \wedge, *, 0, 1\rangle$, $\langle \vee, *, ^{\backslash}, ', 0, 1\rangle$ *or* $\langle \vee, \wedge, *, ^{\backslash}, ', 0, 1\rangle$. *For each* $\sigma$, *the system* $\mathcal{FL}_{\sigma}[\Psi]$ *is the* $\Psi$*-fragment of* $\mathcal{FL}_{\sigma}$.

*Proof.* We want to see that, for every $\Phi \cup \{\varsigma\} \subseteq Seq_{\Psi}^{\omega \times \{0,1\}}$,

$$\Phi \vdash_{\mathbf{FL}_{\sigma}} \varsigma \quad \text{iff} \quad \Phi \vdash_{\mathbf{FL}_{\sigma}[\Psi]} \varsigma.$$

Let $\tau$ be the translation of Definition 30. Then we have the following chain of equivalences:

$$\Phi \vdash_{\mathbf{FL}_{\sigma}} \varsigma \quad \text{iff} \quad \tau(\Phi) \vDash_{\mathbb{FL}_{\sigma}} \tau(\varsigma) \quad \text{iff} \quad \tau(\Phi) \vDash_{\mathbb{K}_{\sigma}[\Psi]} \tau(\varsigma) \quad \text{iff} \quad \Phi \vdash_{\mathbf{FL}_{\sigma}[\Psi]} \varsigma.$$

The first equivalence is obtained by applying Theorem 21; the second is due to the fact that, in each case, the class $\mathbb{K}_{\sigma}[\Psi]$ is the class of all the $\Psi$-subreducts of the class $\mathbb{FL}_{\sigma}$ (Theorem 24); and the third is a consequence of Theorem 21 □

Now we arrive to the main result of the paper; the result concerning four fragments without implications of the full Lambek logic **HFL** and of all its substructural extensions.

**Corollary 31.** *Let $\Psi$ be one of the languages:*

$$\langle \vee, *, 0, 1\rangle, \langle \vee, \wedge, *, 0, 1\rangle, \langle \vee, *, \backslash, /, 0, 1\rangle, \langle \vee, \wedge, *, \backslash, /, 0, 1\rangle.$$

*For each $\sigma$, the external system $\mathfrak{e}\mathcal{FL}_\sigma[\Psi]$ associated with $\mathcal{FL}_\sigma[\Psi]$ is the $\Psi$-fragment of the external system $\mathfrak{e}\mathcal{FL}_\sigma$ associated with $\mathcal{FL}_\sigma$.*

*Proof.* Applying the definition of external system and Theorem 25 we have, for every $\Gamma \cup \{\varphi\} \subseteq Fm_\Psi$:

$$\Gamma \vdash_{\mathfrak{e}\mathcal{FL}_\sigma} \varphi \quad \text{iff} \quad \{\emptyset \Rightarrow \gamma : \gamma \in \Gamma\} \vdash_{\mathbf{FL}_\sigma} \emptyset \Rightarrow \varphi \quad \text{iff}$$
$$\{\emptyset \Rightarrow \gamma : \gamma \in \Gamma\} \vdash_{\mathbf{FL}_\sigma[\Psi]} \emptyset \Rightarrow \varphi \quad \text{iff} \quad \Gamma \vdash_{\mathfrak{e}\mathcal{FL}_\sigma[\Psi]} \varphi. \qquad \square$$

## 14. Equivalence between $\mathcal{FL}_\sigma$ Systems and their Associated External Systems.

Next, using the algebraization results for the systems $\mathcal{FL}_\sigma$, we will show that every one of these Gentzen systems is equivalent to its associated external system.[18]

**Theorem 26.** *Let $\sigma \leq ew_lw_rc$. The system $\mathcal{FL}_\sigma$ is equivalent to its associated external system $\mathfrak{e}\mathcal{FL}_\sigma$.*

*Proof.* We define the translations $\tau'$ from $\mathfrak{L}$-sequents to $\mathfrak{L}$-formulas and $\rho'$ from $\mathfrak{L}$-formulas to $\mathfrak{L}$-sequents in the following way:

$$\tau'(\varphi_0, ..., \varphi_{m-1} \Rightarrow \varphi) := \begin{cases} \{\varphi_{m-1}\backslash(\varphi_{m-2}\backslash(\dots\backslash(\varphi_0\backslash\varphi)\dots))\}, & \text{if } m \geq 1, \\ \{\varphi\}, & \text{if } m = 0, \end{cases}$$

$$\tau'(\varphi_0, ..., \varphi_{m-1} \Rightarrow \emptyset) := \begin{cases} \{\varphi_{m-1}\backslash(\varphi_{m-2}\backslash(\dots\backslash(\varphi_0\backslash 0)\dots))\}, & \text{if } m \geq 1, \\ \{0\}, & \text{if } m = 0, \end{cases}$$

$\rho'(\varphi) := \{\emptyset \Rightarrow \varphi\}$.

To see that $\mathcal{FL}_\sigma$ and $\mathfrak{e}\mathcal{FL}_\sigma$ are equivalent, we will prove that the following conditions are satisfied (see Section 2.8):

*a*) For every $\Gamma \cup \{\varphi\} \subseteq Fm_{\mathfrak{L}}$, $\Gamma \vdash_{\mathfrak{e}\mathcal{FL}_\sigma} \varphi$ iff $\rho'[\Gamma] \vdash_{\mathbf{FL}_\sigma} \rho'(\varphi)$,

*b*) For every $\varsigma \in Seq_{\mathfrak{L}}^{\omega \times \{0,1\}}$, $\varsigma \dashv\vdash_{\mathbf{FL}_\sigma} \rho'\tau'(\varsigma)$.

*a*): Given that $\vdash_{\mathfrak{e}\mathcal{FL}_\sigma}$ and $\vdash_{\mathbf{FL}_\sigma}$ are finitary, we will restrict ourselves to finite sets of formulas without loss of generality. Suppose $\Gamma = \{\varphi_0, ..., \varphi_{m-1}\}$, with $m \in \omega$. Then we have:

$$\Gamma \vdash_{\mathfrak{e}\mathcal{FL}_\sigma} \varphi \quad \text{iff} \quad \{\emptyset \Rightarrow \varphi_0, ..., \emptyset \Rightarrow \varphi_{m-1}\} \vdash_{\mathfrak{e}\mathcal{FL}_\sigma} \emptyset \Rightarrow \varphi \quad \text{iff} \quad \rho'[\Gamma] \vdash_{\mathbf{FL}_\sigma} \rho'(\varphi).$$

*b*): Let $\varsigma = \Gamma \Rightarrow \Delta$. Let us define $\delta$ as the formula 0 if $\Delta$ is the empty sequence and as the formula $\varphi$ if $\Delta$ is constituted by the formula $\varphi$. If $\Gamma$ is a sequence of $m$ formulas $\varphi_0, \dots, \varphi_{m-1}$ we will use the following abbreviation:

$$\Gamma\backslash\delta := \begin{cases} \varphi_{m-1}\backslash(\varphi_{m-2}\backslash(\dots\backslash(\varphi_0\backslash\delta)\dots)), & \text{if } m \geq 1; \\ \delta, & \text{if } m = 0. \end{cases}$$

[18] In Section 4 above we summarize some Hilbert-style presentations for these external systems.

Using these conventions and the definition of the translations $\rho'$ and $\tau$ we have:

$$\rho'\tau'(\Gamma \Rightarrow \Delta) = \emptyset \Rightarrow \Gamma\backslash\delta.$$

From this, using the fact that, for every $\mathfrak{L}$-formula $\psi$, the sequents $\emptyset \Rightarrow \psi$ and $1 \Rightarrow \psi$ are interderivable, we obtain:

$$\rho'\tau'(\Gamma \Rightarrow \Delta) \dashv\vdash_{\mathbf{FL}_\sigma} 1 \Rightarrow \Gamma\backslash\delta.$$

Let $\rho$ be the translation considered in Definition 30. By Lemma 10 we have that the sequent $1 \Rightarrow \Gamma\backslash\delta$ and the sequents in $\rho(1 \preccurlyeq \Gamma\backslash\delta)$ are interderivable and, therefore:

$$\rho'\tau'(\Gamma \Rightarrow \Delta) \dashv\vdash_{\mathbf{FL}_\sigma} \rho(1 \preccurlyeq \Gamma\backslash\delta).$$

Now observe that in the equational system associated with $\mathbb{FL}_\sigma$, the equations $1 \preccurlyeq \Gamma\backslash\delta$ and $\prod\Gamma \preccurlyeq \delta$ are interderivable (by the Law of Residuation). Thus, since $\mathcal{FL}_\sigma$ and $\langle\mathfrak{L}, \vDash_{\mathbb{FL}_\sigma}\rangle$ are equivalent (Theorem 21), the translations by $\rho$ of these equations are interderivable in $\mathcal{FL}_\sigma$:

$$\rho(1 \preccurlyeq \Gamma\backslash\delta) \dashv\vdash_{\mathbf{FL}_\sigma} \rho(\prod\Gamma \preccurlyeq \delta).$$

But $\rho(\prod\Gamma \preccurlyeq \delta) = \rho\tau(\Gamma \Rightarrow \Delta)$, where $\tau$ is the translation from sequents to equations given in Definition 30 and, therefore, again applying the equivalence of $\mathcal{FL}_\sigma$ and $\langle\mathfrak{L}, \vDash_{\mathbb{FL}}\rangle$ we have: $\rho\tau(\Gamma \Rightarrow \Delta) \dashv\vdash_{\mathbf{FL}_\sigma} \Gamma \Rightarrow \Delta$; and consequently:

$$\rho'\tau'(\Gamma \Rightarrow \Delta) \dashv\vdash_{\mathbf{FL}_\sigma} \Gamma \Rightarrow \Delta.$$ □

Now we give an alternative proof of the well-known result concerning the algebraization of the external system associated with a calculus $\mathbf{FL}_\sigma$ (see for instance Galatos & Ono (2006) ) using Theorem 26 and the algebraization result for $\mathcal{FL}_\sigma$.

**Corollary 32.** *For every sequence, $\sigma$, the Hilbert system $\mathfrak{e}\mathcal{FL}_\sigma$ is algebraizable with the variety $\mathbb{FL}_\sigma$ as its equivalent algebraic semantics and with the set of equivalence formulas: $\varphi\Delta\psi = \{\varphi\backslash\psi, \psi\backslash\varphi\}$; and the defining equation: $1 \vee \varphi \approx \varphi$.*

*Proof.* Fix a sequence $\sigma$ and take the translations $\tau$ and $\rho$ of Theorem 21 and $\tau'$ and $\rho'$ of Theorem 26 We define the translations $\tau'' := \tau\rho'$ from $\mathfrak{L}$-formulas to $\mathfrak{L}$-equations and $\rho'' := \tau'\rho$ from $\mathfrak{L}$-equations to $\mathfrak{L}$-formulas; that is, $\tau''(\varphi) = \{1 \preccurlyeq \varphi\}$ and $\rho''(\varphi \approx \psi) = \{\varphi\backslash\psi, \psi\backslash\varphi\}$. With these translations it is immediately clear that $\mathfrak{e}\mathcal{FL}_\sigma$ is algebraizable with the variety $\mathbb{FL}_\sigma$ as its equivalent algebraic semantics. Moreover, given the definitions of $\tau''$ and $\rho''$, we have that the set of equivalence formulas is: $\{\varphi\backslash\psi, \psi\backslash\varphi\}$; and the defining equation is: $1 \preccurlyeq \varphi$. □

## 15. Non-Equivalence of the Implication-less Fragments Systems and their Associated External Systems.

In this section we show that the classes $\mathring{\mathbb{M}}^{s\ell}_\sigma$, $\mathring{\mathbb{M}}^{\ell}_\sigma$, $\mathbb{PM}^{s\ell}_\sigma$ and $\mathbb{PM}^{\ell}_\sigma$ are not equivalent to any Hilbert system. As a corollary, we obtain that the systems $\mathcal{FL}_\sigma[\Psi]$, where $\Psi$ is one of the four languages considered, are not equivalent to any Hilbert system and therefore, unlike the case of the systems $\mathcal{FL}_\sigma$ and $\mathfrak{e}\mathcal{FL}_\sigma$, they are not equivalent to their associated external systems.

**Theorem 27.** *The varieties $\mathring{\mathbb{M}}^{s\ell}_{\sigma}$, $\mathring{\mathbb{M}}^{\ell}_{\sigma}$ are not the equivalent algebraic semantics to any Hilbert system.*

*Proof.* Let us consider the four-element distributive lattice

$$\mathbf{A} = \langle \{0, a, b, 1\}, \vee, \wedge\rangle,$$

with $0 < a < b < 1$. In (Font & Verdú, 1991, Proposition 2.1) it is proved that the Leibniz operator $\Omega_{\mathbf{A}}$ cannot be an isomorphism between the filters of an arbitrary Hilbert system and the congruences of the algebra $\mathbf{A}$. Now let $\mathbf{A}' = \langle \{0, a, b, 1\}, \vee, \wedge, *, 0, 1\rangle$, where $* = \wedge$; then $\mathbf{A}' \in \mathring{\mathbb{M}}^{\ell}_{\sigma}$. It is easy to check that the proof of Proposition 2.1 inFont & Verdú (1991) applies in our case to show that $\mathring{\mathbb{M}}^{\ell}_{\sigma}$ is not the equivalent algebraic semantics of any Hilbert system. The same proof also works for $\mathring{\mathbb{M}}^{s\ell}_{\sigma}$. □

**Theorem 28.** *The varieties $\mathbb{PM}^{s\ell}_{\sigma}$, $\mathbb{PM}^{\ell}_{\sigma}$ are not the equivalent algebraic semantics of any Hilbert system.*

*Proof.* Let us consider the five-element pseudocomplemented lattice

$$\mathbf{A} = \langle \{0, a, b, c, 1\}, \vee, \wedge, \neg\rangle,$$

with $0 < a < b < c < 1$, and $\neg 0 = 1$, $\neg a = \neg b = \neg c = \neg 1 = 0$. In (Rebagliato & Verdú, 1993, Theorem 3.1) it is proved that the Leibniz operator $\Omega_{\mathbf{A}}$ cannot be an isomorphism between the filters of an arbitrary Hilbert system and the congruences of the algebra $\mathbf{A}$. Now let $\mathbf{A}' = \langle \{0, a, b, c, 1\}, \vee, \wedge, *, \neg, 0, 1\rangle$, where $* = \wedge$; then $\mathbf{A}' \in \mathbb{PM}^{\ell}_{\sigma}$. It is easy to check that the proof of Theorem 3.1 in Rebagliato & Verdú (1993) applies in our case to show that $\mathbb{PM}^{\ell}_{\sigma}$ is not the equivalent algebraic semantics of any Hilbert system. The same proof also works for $\mathbb{PM}^{s\ell}_{\sigma}$. □

**Corollary 33.** *The Gentzen systems*

$$\mathcal{FL}_{\sigma}[\vee, *, 0, 1],\ \mathcal{FL}_{\sigma}[\vee, \wedge, *, 0, 1],\ \mathcal{FL}_{\sigma}[\vee, *, {}^{\backprime}, {}^{\prime}, 0, 1] \text{ and } \mathcal{FL}_{\sigma}[\vee, \wedge, *, {}^{\backprime}, {}^{\prime}, 0, 1]$$

*are not equivalent to any Hilbert system and, therefore, none of these Gentzen systems is equivalent to its associated external system.*

*Proof.* The Gentzen systems

$$\mathcal{FL}_{\sigma}[\vee, *, 0, 1],\ \mathcal{FL}_{\sigma}[\vee, \wedge, *, 0, 1],\ \mathcal{FL}_{\sigma}[\vee, *, {}^{\backprime}, {}^{\prime}, 0, 1], \text{ and } \mathcal{FL}_{\sigma}[\vee, \wedge, *, {}^{\backprime}, {}^{\prime}, 0, 1]$$

are algebraizable with the varieties $\mathring{\mathbb{M}}^{s\ell}_{\sigma}$, $\mathring{\mathbb{M}}^{\ell}_{\sigma}$, $\mathbb{PM}^{s\ell}_{\sigma}$ and $\mathbb{PM}^{\ell}_{\sigma}$, respectively, as their equivalent algebraic semantics (Theorem 21). If these Gentzen systems were equivalent to their respective external systems, then these varieties would be the equivalent algebraic semantics of the corresponding external systems, which contradicts Theorem 28. □

**Theorem 29.** *The Hilbert systems*

$$\mathfrak{e}\mathcal{FL}_{\sigma}[\vee, *, 0, 1],\ \mathfrak{e}\mathcal{FL}_{\sigma}[\vee, \wedge, *, 0, 1],\ \mathfrak{e}\mathcal{FL}_{\sigma}[\vee, *, {}^{\backprime}, {}^{\prime}, 0, 1] \text{ and } \mathfrak{e}\mathcal{FL}_{\sigma}[\vee, \wedge, *, {}^{\backprime}, {}^{\prime}, 0, 1]$$

*are not protoalgebraic.*

*Proof.* Let $\Psi$ be one of the four languages considered. Each system $\mathfrak{e}\mathcal{FL}_\sigma[\Psi]$ is a subsystem of $\mathfrak{e}\mathcal{FL}_{ewc}[\Psi]$ which, by Corollary 31, is the $\Psi$-fragment of $\mathfrak{e}\mathcal{FL}_{ewc}$; that is, of the intuitionistic logic presented in the language $\mathfrak{L}$ of **FL**. As the protoalgebraicity is monotonic, if we have that the systems $\mathfrak{e}\mathcal{FL}_{ewc}[\Psi]$ are not protoalgebraic, then we will be able to conclude that none of the systems $\mathfrak{e}\mathcal{FL}_\sigma[\Psi]$ are protoalgebraic.

The system $\mathfrak{e}\mathcal{FL}_{ewc}[\vee, *, 0, 1]$ is the $\langle \vee, *, 0, 1\rangle$-fragment of $\mathfrak{e}\mathcal{FL}_{ewc}$; that is, of the intuitionistic logic that (given the presence of the rules of left weakening and contraction in the corresponding calculus) is a notational copy of the $\langle \vee, \wedge, 0, 1\rangle$-fragment of the classical logic; and it is known that this fragment is not protoalgebraic (see Font & Verdú (1991)). Given that the behavior of the connectives $\wedge$ and $*$ is the same in this context, it is easy to see that $\mathfrak{e}\mathcal{FL}_{ewc}[\vee, \wedge, *, 0, 1]$ is not protoalgebraic either.

The system $\mathfrak{e}\mathcal{FL}_{ewc}[\vee, *, ^{\backprime}, ', 0, 1]$ is definitionally equivalent to $\mathfrak{e}\mathcal{FL}_{ewc}[\vee, *, \neg, 0, 1]$ and the system is a notational copy of the $\langle \vee, \wedge, \neg, 0, 1\rangle$-fragment of the intuitionistic logic, which is known not to be protoalgebraic (see (Blok & Pigozzi, 1989, Theorem 5.13)). It is easy to see that $\mathfrak{e}\mathcal{FL}_{ewc}[\vee, \wedge, *, \neg, 0, 1]$ is not protoalgebraic either. □

Although they are non-protoalgebraic, these systems have an algebraic semantics, as we now show.

**Theorem 30.** *For every $\sigma$, $\mathring{\mathbb{M}}^{s\ell}_\sigma$ ($\mathring{\mathbb{M}}^{\ell}_\sigma$, $\mathbb{PM}^{s\ell}_\sigma$, $\mathbb{PM}^{\ell}_\sigma$) is an algebraic semantics for*

$$\mathfrak{e}\mathcal{FL}_\sigma[\vee, *, 0, 1](\mathfrak{e}\mathcal{FL}_\sigma[\vee, \wedge, *, 0, 1], \mathfrak{e}\mathcal{FL}_\sigma[\vee, *, ^{\backprime}, ', 0, 1], \mathfrak{e}\mathcal{FL}_\sigma[\vee, \wedge, *, ^{\backprime}, ', 0, 1])$$

*with the defining equation: $1 \preccurlyeq p$.*

*Proof.* As a consequence of Theorem 21, we have that for every $\Sigma \cup \{\varphi\} \subseteq Fm_{\langle \vee, *, 0, 1\rangle}$,

$$\{\emptyset \Rightarrow \psi : \psi \in \Sigma\} \vdash_{\mathbf{FL}_\sigma[\vee, *, 0, 1]} \emptyset \Rightarrow \varphi \quad \textit{iff} \quad \{\tau(\emptyset \Rightarrow \psi) : \psi \in \Sigma\} \models_{\mathring{\mathbb{M}}^{s\ell}_\sigma} \tau(\emptyset \Rightarrow \varphi).$$

That is,

$$\Sigma \vdash_{\mathfrak{e}\mathcal{FL}_\sigma[\vee, *, 0, 1]} \varphi \quad \textit{iff} \quad \{1 \preccurlyeq \psi : \psi \in \Sigma\} \models_{\mathring{\mathbb{M}}^{s\ell}_\sigma} 1 \preccurlyeq \varphi.$$

The other three cases are proved analogously. □

These results were obtained for the case $\sigma = ew$ in Adillon & Verdú (2002); Bou et al. (2006).

## 16. Some Results concerning Decidability.

In this section, we give some decidability results for some of the fragments considered. When $\sigma \leq w_l$, it is well known that the variety $\mathbb{FL}_\sigma$ has the finite embeddability property (FEP). Using this fact and the results regarding subreducts in Part 4, we show that, when $\sigma \leq w_l$, the varieties $\mathring{\mathbb{M}}^{s\ell}_\sigma$, $\mathring{\mathbb{M}}^{\ell}_\sigma$, $\mathbb{PM}^{s\ell}_\sigma$ and $\mathbb{PM}^{\ell}_\sigma$ have the FEP. This fact allow us to prove decidability results for the corresponding systems, $\mathcal{FL}_\sigma[\Psi]$ and $\mathfrak{e}\mathcal{FL}_\sigma[\Psi]$.

**Theorem 31.** (Cf.(Galatos et al., 2007b, Theorem 6.46)) *Let $\sigma$ be such that $w_l \leq \sigma$. The variety $\mathbb{FL}_\sigma$ has the FEP.*

**Theorem 32.** *Let $\sigma$ be such that $w_l \leq \sigma$. The varieties $\mathring{\mathbb{M}}^{s\ell}_\sigma$, $\mathring{\mathbb{M}}^{\ell}_\sigma$, $\mathbb{PM}^{s\ell}_\sigma$ and $\mathbb{PM}^{\ell}_\sigma$ have the FEP.*

*Proof.* It is enough to prove the first part. Let $\mathbf{A}$ be any algebra in $\mathring{\mathbb{M}}^{s\ell}_{\sigma}$ and let $\mathbf{B}$ be a finite partial subalgebra of $\mathbf{A}$. By Theorem 24, $\mathbf{A}$ is embeddable in an $\mathbb{FL}_{\sigma}$-algebra, $\mathbf{A}'$. Let $i$ be such an embedding. Now we have that $i[\mathbf{B}]$ is a finite partial subalgebra of $\mathbf{A}'$. By Theorem 31, we have that $\mathbb{FL}_{\sigma}$ has the FEP. Therefore, $\mathbf{A}'$ can be embedded in a finite $\mathbb{FL}_{\sigma}$-algebra, $\mathbf{D}$. Let $h$ be this embedding and let $\mathbf{D}'$ be the $\langle \vee, *, 0, 1\rangle$-reduct of $\mathbf{D}$. $\mathbf{D}'$ is a finite $\mathring{\mathbb{M}}^{s\ell}_{\sigma}$-algebra and the map $h \circ i$ is an embedding from $\mathbf{B}$ into $\mathbf{D}'$. A similar argument runs for $\mathring{\mathbb{M}}^{\ell}_{\sigma}$, $\mathbb{PM}^{s\ell}_{\sigma}$, and $\mathbb{PM}^{\ell}_{\sigma}$. □

**Corollary 34.** *Let $\sigma$ be as in Theorem 32 The quasi-equational (and universal) theory of each one of the varieties $\mathring{\mathbb{M}}^{s\ell}_{\sigma}$, $\mathring{\mathbb{M}}^{\ell}_{\sigma}$, $\mathbb{PM}^{s\ell}_{\sigma}$ and $\mathbb{PM}^{\ell}_{\sigma}$ is decidable.*

**Remark 19.** *Observe that the method used to obtain the previous result can only be applied to the case of the varieties $\mathring{\mathbb{M}}^{s\ell}_{\sigma}$, $\mathring{\mathbb{M}}^{\ell}_{\sigma}$, $\mathbb{PM}^{s\ell}_{\sigma}$ and $\mathbb{PM}^{\ell}_{\sigma}$ falling under the scope of Theorem 32 For example, it is well known that the varieties $\mathbb{FL}$, $\mathbb{FL}_e$ and $\mathbb{FL}_{w_r}$ do not have the FEP (see (Galatos et al., 2007b, Theorem 6.56)).*

**Corollary 35.** *Let $\sigma$be as in Theorem 32 The Gentzen systems*

$$\mathcal{FL}_{\sigma}[\vee, *, 0, 1], \mathcal{FL}_{\sigma}[\vee, \wedge, *, 0, 1], \mathcal{FL}_{\sigma}[\vee, *, {}^{\backprime}, {}', 0, 1], \text{ and } \mathcal{FL}_{\sigma}[\vee, \wedge, *, {}^{\backprime}, {}', 0, 1]$$

*are decidable, i.e., their sets of entailments of the form $\{\Gamma_i \Rightarrow \Delta_i : i \in I\} \vdash \Gamma \Rightarrow \Delta$, with $I$ finite, are decidable.*

*Proof.* It is an immediate consequence of the algebraization of the Gentzen systems (Theorem 21) and Corollary 32. □

**Corollary 36.** *Let $\sigma$ be as in Theorem 32 The external systems*

$$\mathfrak{e}\mathcal{FL}_{\sigma}[\vee, *, 0, 1], \mathfrak{e}\mathcal{FL}_{\sigma}[\vee, \wedge, *, 0, 1], \mathfrak{e}\mathcal{FL}_{\sigma}[\vee, *, {}^{\backprime}, {}', 0, 1], \text{ and } \mathfrak{e}\mathcal{FL}_{\sigma}[\vee, \wedge, *, {}^{\backprime}, {}', 0, 1]$$

*are decidable, i.e., their entailments of the form $\Gamma \vdash \varphi$, with $\Gamma$ finite, are decidable.*

*Proof.* By Theorem 30, these Hilbert systems have, respectively, the varieties $\mathring{\mathbb{M}}^{s\ell}_{\sigma}$, $\mathring{\mathbb{M}}^{\ell}_{\sigma}$, $\mathbb{PM}^{s\ell}_{\sigma}$ and $\mathbb{PM}^{\ell}_{\sigma}$ as their algebraic semantics. The result is an immediate consequence of this last fact and Corollary 32. □

## 17. Future Work.

The case of omplication-less in languages without additive connectives (i.e., disjunction and conjunction) but with fusion will be studied in a forthcoming paper. The work will be undertaken within the framework of the notion of order-algebraizability Raftery (2006). We also plan to extend the work in the current paper to the fragments without implications of the non-associative full Lambek calculus.

## 18. Funding.

The authors acknowledge support from the Spanish MICINN via projects MTM 2011-25747, TIN2012-39348-C02-01 and TASSAT TIN2010-20967-C04-01, as well as support from the Generalitat de Catalunya via grant 2009-SGR-1433.

## 19. Acknowledgments.

The authors thank Félix Bou for his helpful comments and suggestions.

## References

Adillon, R. (2001). *Contribució a l'estudi de les lògiques proposicionals intuïcionistes i de Gödel sense contracció, i dels seus fragments.* Ph. D. Dissertation, University of Barcelona.

Adillon, R., García-Cerdaña, A., & Verdú, V. (2007). On three implication-less fragments of t-norm based fuzzy logics. *Fuzzy Sets and Systems* **158**, 2575–2590.

Adillon, R., García-Cerdaña, À., & Verdú, V. (2007). On three implication-less fragments of t-norm based fuzzy logics. *Fuzzy Sets and Systems* **158**(23), 2575–2590.

Adillon, R., & Verdú, V. (2000). On a contraction-less intuitionistic propositional logic with conjunction and fusion. *Studia Logica, Special Issue on Abstract Algebraic Logic* **65**(1), 11–30.

Adillon, R., & Verdú, V. (2002). On a substructural Gentzen system, its equivalent variety semantics and its external deductive system. *Bulletin of the Section of Logic* **31**(3), 125–134.

Avron, A. (1988). The semantics and proof theory of linear logic. *Theoretical Computer Science* **570**, 161–184.

Balbes, R., & Dwinger, P. (1974). *Distributive lattices.* Columbia, Mo.: University of Missouri Press.

Birkhoff, G. (1973). *Lattice Theory* (3rd. ed.), Volume XXV of *Colloquium Publications.* Providence: American Mathematical Society. (1st. ed. 1940).

Blok, W. J., & Pigozzi, D. (1989, January). *Algebraizable logics*, Volume 396 of *Mem. Amer. Math. Soc.* Providence: A.M.S.

Bou, F., García-Cerdaña, À., & Verdú, V. (2006). On two fragments with negation and without implication of the logic of residuated lattices. *Archive for Mathematical Logic* **45**(5), 615–647.

Cintula, P., Hájek, P., & Horčík, R. (2007). Formal systems of fuzzy logic and their fragments. *Annals of Pure and Applied Logic* **150**(1–3), 40–65.

Davey, B. A., & Priestley, H. A. (1990). *Introduction to lattices and order.* Cambridge: Cambridge University Press.

Dubreil-Jacotin, L., Lesieur, L., & Croisot, R. (1953). *Leçons sur la théorie des treillis, des estructures algébriques ordonées et des treillis géométriques.* Gauthier-Villars, Paris.

Dyrda, K., & Prucnal, T. (1980). On finitely based consequence determined by a distributive lattice. *Bulletin of the Section of Logic* **9**, 60–66.

Esteva, F., Godo, L., & García-Cerdaña, À. (2003). On the hierarchy of t-norm based residuated fuzzy logics. In Fitting, M. & Orłowska, E., editors, *Beyond two: theory and applications of multiple-valued logic*, Volume 114 of *Studies in Fuzziness and Soft Computing*, pp. 251–272. Heidelberg: Physica.

Font, J. M., Guzmán, F., & Verdú, V. (1991). Characterization of the reduced matrices for the $\{\wedge, \vee\}$-fragment of classical logic. *Bulletin of the Section of Logic* **20**, 124–128.

Font, J. M., & Verdú, V. (1991). Algebraic logic for classical conjunction and disjunction. *Studia Logica, Special Issue on Algebraic Logic* **50**, 391–419.

Galatos, N., Jipsen, P., Kowalski, T., & Ono, H. (2007a). *Residuated Lattices: An Algebraic Glimpse at Substructural Logics*, Volume 151 of *Studies in Logic and the Foundations of Mathematics.* Elsevier.

Galatos, N., Jipsen, P., Kowalski, T., & Ono, H. (2007b). *Residuated Lattices: an algebraic glimpse at substructural logics*, Volume 151 of *Studies in Logic and the Foundations of Mathematics.* Amsterdam: Elsevier.

Galatos, N., & Ono, H. (2006). Algebraization, parametrized local deduction theorem and interpolation for substructural logics over FL. *Studia Logica* **83**, 1–32.

Galatos, N., & Ono, H. (2010). Cut elimination and strong separation for substructural logics: An algebraic approach. *Annals of Pure and Applied Logic* **161**(9), 1097–1133.

Gil, A. J., Torrens, A., & Verdú, V. (1997). On Gentzen systems associated with the finite linear MV-algebras. *Journal of Logic and Computation* **7**(4), 473–500.

Gottwald, S. (2001). *A treatise on many-valued logics*, Volume 9 of *Studies in Logic and Computation*. Baldock: Research Studies Press.

Gottwald, S., García-Cerdaña, À., & Bou, F. (2003). Axiomatizing monoidal logic. A correction to: "A treatise on many-valued logics". *Journal of Multiple-Valued Logic and Soft Computing* **9**(4), 427–433.

Grätzer, G. (1979). *Universal algebra* (Second Edition ed.). Berlin: Springer-Verlag.

Higgs, D. (1984). Dually residuated commutative monoids with identity element as least element do not form an equational class. *Mathematica Japonica* **29**(1), 69–75.

Höhle, U. (1995). Commutative, residuated *l*-monoids. In Höhle, U. & Klement, E. P., editors, *Non-classical logics and their applications to fuzzy subsets (Linz, 1992)*, Volume 32 of *Theory Decis. Lib. Ser. B Math. Statist. Methods*, pp. 53–106. Dordrecht: Kluwer Acad. Publ.

Hosoi, T. (1966a). Algebraic proof of the separation theorem on classical propositional calculus. *Proceedings of the Japan Academy* **42**, 67–69.

Hosoi, T. (1966b). On the separation theorem of intermediate propositional calculi. *Proceedings of the Japan Academy* **42**, 535–538.

Jipsen, P., & Tsinakis, C. (2002). A survey of residuated lattices. In Martinez, J., editor, *Ordered Algebraic Structures*, pp. 19–56. Dordrecht: Kluwer Academic Publishers.

Lakser, H. (1973). Principal congruences of pseudocomplemented distributive lattices. *Proceedings of the Americal Mathematical Society* **37**, 32–37.

Ono, H. (1990). Structural rules and a logical hierarchy. In Petkov, P. P., editor, *Mathematical logic, Proceedings of the Heyting'88 Summer School*, pp. 95–104. New York: Plenum.

Ono, H. (1993). Semantics for substructural logics. In Došen, K. & Schroeder-Heister, P., editors, *Substructural logics*, pp. 259–291. Oxford University Press.

Ono, H. (1998). Proof-theoretic methods for nonclassical logic - an introduction. In Takahashi, M., Okada, M., & Dezani-Ciancaglini, M., editors, *Theories of Types and Proofs*, MSJ Memoirs 2, pp. 207–254. Mathematical Society of Japan.

Ono, H. (2003a). Closure operators and complete embeddings of residuated lattices. *Studia Logica* **74**(3), 427–440.

Ono, H. (2003b). Substructural logics and residuated lattices - an introduction. In Hendricks, V. F. & Malinowski, J., editors, *50 Years of Studia Logica*, Volume 21 of *Trends in Logic—Studia Logica Library*, pp. 193–228. Dordrecht: Springer.

Ono, H., & Komori, Y. (1985). Logics without the contraction rule. *The Journal of Symbolic Logic* **50**(1), 169–201.

Porębska, M., & Wroński, A. (1975). A characterization of fragments of the intuitionistic propositional logic. *Reports on Mathematical Logic* **4**, 39–42.

Raftery, J. (2006). Correspondences between gentzen and hilbert systems. *The Journal of Symbolic Logic* **71**(3), 903–957.

Rebagliato, J., & Verdú, V. (1993). On the algebraization of some Gentzen systems. *Fundamenta Informaticae, Special Issue on Algebraic Logic and its Applications* **18**, 319–338.

Rebagliato, J., & Verdú, V. (1994). A finite Hilbert-style axiomatization of the implication-less fragment of the intuitionistic propositional calculus. *Mathematical Logic Quarterly* **40**, 61–68.

Rebagliato, J., & Verdú, V. (1995, June). Algebraizable Gentzen systems and the deduction theorem for Gentzen systems. Mathematics Preprint Series 175, University of Barcelona. `http://www.ub.edu/plie/personal_PLiE/verdu_HTML(2)/docs/ReV95-p.pdf`.

Troelstra, A. S. (1992). *Lectures on Linear Logic*. Number 29 in CSLI Lecture Notes. Stanford: CSLI Publications. Distributed by University of Chicago Press.

van Alten, C. J., & Raftery, J. G. (2004). Rule separation and embedding theorems for logics without weakening. *Studia Logica* **76**(2), 241–274.

[1]Institut d'investigació en Intel·ligència Artificial (IIIA-CSIC). Consell Superior d'Investigacions Científiques. Campus de la UAB, E-08193 Bellaterra, Catalonia (Spain).

[2] Departament de Tecnologies de la Informació i les Comunicacions. Universitat Pompeu Fabra. Tànger 122-140, E-08018 Barcelona, Catalonia (Spain).

[3]Departament de Probabilitat, Lògica i Estadística. Universitat de Barcelona. Gran Vía de Les Corts Catalanes 585, E-08007 Barcelona, Catalonia (Spain).,